\documentclass{amsart}

\usepackage{amsmath, amsthm}
\usepackage{amssymb, latexsym}
\usepackage{amsfonts}

\newtheorem{Theorem}[equation]{Theorem}
\newtheorem{theorem}[equation]{Theorem}
\newtheorem{Proposition}[equation]{Proposition}
\newtheorem{proposition}[equation]{Proposition}

\newtheorem{lemma}[equation]{Lemma}
\newtheorem{Lemma}[equation]{Lemma}
\newtheorem{Corollary}[equation]{Corollary}
\newtheorem{corollary}[equation]{Corollary}

\newtheorem{Remark}[equation]{Remark}
\newtheorem{exam}[equation]{Example}

\theoremstyle{definition}
\newtheorem{definition}[equation]{Definition}

\newcommand{\eqn}[1]{(\ref{#1})}
\newcommand{\beql}[1]{\begin{equation}\label{#1}}
\newcommand{\eeq} {\end{equation}}

    \font\Aaa=msam10

\newcommand\qedd{\hbox{~~\Aaa\char'003}}

\font\Bbb=msbm10
\newcommand\semi{\hbox{\Bbb o}}

\newcommand\Z{\hbox{\Bbb Z}}
\newcommand\Q{\hbox{\Bbb Q}}
\newcommand\F{\hbox{\Bbb F}}

\font\Bbcc=msbm7
\newcommand\Fb{\hbox{\Bbcc F}}

\numberwithin{equation}{section}

\newcommand\Proof {\noindent{\em Proof.} }
\newcommand\M{{\mathcal M}}

\DeclareMathOperator{\Sym}{Sym}

\newcommand\e{{ \epsilon }}
\newcommand\g{\gamma }

\newcommand\z{{ \zeta }}

        \DeclareMathOperator{\PSL}{PSL}
        \DeclareMathOperator{\SL}{SL}
        \DeclareMathOperator{\St}{St}
        
        \DeclareMathOperator{\SU}{SU}
        
        \DeclareMathOperator{\PSU}{PSU}
        
        \DeclareMathOperator{\PGL}{PGL}
        \newcommand\BCRW{{\rm BCRW}}
        
        \DeclareMathOperator{\PSp}{PSp}
        \DeclareMathOperator{\PO}{P\Omega}
        \DeclareMathOperator{\Spp}{Sp}
        
        \DeclareMathOperator{\GL}{GL}

        \DeclareMathOperator{\AGL}{AGL}

        \DeclareMathOperator{\diag}{diag}
        \newcommand{\Sz}{{\rm Sz}}

\DeclareRobustCommand{\SkipTocEntry}[4]{}
\begin{document}

\title[Presentations of finite simple groups]
{Presentations of finite simple groups:  \\a quantitative approach}

\thanks{The authors were partially supported by
        NSF grants DMS 0140578,
          DMS 0242983, DMS 0600244  and   DMS~0354731.    The  authors
are grateful for the support and hospitality of   the Institute for
Advanced Study, where this research was carried out.
The research by the fourth author also was supported by the Ambrose Monell
   Foundation and the Ellentuck Fund.}

       \author{R. M. Guralnick}
       \address{Department of Mathematics, University of Southern California,
       Los Angeles, CA 90089-2532, USA}
       \email{guralnic@usc.edu}

       \author{W. M. Kantor}
       \address{Department of Mathematics, University of Oregon,
       Eugene, OR 97403, USA}
       \email{kantor@math.uoregon.edu}

       \author{M. Kassabov}
       \address{Department of Mathematics, Cornell University, Ithaca, NY 14853-4201, USA}
       \email{kassabov@math.cornell.edu}

       \author{A. Lubotzky}
       \address{Department of Mathematics, Hebrew University, Givat Ram, Jerusalem 91904, Israel}
       \email{alexlub@math.huji.ac.il}


\begin{abstract}
There is a constant  $C_0$ 
such that  all nonabelian finite
simple groups of rank $n$ over $\Fb_q$,
with the possible exception of the Ree groups $^2G_2(3^{2e+1})$,  have 
presentations with  at most $C_0$  generators and relations and
total length at most 
$C_0(
\log n +\log q)$. As a corollary, we deduce a conjecture of Holt: there is
a constant
$C$ such that   $\dim   H^2(G,M)
\leq C\dim   M$  for every finite simple group $G$, every
 prime $p$ and
every irreducible
 $\Fb_p  G $-module $M$.
   
\end{abstract}

\maketitle

\tableofcontents

\section {Introduction}
\label {Introduction} 
In this paper we study presentations of finite simple groups $G$ from 
a quantitative point of view.
%
Our main result  provides  unexpected answers to the
following questions: how many relations are needed to define $G$, and
how short can these relations be?

The classification of the finite simple groups states that every
nonabelian finite simple group is alternating, Lie type, or one of 26
sporadic groups.  The latter are of no relevance to our asymptotic
results.
Instead, we will primarily deal with groups of Lie type, which  have a 
(relative) {\em rank}
$n$ over a field
$\F_q$.  In order to  keep our results uniform, we  
view the alternating group $A_n$  
and symmetric group $S_n$ as
     groups of rank $n-1$ over ``the field
$\F_1$  with 1 element'' 
\cite{Tits_order_1}. With this in mind, we will prove the following

{
\renewcommand{\theequation}{A}
\begin{theorem}
\label{A}
All  nonabelian finite simple
groups  of rank $n$ over $\F_q,$  with the\break
 possible exception of the
Ree groups $^2G_2(q),$  have presentations with at most $C_0$  
generators and relations and total  length   at most $C_0(\log n+\log q) ,$
for a constant $C_0$ independent of $n$ and $q$.%
\footnote{Logarithms are  to the base 2.} 
\end{theorem}
\addtocounter{equation}{-1}
}

We estimate (very crudely) that  
  $C_0<1000$; this  reflects  the 
explicit and 
constructive nature  of  our presentations.
 The theorem also holds for all perfect central extensions of the  stated
groups.

  The theorem is interesting in several ways.   It already seems quite
surprising that the alternating and symmetric groups have {\em bounded
presentations}, i.e., with a bounded number of generators and relations
independent of the size of the group,   as in the theorem (this
possibility
    was recently inquired about  in
\cite[p.~281]{CHRR}).  This was even less expected for the groups of Lie
type such as
$\PSL(n,q)$.
   Mann \cite{Mann} conjectured that {\em
every finite simple group $G$ has
a presentation with  $O(\log|G|)$ relations}
(this is discussed at length in \cite[Sec.~2.3]{LS}).
   As a simple\break
 application of
\cite{BGKLP},  this was proved in  \cite{Mann}     with
the possible exception of the twisted rank 1 groups (the Suzuki groups $^2\!
B_2(q)=\Sz(q)$, Ree groups
$^2G_2(q)=R(q)$ and   unitary groups $^2\! A_2(q)=\PSU(3,q)$; in fact,  by
\cite{Suz,HS}  Mann's  result also holds  for  $\Sz(q)$
   and    $ \PSU(3,q)$).
Theorem~\ref{A}
 clearly goes much further than this result,   providing an {\em
absolute} bound on the number of relations. 

The {\em length}  of a presentation is the sum of the number of generators
and the  lengths of the relations as words in  both the generators and their
inverses.
   A presentation for a  group $G$ in the  theorem  will be
called {\em  short} if its length is $O(\log n+\log q) $.
(A significantly different   definition of ``short'', used in
\cite{BGKLP}, involves a bound
$O((\log |G|)^c)$, i.e., $O((n^2\log q)^c)$ if $q>1$.)
When $q>1$, the standard Steinberg
presentation \cite{St1,St2} for  a group in the theorem has length $O(n^4
q^2)$ (with a slightly different bound for some of the twisted groups).
   In \cite{BGKLP}   (when combined with
\cite{Suz,HS} for the cases  $\Sz(q)$
   and    $ \PSU(3,q)$),  the Curtis-Steinberg-Tits presentation
\cite{Cu},\,\cite[Theorem 13.32]{Tits_book} 
was used to prove   that, once again with the possible\break
 exception of
the   groups $^2G_2(q)$, all  finite simple groups $G$  have
presentations of   length
$O((\log  |G|)^2)$ (i.e., $O(( n^2\log q )^2 )$  if $q>1$) -- in fact of
length  $O(
\log |G|)$ for most families of simple groups.   Theorem~\ref{A} 
  clearly also
improves this   substantially.
Note that our $O(\log (nq))$ bound on length is optimal in terms of $n$
and $q$, since there are at least
$cnq/\!\log q\,$ groups 
 of rank at most $n$ over fields
of order at most
$q$
(see Section~\ref{Concluding remarks}).  On the other hand,   whereas our
theorem shows that   nonabelian finite simple groups have  presentations far
shorter than $O(\log |G|)$,
   \cite{BGKLP} (combined with
\cite{Suz,HS})  showed that  every finite group $G$
with no $^2G_2(q)$ composition factor has a presentation of length
$O((\log |G|)^3)$,  where  the constant 3 is best possible.

Theorem~\ref{A} 
 seems counterintuitive, in view of the  standard types
of presentations of simple groups.  We have already mentioned such
presentations for groups of Lie type, but symmetric and alternating
groups are far more familiar.
The  most well-known  presentation  for   $S_n$ is the ``Coxeter
presentation''
\eqn{Coxeter}
(discovered by Moore \cite{Mo}  ten years before Coxeter's birth).
It uses roughly
$n$ generators and $n^2$ relations, and  has  length roughly $  n  ^2 $;
others use
$O( n )$ relations \cite{Mo,Cox}.  We know of no substantially  shorter
presentations in print.    However,
independent of this paper,
\cite{Con} obtained a  presentation of
$S_n$  that is bounded
 and has
 {\em bit-length\/}
$O( \log  n  )$; however,  it is not short in the
sense given above.  (Section~\ref{lengths} contains a definition of
bit-length, together with   a discussion of various notions of lengths of
presentations and their relationships.  For example, it is
straightforward to turn a presentation having  bit-length
$N$  into a presentation having length $4N$, but generally at the
cost of introducing an unbounded number of additional relations.)



Either   bounded {\em or} short presentations for nonabelian  simple groups
go substantially beyond what one might expect.  Obtaining both
simultaneously was a surprise.  By contrast, while abelian simple
groups (i.e., cyclic groups of prime\break
 order) have    bounded
presentations, as well as ones that are short,  they cannot have
presentations satisfying both of these conditions (cf.
Remark~\ref{unbounded cyclic}).  In fact, this\break
 example led us to
believe, initially, that nonabelian simple groups also would not have
short bounded presentations. A hint  that nonabelian finite simple
groups   behave differently from abelian ones came from the Congruence
Subgroup Property, which can be used  to obtain a short bounded
presentation for
$\SL(2,p)$ when $p$ is prime (see Section~\ref{Congruence Subgroup
Property}).  This approach was later replaced by an elementary but
somewhat detailed one (see Section~\ref{PSL(2,p)}), which became
increasingly more intricate as we moved to $\PSL(2,q)$ and other rank 1
groups (Sections~\ref{BCRW trick }-\ref{Presentations for rank 1
groups}).  The
$\PSL(2,p)$ presentation also was  used to obtain a bounded and short
presentation for
$S_n$ (Section~\ref{S{n}}); when combined with presentations for rank~1
groups the latter was then used to prove the theorem for all groups of
large rank (Section~\ref{Bounded presentations for groups of Lie type}).
Our arguments could  have been  considerably
simplified if we had been willing to give up one of the two features of
the presentations in Theorem~\ref{A}.   
 A modified version of Theorem~\ref{A}  is given in \cite{GKKL2} 
 that uses significantly fewer relations in presentations that  are no longer short.

In order to emphasize our special interest in bounded presentations,
we note the following consequence of the above theorem, which follows
immediately from Lemma~\ref{d generators} together with the fact that,
by the classification of the finite simple groups,
every such group can be generated by two elements:

{
\renewcommand{\theequation}{A$'$}
\begin{corollary}
\label{A'}
There is a   constant  $C$ such that every finite simple group$,$
with the possible exception of the Ree groups $^2G_2(q),$  has a
presentation having $2$ generators and   at most
$C$ relations. 
\end{corollary}
\addtocounter{equation}{-1}
}

In \cite{GKKL2}  we
show that
$C\le100$ by ignoring presentation length. 
 The groups $^2G_2(q)$  are
obstacles for all of these results:  no short  {\em or}  bounded
presentation is presently known for them.
We can handle these offending exceptions
by changing our focus either to cohomology or  to profinite
presentations.
We say that a finite or profinite group $G$
has a {\em profinite presentation} with
$d$ generators $X$ and $r$ relations $Y$ if there exists a continuous
epimorphism  $\hat F_d \rightarrow G$, where $\hat F_d$ is the
free profinite group on a set $X$ of $d$ generators, whose kernel
is the minimal closed normal subgroup of $\hat F_d$ containing
    the $ r$-element subset $Y$.  Then
in the category of profinite presentations there are no exceptions:
{
\renewcommand{\theequation}{B}
\begin{theorem}
\label{C}
\label{B}
There is a constant $C$ such that
every finite simple group $G $
has a {\em profinite} presentation
with $2$ generators and at most   $ C $ relations.
\end{theorem}
\addtocounter{equation}{-1}
}

   The quantitative study of profinite presentations of
finite simple groups  was  started in
\cite{Lub1}, motivated by an attempt to prove the Mann-Pyber conjecture
\cite{Mann,Py}, which asserts that {\em the number of normal subgroups
of index
$n$ of the free group
$F_d$ is}   $O ( n^{cd\log n} )$ for some constant $c$.
   Mann \cite{Mann} showed that his conjecture about
    $O(\log |G|)$-relation presentations
for finite simple groups implies the Mann-Pyber conjecture;   and hence he
proved that conjecture except for the twisted rank 1 groups.

   In  \cite{Lub1}  the
Mann-Pyber conjecture  was proved by using  a weaker version of Mann's
conjecture: \emph{ There is a constant $C$ such that every finite
simple group $G$ has a}  profinite \emph{presentation with $2$
generators and at most $C  \log |G|$ relations.}
The crucial ingredient in the proof of the latter result  was a theorem of
Holt   \cite{Holt}:
{\em
   There is a constant $C$ such that$,$ for every
finite  simple group $G,$ every prime $p$ and every  irreducible
$\F_p  G $-module $M,$}
\begin{equation}
\label{Holt''s Theorem}
\dim   H^2(G,M) \leq C(\log_p |G|_p)\dim   M,
\end{equation}
where $|G|_p$ denotes the $p$-part of the integer $|G|$.
In fact, Holt proved 
 that $C$ can be taken to be 2, using  tools including
standard cohomology methods together with results in  \cite{AG}.
Moreover,    Holt
\cite{Holt}  conjectured the following stronger result, which
we will see   is equivalent to   Theorem B:

{
\renewcommand{\theequation}{B$'$}
\begin{theorem}[Holt's Conjecture for simple groups]
\label{HOLT}
\label{B'}
There is a  constant $C$ such that$,$
for every  
finite simple  group $G$$,$ every prime 
$p$ and every  irreducible  $\F_pG$-module $M ,$ $\dim
H^2(G,M)\leq  C\dim  M$.
\end{theorem} 
\addtocounter{equation}{-1}
}

In a subsequent paper  \cite{GKKL}  we will show that the constant in
Theorems~\ref{B}  and B$'$ can be taken to be less than 20 
 (and is
probably even smaller than this).  The  proof uses an
interplay of cohomological {\em and} profinite presentation  results.

More generally, Holt proved that \eqn{Holt''s Theorem} holds for any
finite group $G$ and any faithful irreducible $\F_p G $-module $M$; his
proof reduced this to the simple group case.
Similarly, using a  different reduction,
in \cite{GKKL}  we use Theorem~\ref{B'}
 to  prove a stronger conjecture made
by Holt:

{
\renewcommand{\theequation}{B$''$}
\begin{theorem}[Holt's Conjecture]
\label{C''}
\label{B''}
There is a  constant $C$ such that$,$
for every
finite  group $G$$,$\! every prime 
$p$ and every faithful irreducible $\F_p G $-module $M ,$  $\dim
H^2(G,M)\break\leq  C\dim  M$.
\end{theorem} 
\addtocounter{equation}{-1}
}

As noted above, Theorem~\ref{A}
 does not require  the classification of the finite
simple groups: it only deals with groups having   rank $n$ over $\F_q$.
Of course, by the\break
 classification this  ignores only the 26 sporadic
simple groups. On the other hand,
Theorems~\ref{B},~\ref{B'} and~\ref{B''}  do use the classification, as does Holt's
proof of \eqn{Holt''s Theorem}.  However, unlike many papers by some of
the authors, no classification-dependent internal structural properties
of these  groups are required for the proof.

We emphasize that our methods  are all essentially
elementary:
the needed structural properties  of the groups are standard, and the only
number theory used is ``Bertrand's Postulate'', a weak precursor to the
Prime Number Theorem.  Finally, we note that most of our presentations
contain standard structural aspects of various simple groups:  we do not
eliminate or otherwise fine-tune relations in a manner standard in this
subject (there are many examples of the latter  process in \cite{CoMo}).

\subsection{Outline of the proofs}
\label{Outline of the proofs and the paper}

In Section~\ref{Theorems B and B'}  we will show that
Theorems~\ref{B} and~\ref{B'} are equivalent, and that they follow from
Corollary~\ref{A'} except for the Ree groups
$^2G_2(q)$.  In order to get a smaller constant in   Theorem~\ref{B},
a detailed estimation of  second cohomology groups is needed  using
\eqn{hat r(G) inequalities}. This  is postponed to a later paper
\cite{GKKL}.


   Sections~\ref{Symmetric groups}-\ref{Theorems A and B}
contain the proof  of Theorem A.
Groups of Lie type are built out of symmetric groups (and related Weyl
groups) together with rank 1 groups of Lie type.  Most of our efforts are
directed towards these two classes of groups.

In  Section~\ref{Symmetric groups}  we  obtain   bounded
presentations of  length
$O( \log n )$ for symmetric and alternating groups.
 Our approach  is based on the groups
$\PSL(2,p)$. Therefore we first have to handle that case  in
Sections~\ref{Congruence Subgroup Property} and \ref{PSL(2,p)}.
A short bounded  presentation for $\PSL(2,p)$
already  requires having the presentation encode exponents up to $p$
    in binary; this is accomplished using a trick from
\cite{BKL}: a group-theoretic version of ``Horner's Rule''  \eqn{base
4}.

We combine our   presentation  for  $\PSL(2,p)$ with  the usual embedding
$\PSL(2,p)<S_{p+1}$  in order to obtain a short bounded presentation
for $S_{ p+2}$  whenever $p>3$ is  prime.  Then  we glue together such
presentations for two copies of $S_{ p+2}$ in order to obtain  a short
bounded presentation for the general symmetric group
$S_n$.%

 We deal with
rank 1 groups (especially $\SL(2,q)$, $\SU(3,q)$ and $\Sz(q)$)
in  Section~\ref{Rank 1 groups}.
For these  groups  we have  to deal with several
obstacles:
\begin{itemize}
\item{}  {\em Borel subgroups  do not have short bounded presentations.}
However,  we need to use   a standard presentation of Steinberg
\cite[Sec.~4]{St2}, based on rank 1  BN-pairs,  that  is  built
from a   presentation for a Borel subgroup $B$ together with other data
(cf. Theorem~\ref{Steinberg presentation rank 1}).
(For example, in the case
$\SL(2,q)$ the subgroup   $B$ is a semidirect product
$\F_q^+\semi \F_q^*$ and hence    does not have a short bounded
presentation;  cf.  Remark~\ref{unbounded cyclic}.)

   As a substitute for  a short bounded presentation  for  $B$, we find one
for an {\em infinite} central extension of
$B$, after which almost all of the center is killed when we include the
remaining ingredients in a     presentation for groups such as $\SL(2,q)$.
\item{} Somewhat long relations, such as $x^p=1$, need to be handled, once
again using Horner's Rule.
\item{} Elementary abelian subgroups of order $q$ and arbitrary rank  have
to be {\em shown to be abelian by using only a bounded number of
relations}. For this purpose
   we need to significantly generalize an idea due to  
Baumslag \cite{Bau} and    Campbell, Robertson
and Williams
\cite{CRW1} so that we can even    handle  nonabelian subgroups of
 unitary and Suzuki groups.
\item{} Longer relations, such as $x^{(q-1)/2}=1$, have to be handled.
This is accomplished primarily  by using the minimal polynomial of a field
element of the required order as a replacement for the exponent.
\end{itemize}
Thus, already in rank 1 we have to put significant effort into handling
the restrictions imposed in Theorem~\ref{A}.

The rank 1 presentations, combined with   the Curtis-Steinberg-Tits
 presentations  \cite{Cu},\,\cite[Theorem~13.32]{Tits_book},
 are used  in
Section~\ref{Fixed rank} to  obtain
short  bounded    presentations for  bounded rank groups.  Bounded (but
not short) presentations for   untwisted groups of bounded rank were
obtained  in
\cite{KoLu}.
Starting with our   rank 1 results,
all bounded rank groups  can be readily  handled somewhat as in
\cite{KoLu} (and in fact  this   is essentially  in \cite{BGKLP}).
Fortunately, the rank 1 group $^2G_2(q)$ does not arise as a Levi subgroup
of any higher rank group.

Finally, in Section~\ref{Bounded presentations for groups of Lie type}
we first consider
$\SL(n,q)$, and then   show how to extend the results from both
$\SL(n,q)$ and  bounded rank  to the universal covers of all groups of
Lie type   (except $^2G_2(q)$), using the Curtis-Steinberg-Tits
presentation. For all groups of large rank we go to some additional
effort to increase our sets of generators  and relations somewhat in
order to obtain a presentation in which a suitable  element of the center
is a short word in the generators, and hence   the center  can be factored
out within the constraints of Theorem~\ref{A}..

\subsection{The lengths to which we could go}
\label{lengths}

Table~1 contains a summary of the preceding outline, together with
information concerning different notions of  lengths of a presentation
$\langle  X\mid R\rangle$:
\begin {itemize}
\item[]
\begin {itemize}
\item [length]  = {\em word length}: $  |X| \, + $ sum of the lengths of
the words in $R$  within the free group on $X$;
thus, length refers to strings in the alphabet $X\cup X^{-1}$.  This is
the notion of length used in this paper.  Mild variations,  the same as
ours up to a constant,  are   used by practitioners of the art of
computing explicit presentations for groups.  For example,
\cite{CHRR}  just uses the above length-sum, while
\cite{Do} uses $|R|\,+$ the length-sum.

As noted in \cite[pp.~290-291]{Sim}: ``There is no universal agreement
as to when one presentation is simpler than another.  Usually it
is considered good if the number of generators can be reduced or the
total length of the relators can be made smaller.... The real challenge
is to reduce both the number of generators and the total length of the
presentation.''   A similar notion  appears    in the discussion
of length in \cite[p.~184]{HEO}: ``As a measure of the complexity of a
group presentation, we formally define its size to be the triple
$(|X|,|R|,l)$, where $l$ is the total relator length."
\item [ plength]  = {\em power length}: $|X| \, +$ total
   number of powers appearing  in $R$.
\item [blength]  = {\em bit-length}: total number of  bits required  to
write the presentation, used in  \cite{BGKLP}  and    \cite{Con}.
In particular, all exponents are encoded as binary strings, the sum of
whose lengths enters into the blength, as does the space required  to
enumerate the list of  generators and relations. This definition seems
especially suitable from the point of view of Computer Science.
  It also makes results similar to those given here significantly
easier to prove (see \cite{GKKL2}).
\end {itemize}
\end {itemize}

\smallskip
Some of our presentations can be modified  so as to have
these  {\em other}  lengths,
   as indicated in
Table~1, where we write $q=p^e$ with $p$ prime. However, we definitely do
not
   provide presentations
having all of the indicated types of lengths simultaneously. 

\begin{table}[b]
\label{Sections and lengths}
\begin{center}
\caption{Different lengths}
            $ \begin{array}{|c|c|c|c|c|c|} \hline
$Sec.$ & $groups$&\#$ relations$& $length$ &$ plength$& $bit-length$ \\
\hline 
 3  &  S_n,A_n & O(1) & O(\log n) & O(1)  & O(\log n)\\
4  & $rank$~ 1  & O(1) & O(\log q) & O(e)
    & O(\log q)

\\
5 & $rank$~ n  &O(n^2)& O(n^2\log q) & O(n^2 + e)  & O(n^2 \log n \log q
)
\\
6.1 & \PSL(n,q) &  O(1)& O(\log n  +\log q ) & O( e )   &
O(\log n  + \log  q)\\
6.2& $large$ ~n  & O(1)& O(\log n  +\log q ) &
O(e)  & O(\log n  +\log q )\\
       \hline
             \end{array}$
\end{center}
\end{table}

\smallskip
\smallskip

{\noindent \bf Example\,1:}
$\langle  x \mid  x^n \rangle$  has length $1+n$, plength 2, and
blength
$ 2+[\log n]$.
\smallskip

{\noindent \bf Example\,2:}
If $k=\sum_0^m a_i2^i$ with
$a_i\in\{0,1\}$, the  power $x^k$ can be replaced by
$ \prod _0^m (x^{2^i})^{a_i} =\prod _0^m  x_i ^{a_i}$ for  new
generators  $x_i$ and  relations
$x=x_0, x_i=x_{i-1}^2\break (1\le i\le m) $.

   This multiplies bit-length   by at most a factor of 4 in order to
obtain length.  Thus, we see that  length and bit-length
are essentially the same, {\em except if presentations are to  be
bounded$,$ as in this paper.}  For bounded presentations, bit-length can be
far smaller than length.
Clearly, plength is always at most length.  We also note that the
computer algebra system Magma stores  relations as written: it does not
introduce new relations (as in this example) unless asked to do so.

\smallskip

This example shows that a cyclic group of order $n$ has a presentation of
length $O(\log n)$.  However, there is no such short presentation with a
{\em bounded} number of relations:

\Remark
\label{unbounded cyclic}
   Each family of cyclic groups $C_n$ whose orders are   unbounded 
cannot have bounded presentations of length
$O(\log n)$.
\rm To see this, assume that $C_n$ has a
presentation with a bounded number of generators  $x_1,\dots,x_d$ and
  relations
$w_1, \dots,w_{r }$.  Include the  bounded number of additional relations
$[x_i,x_j]=1$, thereby turning our presentation into one in the category
of abelian groups.  Then every $w_i$ can be viewed as a vector in the
group $\Z^d$ with basis  $x_1,\dots,x_d$.  We can renumber the $w_i$
so that $w_1,\dots,w_{d}$ generate a subgroup of finite index, at least
$n$, in
$\Z^d$.   Then the vectors $w_1,\dots,w_{d}$ are the rows of   a matrix
whose determinant is at least
$n$, and which is also a polynomial of degree $d$ in the entries of the
$w_i$.  Hence, at least one of the entries of some $w_i$ is at least
$n^{1/d}/d!$, which means that the corresponding exponent of some $x_j$
appearing 
in $w_i$ is at least that large. Consequently,  the length of the
presentation is at least
$n^{1/d}/d!$,  and hence is certainly not $O(\log n)$ since $d$ is
bounded.

It follows immediately  that various other groups cannot have short
bounded presentations.  For example, {\em the
$1$-dimensional affine group
$\AGL(1,q)$ cannot have a
bounded presentation of length $O(\log q)$.}  (Namely,  if we mod out
by the commutators of the generators in a presentation for this group
then  we obtain a presentation for a cyclic group of order $q-1$.)

\subsection{Historical remarks}
\label{Historical remarks}

There seems to have been an effort, at least in the 1930's but also as
far back as 1897, to obtain presentations for symmetric and alternating
groups involving as few generators and relations as possible.
For example, Moore
\cite{Mo} gave not only the Coxeter presentation \eqn{Coxeter} for
$S_n$  but also one  with 2 generators and approximately $n$
relations.
Many such presentations  are reproduced in  \cite[Sec.~6.2]{CoMo}.
In \cite[p.~281]{CHRR}  we find the question:  ``Is there a
two-generator presentation for $A_n$ with $k$ relators, where
$k$ is independent of $n$?''  Our Theorem A, as well as  \cite{Con},
    answer  this question.
We note that \cite {CHRR} lists various presentations of small simple
groups, using as one of the  criteria for ``niceness'' the short lengths
of all relations.

There has been a great deal of research on presentations for relatively
small simple groups, and for  small-dimensional quasisimple groups.
Extensive tables are given in
\cite[pp.~134--141]{CoMo}.  Special emphasis is given to $\PSL(2,p)$ and
$\PGL(2,p)$ in \cite[Sec.~7.5]{CoMo}.  An indication  of the large amount
of more recent work of this sort can be found in the references in
\cite{CCHR,CHRR,CR1,CRW1,CRW2}.

   On the other hand,
there are  few general references containing  presentations for
groups of Lie type; cf. \cite[Theorem~13.32]{Cu,St1,St2,Tits_book}
 --
described in
   \cite[Sec. 2.9]{GLS}   -- and their consequences \cite{BGKLP,KoLu}.
However, there are a few specific groups for which presentations
have been published; \cite{CRW2} and \cite{CHLR} are  typical.
The only reference containing a hint in the direction of 
Corollary~\ref{A'}   is
\cite{Wil}, which contains a more precise conjecture than the above
question in \cite{CHRR},  namely, that    the universal central
extension of every finite simple  group has a presentation with 2
generators and 2 relations.

Finally, we mention a recent, very different type of presentation for
some groups of Lie type, 
 initiated in \cite{Ph} and continued very
recently in \cite{BS,GHNS,GHS}.
While these methods can be used to obtain   presentations of unitary
groups with fewer relations than here   (see
\cite{GKKL2}), in this paper we prefer to deal with all classical groups
in a more  uniform manner.

\section{Elementary preliminaries}
\label{Elementary preliminaries}

\subsection{Presentations}
We begin with some elementary observations concerning presentations.
Functions will always act on the left, and we use the notation
$g^h=h^{-1}gh$ and $[g,h]=g^{-1}h^{-1}gh =g^{-1}g^h$.

Let $d(G)$ denote the minimum number of elements needed to generate $G$.

\begin{lemma}
\label{d generators}
If $G=  \langle D  \rangle$ is a finite group  having  a presentation
$\langle X \mid  R \rangle ,$   then $G$ also has   a presentation
$ \langle D'\mid R'\rangle $  such that $|D'|=|D|$ and the natural map
$F_{D'}\to G$  sends  $D'$ bijectively to $D,$ and $|R'|= |D|+|R|.$
In particular$,$  $G$ has  a presentation with  $d(G)$ generators and
$d(G)+|R|$ relations.
\end{lemma}

\Proof
Write each $x\in X$ as a word $v_x(D)$ in $D$ and each $d\in D$ as a
word
$w_d(X)$ in $X$. Let $V(D)=\{v_x(D) \mid x\in X\}$.  Let $d\mapsto d'$
be a bijection from $D$ to  $D'$.

We claim that   $ \langle D'  \mid   d'=w_d(V(D')), \, r(V(D')) =  1,
d\in D, r\in R \rangle  $ is a presentation for   $G$.  For, let
$\tilde G$ be the presented group, and let $\pi\colon F_{D'} \to \tilde
G$ be the natural surjection.

If  $H:=\langle   \pi (V(D'))\rangle  \le\tilde  G$, then $H=\tilde G$
since
$\tilde G$ is generated by the elements $\pi(d') =\pi(w_d( V(D') )  \in
H$, $d\in D$.

Also,  $H$ is a homomorphic image of $G$   since $H$ is
generated by elements $ \pi(v_x(D' ))$, $x\in X$, satisfying the defining
relations for $G$:  if $r\in R$ then    $r(\pi( v_x(D' ) ))=  1 $ holds  in
$\tilde G =H$.

Finally, $G$ is a homomorphic image of $\tilde G$  since the map
$d'\mapsto d$ sends the generators of
$\tilde G$  to generators of $G$ satisfying the same relations as those
satisfied  in  $\tilde G$  by    $\pi(D')$.
   \qedd

\begin{lemma}
\label{Schreier}
Let  $H$ be a subgroup of index $m$ in $G$. 
\begin{itemize}
\item[(i) ]
\begin{itemize}
\item[(a) ]  If $G$ has a presentation  with $x$ generators and $r$
relations$,$ then $H$ has a presentation with 
 $m(x-1)+1$
generators and
$m r$ relations.
\item[(b) ]  If the presentation in {\rm (a)} has length $l,$ then the
corresponding presentation of $H$ has length at most $ml$.
\item[(c) ]   
The length of an element of $H$ with respect to  the above generators of H
 is at most its length with respect to  the original generators of $G$.

\end{itemize}
\item[(ii) ] If $H$  has a presentation with $y$ generators and $s$
relations$,$ then
$G$ has a presentation with $d(G)$ generators and $md(G)+s$  relations
$($independent of~$y)$.
\end{itemize}
\end{lemma}
\Proof  (i) Part (a)  follows from  the standard
Reidemeister-Schreier algorithm
\cite[Secs.~6.3,$\,$6.4]{MKS}:  the given presentation for $G$
produces  an explicit presentation for $H$, with the stated numbers of
generators and relations.  Moreover, for (b) note that each relation of
$G$ of length  $l$ gives rise to $m$ relations for $H$, each of length $l$
in the generators of $H$ (cf. \cite[p.~184]{HEO}).
Finally, for (c) see \cite[p.~58]{Ser}.

(ii)   Let $\pi\colon F_{d(G )} \to G$ be a  surjection with kernel $N$.
As in (i),  $L =  \pi^{-1}(H)  $  is a free group on a set  $X$  of
$n=m(d(G)-1)+1$ generators.  By Lemma~\ref{d generators}, $H$ has a
presentation using the set $\pi (X)$ of these $n$ generators  together
with            $n+s$  relations.   Consequently, $N$ is the normal
closure in $L$ -- and  hence also in  $F_{d(G) }$ -- of a set of
$ n+s  \le md(G)+s$ elements.  \qedd

   \smallskip
We note that we do not know a length version of part (ii) analogous to
(b) or (c).

In the future we will want to identify subgroups of our target group
with subgroups of the group given by a presentation:

\begin{lemma}
\label{It's a subgroup}
\label{it's a subgroup}
\label{It's a group}
Let $\pi\colon F_{X\cup Y} \to G=\langle  X,Y\mid R,S\rangle$ and
$\lambda\colon F_X\to H=\langle  X\mid R\rangle $ be the natural
surjections$,$ where $H$ is finite.~Assume that $\alpha \colon G\to G_0$ is
a homomorphism such that $ \alpha \langle \pi(X )  \rangle\cong H$.~Then
$\langle \pi(X ) \rangle \cong H$.%
\end{lemma}
\Proof
Clearly $\alpha $ sends  $ \langle  \pi(X )  \rangle $  onto
$\alpha ( \langle \pi(X ) \rangle)\cong H $.
If $r\in R$ then $r(\pi(X )  )=1$, so that        $\lambda $ induces a
surjection  $H=\langle  X\mid R\rangle \to\langle \pi(X ) \rangle$.
\qedd

   \smallskip
  From now on we   will be less careful than   in the preceding three
lemmas: we will usually identify the  generators   in a presentation with
their images in the presented group.

\begin{Proposition}
\label{Alex'}
Suppose that $G$ has  a presentation
$\langle X\mid R\rangle $ in which $R$ has total length $L$.
If   $\hat G $   is a perfect group such that
   $\hat G/Z=G $ with  $Z\le Z(\hat G)$ of  prime
order  $p ,$ then $\hat G$ has a presentation
$\langle\hat X\mid \hat R\rangle $
   such that $|\hat X|= |X|+1,$ $ |\hat R|= |X| + |R|+1,$ the  length of
$\hat R$ is less than
$4|X|+ 2L +p|R|,$  and   $\hat X$ contains a generator of $Z$.
\end{Proposition}

\Proof
Since $\hat G/Z=G$,  we can lift the surjection $ \phi\colon F=F_X\to G$,
with kernel $N:= \langle  R^F\rangle $,
to a map
$\hat  \phi\colon F\to \hat G$ (just send each  $x\in X$ to a lift
of  $\phi(x)$  in $\hat G$); since $Z$ is contained in the Frattini
subgroup of $\hat G$, $ \hat \phi$ is surjective.
     Then  $K:=\ker \hat \phi<N  $  and
   $ | N /K|  =p$.

Let $\hat X=X\dot\cup \{y\}$.
We may assume that
$r_1\in R$ is not in $K$, so that $N=K  \langle r_1 \rangle  $.
    For each $r\in R$  choose $e_r$ such that
$0\le e_r<p$ and  $s_r(r_1):=rr_1^{e_r} \in K$.
Then $K_1 :=\langle \hat R^F \rangle\le K$, where
$\hat R:=\{ yr_1^{-1}, [y,x],  s_r(y)  \mid  x\in X, r\in R \}$.

We claim that  $K_1 = K$.
First note that $\langle r_1,K_1\rangle  $ is normal in $F$ since
$[r_1,X]\subseteq K_1$.
Then  $ N = \langle r_1,K_1\rangle  $  since each $r=s_r(y)r_1^{-e_r}\in
\langle r_1,K_1\rangle  $, 
 so that 
   $ N /K_1 = \langle r_1 K_1\rangle  $ with 
$r_1^p =s_{r_1}(y)\in K_1$.  Since $K_1\le K$ and  $ | N /K|  =p$, it
follows that $K_1=K$, as claimed.

Clearly $y$ represents a generator for $Z$.
The total length of $\hat R$  is at most\break
$(L+1 ) +4|X| +L+( p-1) |R|  .$~\qedd

\subsection{Some elementary presentations}
\label{Some elementary presentations}
   The most familiar presentation for the symmetric group $S_n$ is the
``Coxeter presentation''  \cite{Mo}:
   \begin{equation}
   \label{Coxeter}
\begin{array}{lll}
 S_{n}=\langle x_1, \ldots, x_{n-1}
& \hspace{-8pt} |& \hspace{-12pt}
 x_i^2= (x_ix_{i+1})^3=
(x_ix_j)^2=1     \vspace{2pt}
\\ 
  &&\hspace{-12pt}   \mbox{for all possible $i$     and for
$i+2\le j\le n-1$}  \rangle,
\end{array}   \quad \quad
   \end{equation}
   based on the transpositions $(i,i+1)$.
We will  often use the following
    presentation  based on the transpositions $(1,i)$:
   \begin{equation}
   \label{Burnside} 
   S_{n}=
\langle  x_2, \ldots, x_n \mid  x_i^2=
(x_ix_{j})^3= ( x_ix_{j}x_ix_{k})^2=1 
    \mbox{ for distinct $i,j,k$}\rangle . 
   \end{equation}
   This presentation is due to   Burnside \cite[p.~464]{Bur} and Miller
\cite[p.~366]{Mil} in 1911; Burnside describes it as  probably ``the most
symmetrical form into which the\break
 abstract definition [of
$S_{n}$]  can be thrown",  an idea that is crucial for our use of it in
the next section.    Carmichael \cite[p.~169]{Carm2} observed that this
presentation  can  be considerably shortened:  only  the relations
$( x_ix_{j}x_ix_{k})^2=1$
   with
$j=i+1$ need to be used.
   Before we learned of  \eqn{Burnside}, we had been using the following
presentation  we had found, based on the same transpositions, which may
be even more symmetrical than the preceding one:
   $$ S_{n }=\langle x_2, \ldots, x_n \mid x_i^2= (x_ix_j)^3=
( x_i x_j x_k)^4=1
    \mbox{~  for distinct $i,j,k$} \rangle.%
$$
 For alternating groups, see \cite{GKKL2}.

There is   an analogue   of   \eqn{Burnside} for linear
groups.   The usual Steinberg presentation \cite{St1}  for
$\SL(n,p)$, $n\ge 3$,  is as follows (where $e_{ij}$  represents  the
elementary matrix with 1 on the diagonal,
$i,j$-entry 1 and all other entries 0):
\begin{equation}
\label{Elementary SLn}
\hspace{-28pt}
\begin{array}{ll}
\mbox{{\bf Generators:}}\hspace{-5pt}&
\mbox{$e_{ij}$  for   $1\leq i \not = j \leq n$.} \vspace{2pt}
 \\
\mbox{{\bf Relations:}}\hspace{-5pt}&
  \mbox{$e_{ij}^p =1,[e_{ij},e_{ik}]=1
,[e_{ij},e_{kj}]=1,[e_{ij},e_{kl}]=1,$ }  \vspace{2pt}
\\
\hspace{-5pt}&
\mbox{$[e_{ij}, e_{jl}]=e_{il} ~ $  for  all
distinct $~i,j,k,l$.}%
\end{array}
    \end{equation}

\noindent The following variation, suggested by \eqn{Burnside},
allows us to use fewer generators:


   \begin{Proposition}
\label{Martin's Steinberg}
\label{SLn(q)}
$\SL(n,p)$ has the following presentation  when  $n\ge4$$:$
\smallskip\rm

{\noindent \bf Generators:} $e_{1j}$ and $e_{j1}$ for $2 \leq j \leq n$.

\smallskip
{\noindent \bf Relations:}
Whenever  $i,j,k$    are distinct and not $1,$
$e_{1j}^p =e_{i1}^p =[e_{i1},e_{1j}]^p=1 $,
$[e_{1j},e_{1k}]=1 $,
$[e_{i1},e_{k1}]=1 $,
   $[[e_{i1},e_{1j}],e_{k1}]=[[e_{i1},e_{1j}],e_{1k}]=1$,
       $[[e_{i1},e_{1j}],e_{j1}]=e_{i1}$ and
$[e_{1i},[e_{i1},e_{1j}]]=e_{1j}$.
   \end{Proposition}

\Proof
   Let $e_{ij}:=[e_{i1},e_{1j}]$ for $i\ne j$.  Note that
$[e_{i1},e_{1j}]$ commutes with $[e_{k1},e_{1l}]$ because it commutes with
$e_{k1}$ and $e_{1l}$.  It follows easily that
\eqn{Elementary SLn} holds.
\qedd

\medskip

  This motivates  presentations in Theorems~\ref{Bounded SLn}
and
\cite{GKKL2}.

\subsection{Gluing}
\label{Glueing, part 1}
The following lemmas   contain the basic idea behind our methodology.
We start with  a subgroup
$T$ -- with a suitable presentation~-- of  our target group $G$, acting
on some set with boundedly many orbits.  We use this to prove that   a
small number of relations -- each arising from  one of these orbits --
implies a very large number of relations.

   \begin{lemma}
\label{3-transitive}
    Let  $T= \langle D\mid R \rangle$ be a presentation of a  transitive
permutation group acting on   $\{ 1,2,\dots,n\}$
such that
$D_1 \subset D$ projects onto generators of
the stabilizer   $T_1$ of the point $1$.
   Assume that every  ordered triple of distinct points can be
sent to $(1,2,3)$ or $(1,2,4)$ by $T$. Let $v_i $  be words in  $D$
mapping onto permutations sending
$ i\mapsto 1$ for  $i= 2,3, 4$. Then the following
is a presentation for a semidirect product $T\ltimes S_{n+1}$
$($where $S_{n+1}$ acts on $\{ 0,1,\dots, n\})$$:$
\smallskip\rm

{\noindent \bf Generators:} $D$ and $z$ $($viewed as the transposition
$(0,1))$.
\smallskip\rm

{\noindent \bf Relations:} $R$, $z^2 =1$, $[z,D_1] =1$,
$[z, v_2 ]^3=1$, \vspace{2pt}

\hspace{43pt}$([z,v_2][z,v_3])^2=
([z,v_2][z,v_4])^2=1$.
   \end{lemma}

   \Proof
Let $G$ be the group defined by the above presentation.
In view of the natural map $G\to T\ltimes S_{n+1}$ we can apply
Lemma~\ref{It's a subgroup} to  identify $T$ with
$\langle D \rangle $.    That map  sends $z^T$ onto
the set   of all transpositions $(0,i)$, $1\le i\le n$. Since $[z,D_1 ]
=1$ it follows that $T$ acts on $z^T$ as it does on  $\{ 1,2,\dots,n\}$.

The relations $z^2 = [z,v_2]^3= ([z,v_2][z,v_3])^2=
([z,v_2][z,v_4])^2=1$, together with the assumed  transitivity of
$T$, imply that
$
x^2=(xy)^3=(xyxw )^2=1
$
for any distinct $x,y,w\in z^T$.
By \eqn{Burnside},
    $N:=\langle z^T \rangle\cong S_{n+1}$.
 
Clearly $G=TN$.  Since    $G$ maps onto $T \ltimes S_{n+1}$, we have
$G\cong T \ltimes S_{n+1}$.  \qedd

\medskip
Note that,  in fact, $G=C_G(S_{n+1} )\times S_{n+1}\cong T
   \times S_{n+1}.$

   \begin{lemma}
\label{SLn(q)x}
\label{3-transitive part 2}
   Let $T=\langle D \mid R \rangle, T_1, D_1, v_i,$ be as in  {\rm
Lemma~\ref{3-transitive}.}
   Then the following is a presentation for
a semidirect product $T \ltimes \SL({n+1}, p)$
$($where $T$ permutes the last $n$ rows and  columns$)$$:$
\smallskip
\rm

{\noindent\bf Generators:} $D$, $e$ and $f$ $($viewed as $e_{01}$ and
$e_{10} $, cf. Proposition~\ref{Martin's Steinberg}).
\smallskip\rm

{\noindent\bf Relations:} $R$, $[e,D_1]=[f,D_1]=1 , $
    and both $e, e^{v_2}, e^{{v_3}},f,f^{v_2}, f^{{v_3}}$
and  $e, e^{v_2}, e^{{v_4}},$ $f,f^{v_2}, f^{{v_4}}$
   satisfy the relations for  $\SL(4,p)$  in {\rm
Proposition~\ref{Martin's Steinberg} (where each of these ordered
6-tuples plays the role  of
$e_{12},e_{13},e_{14},e_{21},e_{31},e_{41}$).}
   \end{lemma}
\Proof
Let $G$ be the group defined by the above presentation.
By  Lemma~\ref{It's a subgroup}, we can identify $T$ with the
subgroup
$\langle D
\rangle $  of  $G$.

The natural  map  $G\to T \ltimes \SL(n+1,p)$ sends $ e ^T \cup f^T $
onto the set of  all matrices $e_{1i}, e_{i1}$, $2 \le i\le n$.  Since
$[\{e, f\},D_1]=1$,  it follows that
$ |e ^T \cup f^T|=2n $ and that $T$ acts on both
$ e ^T $ and $ f^T $ as it does on  $\{ 1,2,\dots,n\}$.

The relations on
$e, e^{v_2}, e^{{v_3}},f,f^{v_2}, f^{{v_3}}$
and  $e, e^{v_2}, e^{{v_4}},f,f^{v_2}, f^{{v_4}}$,
together with the assumed  transitivity of $T$, imply that
the hypotheses of Proposition~\ref{Martin's Steinberg}  hold for
$N := \langle    e^G, f^G\rangle $.  Then $N\cong \SL({n+1},p)$.
As above,  $G$  is a semidirect product of $N$ and $T$.  \qedd

\section{Symmetric and alternating groups}
\label{Symmetric groups}
\label{Symmetric and alternating groups}
In this section we will show that $S_n$ and $A_n$ have
bounded presentations   of length $O(\log n)$.
By Lemma~\ref{Schreier}(i), it suffices only to consider the symmetric
groups.

We start by obtaining short bounded presentations for  $\PSL(2,p)$
(in Sections~\ref{Congruence Subgroup Property} and \ref{PSL(2,p)}).   We
use these  in order to  deal with the case $n =p+2$ with
$p>3$ prime, and then deduce the general $S_n$  from this case.

\subsection{\boldmath$ \SL(2,p)$ and the
 Congruence Subgroup
Property}
\addcontentsline{toc}{subsection}{\protect\tocsubsection{}{\thesubsection}{$
\SL(2,p)$ and the   Congruence Subgroup Property}}
\addtocontents{toc}{\SkipTocEntry}

   \label{Congruence Subgroup Property}

It is not at all obvious that even the groups $\PSL(2,p)$ have  short
   bounded presentations.  Our first such presentation is based on results
concerning arithmetic groups. 
The idea of using the Congruence Subgroup Property to get bounded
presentations of finite groups has already appeared in \cite{BM,Lub3}.
Here we show that the presentations can be also made to be short.
In the next section we will give a
more explicit and far more elementary presentation.

\begin{Proposition}
\label{CSP}
If $p$ is  an odd prime then $\PSL(2,p)$ has a short bounded presentation.
\end{Proposition}
    \Proof
We sketch a proof  using the Congruence
Subgroup Property.  By \cite{Ser1,Ser2}
and \cite[Theorem~2.6]{Mar},
$\SL(2,\Z[1/2] ) $ is finitely presented, has the Congruence Subgroup
Property  and all its non-central normal  subgroups have finite index.
In fact,  $\SL(2,\Z[1/2] ) $ has a finite presentation  based on
its  generators
\begin{equation}
\label{generators}
\mbox{$u=
   \begin{pmatrix} 
        1 & 1 \\
        0 & 1 \\
     \end{pmatrix}
   $,
   $~t=
   \begin{pmatrix}
       \,\,\,\,0 & 1 \\
        -1 & 0 \\
     \end{pmatrix}
~ $
   and $~h_2=
   \begin{pmatrix} 
        1/2 & 0 \\
        0 & 2 \\
     \end{pmatrix}
   $. }
\end{equation}
Adding the additional relation
$u^p=1$  and using the above properties of normal subgroups   produces  a
bounded presentation for
$\SL(2,p)$ for $p>2.$   While
this presentation appears not to be short,  it can be converted to
a bounded presentation of  length
$O(\log p)$ by using  a group-theoretic version of {\em
Horner's Rule\/} \cite[p.~512]{BKL}:
\begin{equation}
\label{base 4}
\begin{array}{lll}
\mbox{If  $2< m \le p$  then $u^m = \prod_i(u^{m_i})^{h_2^i}=
u^{m_0} h_2 ^{-1}u^{m_1} h_2^{-1} u^{m_2} \cdots h_2^{-1} u^{m_k}
h_2^{k}$}  \vspace{2pt}
\\
\mbox{has length $O(\log p),$ where $k=[ \log_4 m]$   and  $m
= \sum _{i=0}^{k}
   m_i4^i $  in base $4$.}
\end{array}
\end{equation}
Namely, since $u^{h_2}=u^4$ in \eqn{generators},
\begin{equation}
\label{collapse}
u^m= u^{m_0}   \cdots (u^{m_k})^{h_2^k}=  u^{m_0}
{h_2^{-1}} u^{m_1}  {h_2}\cdot {h_2^{-2}} u^{m_2}  {h_2^2} \cdot
{h_2^{-3}}\cdots {h_2}^{-k} u^{m_k} {h_2^k},
\end{equation}
 which collapses as in \eqn{base 4}.

This produces a short bounded presentation of $\SL(2,p)$. 
In order to obtain one for 
$\PSL(2,p)$, we write $-1$ as a short word as in \cite{BKL} 
(cf. Lemma~\ref{PSL2 length}) and add one
further relation to kill it. 
 \qedd
 
\smallskip

This method of {\em converting to a short bounded presentation}
will be used often in the rest of this paper.

   We note that the argument used in
the proposition can be greatly generalized:   the Congruence Subgroup
Property (combined with Margulis' normal subgroup theorem) can be used in
a similar manner to provide   short bounded presentations   for various
other families of finite simple groups of fixed rank over suitable fields
(cf. \cite{Lub3}). 

Sunday \cite{Sun} used the presentation in
\cite{BM}, which was obtained using  the\break
 Congruence
Subgroup Property as in the above proposition, in order to produce
the following    presentation for
$\PSL(2,p)$    of bounded plength  and  bit-length $O(\log p)$,  though
not of bounded length:
$$
\langle  u, t \mid  u^p = 1, t^2 = (ut)^3 ,  ( u^4 t u^{(p+1)/2} t )^2 = 1 \rangle .
$$
There is no presentation for this group  with smaller $|X|+|R|$; by
\cite{Sch}, this follows from the fact that  the Schur multiplier of
$\PSL(2,p)$ has order 2. In view of
\eqn{base 4},   this presentation can be modified to a  short  bounded
presentation by including $h_2$ as a new generator,  and using the
additional relations
$u^{h_2}=u^4$  and $h_2  t= u^{\bar 2}  (u^{  2})^t u^{\bar 2}  $,
where $\bar 2 =(p+1)/2$ is the  integer   between 1 and $ p$
such that $2 \cdot  \bar 2   \equiv 1~$(mod $ p)$.

   \subsection{{\boldmath
$\PSL(2,p)$} }
\addcontentsline{toc}{subsection}
{\protect\tocsubsection{}{\thesubsection}{$\PSL(2,p)$}}
\addtocontents{toc}{\SkipTocEntry}

   \label{PSL(2,p)}
We have separated these groups from the general case because we need them
for symmetric groups   and they are less
involved than the general
$\PSL(2,q)$.  The latter groups require a lot more preparation (in
Sections~\ref{BCRW trick } and \ref{Abelian p-groups}),
due to headaches caused by field extensions.
   The ideas used in this section reappear in
Sections~\ref{Central extensions of Borel subgroups} and  \ref{new
PSL(2,q)} with many more bells and whistles.



Fix a prime $p > 3$ and a generator $j$  of $\F_p^*$.
   Let $\bar 2$ be as above,  let   $ \bar j$ be  the integer
between 1 and
$ p$ such that $ j \cdot  \bar j \equiv 1~$(mod $ p)$,
and define
\begin{equation}
\label{more generators}
\mbox{
$h=
   \begin{pmatrix}
       \bar j & 0 \\
        0 &  j \\
     \end{pmatrix}
   $. }
\end{equation}

\begin{Theorem}
\label{PSL2 short bounded}
If $p>3$ is prime then  $\PSL(2,p) $ has the following  
presentation of length
$O(\log p) $     with $4 $ generators and $8$ relations$,$  where 
integer exponents  are   handled using using
$\eqn{base 4}$$:$\rm
\smallskip

{\noindent \bf Generators:}  $u,  h_2, h, t$
$($which  we think of as the matrices in \eqn{generators}
and  \eqn{more generators}).

\smallskip
{\noindent\bf Relations:}
\begin{itemize}
\item []
\begin{enumerate}
\item$ u^p=1$.   
%
%
\item
     $ u^{h_2}=u^4 $.   
%
\item  $ u^h=u^{j^2} $.
\item  $t^2=1.$
%
%
%
\item  $ h^t = h^{-1}$.
\vspace{2pt}
\item $ t =u u^t u$, $h_2  t= u^{\bar 2}  (u^{  2})^t u^{\bar 2} $,  $h t=
u^{\bar j}  (u^{  j})^t u^{\bar j}  $.
\end{enumerate}
\end{itemize}
\end{Theorem}

\Proof
  Let  $G$  be the group defined by this presentation.  Then
$\PSL(2,p)$  is an epimorphic image of $G$.
By  (6), $G=\langle  u , u^t  \rangle$.
Then $G $ is perfect
since  $u^3 = [u,h_2] \equiv 1$ (mod $G'$) by (2), while  $u^p=1$.

We claim that  $Z:=  \langle  z \rangle\le Z(G)$, where
    $z: = h^{(p-1)/2}$. For,  if we write    $j^{p-1} = 1 + mp$  for some
integer $m$, then (3) implies that $u^z = u^{(j^2)^{(p-1)/2}} = u
u^{mp} = u $ by (1). Thus   $[u,z]=1$.  Similarly,  $[  u^t,z]=1 $
since  $(u^t)^{h^{-1}}=(u^t)^{j^2}$ by (4) and (5).

By (1) and (3), in $G/Z$ the subgroup generated by $u$ and $h$ is
isomorphic to a Borel subgroup  of   $\PSL(2,p)$, while
the subgroup generated by $t$  and $h$ is    dihedral of order $ p-1 $
by (3), (4) and (5). Also, by  (6),   if   $~k=1~ $   or   $ ~ j $ then  
there are integers $k' , k'' $ and $k'''$    such that
$$
\mbox{$(u^k)^t=u^{k'} h^{k''} t u^{k'''} .$
      }
$$
\noindent
Conjugating these relations using  $\langle h \rangle $   produces
the same type of relation  whenever
$ 1\le k\le p-1$ (by (3),  $\langle h \rangle $ acts on the nontrivial
powers of $u$  with only  two  orbits, with orbit  representatives  $u$
and $u^j).$

Now we have verified the standard Steinberg presentation for $\PSL(2,p)$
\cite[Sec.~4]{St2}.  (That presentation, as well as  more general  ones
that are  also based on groups with a BN-pair of rank 1, will be very
useful for us; see the beginning of Section~\ref{Presentations for rank 1
groups} for more details.) Thus,
$G/Z \cong \PSL(2,p) $, and   $G$ is a perfect central extension
of  $\PSL(2,p)$.  The only perfect central extensions of $\PSL(2,p)$ are
itself and $\SL(2,p)$, so that
   $G\cong \PSL(2,p)  $  by  (4).  \qedd

\smallskip

\Remark
\label{shorter PSL(2,p)}\rm  
Changing (4)  to
$[t^2,u]=1$ produces a presentation for $ \SL(2,p).$

The  presentation  for
$ \PSL(2,p) $ at the end of Section~\ref{Congruence Subgroup Property}
uses fewer relations.

\medskip

For later use,    the following
      observation is crucial
\cite[p.~512]{BKL}:

\begin{Lemma}
\label{PSL2 length}
If $ p>3$ then every element of  
 $\PSL(2,p) $ and  $\SL(2,p) $ can be
written as a word of length
$O(\log p)$ in the generating set $\{u,  h_2, h, t\}$ used above.
\end{Lemma}

\Proof Every element  can be written as the product of a
bounded number  of powers of $u$ and $ u^t$  (cf.
\cite[p.~512]{BKL}).  Now use \eqn{base 4}  and relation (2).~\qedd

\medskip
This lemma also holds for the presentation in  Proposition~\ref{CSP}.

   \subsection{{\boldmath
$S_{p+2}$} }
\addcontentsline{toc}{subsection}
{\protect\tocsubsection{}{\thesubsection}{$S_{p+2}$}}
\addtocontents{toc}{\SkipTocEntry}

\label{S{p+2}}

   ~ In this section we will obtain  a short bounded presentation
   for  $S_{p+2} $       when $p>3$   is  prime.

Start with a short   bounded presentation of
$\PSL(2,p)=\langle X\mid R \rangle $  such as the one  in
Theorem~\ref{PSL2 short bounded},  so that $X$  {\em consists of   elements
$u ,h,h_2, t$ corresponding to} \eqn{generators}
and  \eqn{more generators}.

In the action of  $\PSL(2,p)$  on the projective line
$\{0,1,\dots, p-1,\infty\}$, 
 we identify $a\in \F_p$ with the 1-space $\langle (a,1) \rangle $
of $\F_p^2$, so
that
$$u= (0,1,\dots, p-1) $$
and $t$ interchanges  $0$ and  $\infty$.  Then
$$B:=\langle   u^t,h \rangle  $$
is a Borel
subgroup fixing
$0$ while $\langle   u,h \rangle  $ fixes  $\infty$.~There are exactly two
orbits of $\PSL(2,p)$ on ordered triples of points, with representatives
$(0,-1,\infty )$ and
$(0,-s,\infty )$ for any non-square $s$   in $\F_p^*$.

We view  $S_{p+2}$ as acting on the set  $\{ \star   , 0, 1,\dots,
p-1,\infty \}$, with  $\PSL(2,p)$  fixing~$\star$.

\begin{Lemma}
\label{long Sp+2}
   $\!\!S_{p+2}$
has  the following   presentation with $5$ generators and $15$
relations$:$

 \rm
\smallskip
{\noindent \bf Generators:} $X$ and $z$ (viewed as the transposition
$(\star   ,0)$).

{\noindent\bf Relations:}
\begin{itemize}
\item []
\begin{enumerate}
\item $R$.
\item
$[z, u^t]= [z,h]=1$.
\item $z^2=1$.
\item
$[z,t]^3=1$.
\item  $([z,t][z,u ])^2= ([z,t][z,u^s])^2=1$.
\item $(zu)^{p+1} =1$.
\end{enumerate}
\end{itemize}
\end{Lemma}
\Proof
Let $G$ be the group defined by the above presentation. It is easy to see
that  $G$ maps onto
$S_{p+2}$.
(The various relations reflect the facts that, when viewed as
permutations,
$z$ commutes with elements having   support   disjoint from
$\{\star,0\}$;    the supports of $z$ and $t$ have just one common
point, as do the supports of $z$ and $u$; and
$[(\star,0)(\star,\infty)(\star,0)(\star,-1)]^2=1=
[(\star,0)(\star,\infty)(\star,0)(\star,-s)]^2$.)

By Lemma~\ref{it's a subgroup}, $G$ has a subgroup
we can identify with 
$T=\langle X\rangle \cong\PSL(2,p)$. By (2), $z$ commutes with the stabilizer
 $B=\langle   u^t,h \rangle$  in $T$
of 0.

By (2)--(5), we can use Lemma~\ref{3-transitive} (with $t,u^{} ,u^{s}$ 
playing the roles of
$v_2,v_3, v_4$:  they send $\infty\mapsto 0,-1\mapsto 0$  and  $-s\mapsto
0,$ respectively).  Thus
$N:=\langle  z^T\rangle  \cong S_{p+2}$.  Also,
   $N\unlhd \langle T,z\rangle =G$.

Now $u = u^{p+1}\equiv (zu)^{p+1}=1$ (mod $N$)   by (6), so  that    $G/N$
is a homomorphic image of $\PSL(2,p)$ such that $u$ is mapped to $1$.
Then $G/N=1$
   since  $\PSL(2,p)$ is simple.~\qedd

\Remark\rm
This is a bounded presentation with bounded plength and
  bit-length  $O(\log p)$. This presentation can be modified to
a presentation of length $O(\log p)$ with  $O(\log p)$  relations
(cf. Section~\ref{lengths}, Example 2).
However,  while it is bounded  
{\em it is not  short\/}  due to the long
relation
   $(zu)^{p+1}=1$ (note that  $u^s$ can be made short using \eqn{base 4}).
Nevertheless, we will use the above presentation to obtain   a short bounded
presentation in the next theorem:
we will {\em deduce} this long relation from shorter ones
by using  a second copy of $\PSL(2,p)$.

\Remark
\label{move relations}
\rm
One of our main methods has been to use 
(nearly)  multiply transitive
permutation groups in order to ``move'' some relations 
and hence 
to deduce many
others. For example, Lemma~\ref{3-transitive} is concerned with groups
that are 3-transitive or close to 3-transitive.
 We  will use this idea again in Section~\ref{SL(n,q)} and
\cite{GKKL2}.  Moreover, in \cite{GKKL} 
and \cite{GKKL2} we even use sharply 2-transitive groups for this purpose
rather than more highly transitive groups.

In preliminary versions of the above lemma, when we sought bounded
presentations but had not  yet approached ones that are both bounded and
short, we proceeded somewhat differently.  For example, we used a simple
and standard bounded presentation for the 3-transitive group $\PGL(2,p)$
instead   of using   its subgroup
$\PSL(2,p)$.

An even earlier  approach to symmetric groups started with a much more
complicated special case, but
was   nevertheless interesting.
First  we  obtained  bounded  presentations for
$\SL(k,2)$ for infinitely many $k$, as follows.  The Steinberg
group   $\St_4(\Z\langle x_1 , \ldots, x_d\rangle)$  over the free
associative (non-commutative) ring   on $d$ variables is finitely
presented for every $d\ge 4$  \cite{KM}.
Consequently,  if $R$ is a   boundedly presented quotient ring   of
$\Z\langle x_1 , \ldots, x_d\rangle$,
then $\St_4(R)$ is boundedly presented.
   Note that   the matrix algebra
$M_t(\F_q)$ is boundedly presented (e.g.,   $ M_t(\F_p)  =  \langle A, B \mid
{A^t=B^t=0,} ~{BA+A^{t-1} B^{t-1} =1},~   {p1=0}  \rangle $; alternatively, a
bounded presentation for   the general $ M_t(\F_q)$ can be obtained using
the fact that it is a cyclic algebra).
Therefore, $\St_4(M_t(\F_q)) \cong \St_{4t}(\F_q) \cong\SL({4t},q)$  also is
boundedly presented.
Since $\SL({4t},2)$  has boundedly many orbits on   $4$-sets  of
nonzero vectors,
    we can  proceed in a manner similar to the proof  of
Lemma~\ref{3-transitive} 
in order  to obtain a short presentation for $S_{2^{4t}-1}$ by
using the presentation  \eqn{Coxeter}
 instead of    \eqn{Burnside}.

\Remark
\label{more elements needed}  \rm
We   {\em  need  the following additional elements
of $\PSL(2,p)$} for future
use:
\begin{itemize}
\item
$k$  of order  $(p+1)/2$, acting with two cycles on the projective line
(this is a generator of a non-split torus),
    and
\item
$l$ such that $1$ and $l$  are representatives of the two
$B,\langle k\rangle$   double cosets in  $\PSL(2,p)$. Thus,
\begin{equation}
\label{double coset}
\langle X \rangle =B \langle k\rangle \cup B l\langle k\rangle.
\end{equation}
\end{itemize}
In order to relate these to our previous generators,
  we use  Lemma~\ref{PSL2 length} in order to
{\em obtain two short relations expressing $ k$ and $l$ in terms of $X$}.

 We emphasize that {\em  the order of $ k$ follows from the  
relations in}   Theorem~\ref{PSL2 short bounded}.
We now use this observation.

\Remark \rm
The element $zu$ acts
as a
$(p+1)$-cycle   $(\star   ,0 )( 0,1,2,\dots, p-1) =   
(\star
,0,1,2,\dots, p-1)$, and
$z^t(zu)z^t$ is the $(p+1)$-cycle $(0,1,2,\dots,p-1,\infty)$.
(Note that we are multiplying permutations from right to left.)
Let
$\sigma$ be a permutation of the points of this  line such that 
$$
k^\sigma =\big(z^t(zu)z^t   \big)^2~ \mbox{ and }
 ~\sigma=(\star)(0)(1,l(0),\dots).
$$
Such a permutation  exists 
 because   $k$  fixes
$\star $ and acts as the product of two disjoint
   $(p+1)/2$-cycles  with representatives $0$ and
$l(0)$, while   $(z^t(zu)z^t)^2$   fixes
$\star$ and  acts as the product of two disjoint
   $(p+1)/2$-cycles  with representatives $0$ and $1$.
(N.B. -- It seems  that
$\sigma$  cannot be expressed as a short word in our generating set, but
this will not cause any difficulties because $\sigma$ will not be used
in any explicit manner.)
Thus we have a relation
$$
(zu)^2 = z^t k^{\sigma^{}} z^t.
$$

With all of this notation, we are ready to prove the following

\begin{Theorem}
\label{Sp+2 Theorem}
   $S_{p+2} $ has a   presentation of length $O(\log
p)$   with $9 $ generators and $2 6$ relations. 
\end{Theorem}

\Proof
We will use two
copies  $\langle X\mid R  \rangle$ and $\langle X'\mid R'  \rangle$
of the presentation  for $\PSL(2,p)$ obtained in
Theorem~\ref{PSL2 short bounded}  and  Remark~\ref{more elements needed}
  (where $X\cap X'=\emptyset$); the natural map $x\mapsto x'$,
$x\in X$,   extends  
to an isomorphism ``prime'':$\langle X\mid R 
\rangle \to  \langle X'\mid R'  \rangle$.
We will show that  $S_{p+2}$ has the following presentation: 

\smallskip
{\noindent \bf Generators:}    $X\cup X'$  and
$z$. (We view $z$ as $(\star,0)$  and   $X'$ as   $X^\sigma$
   in the preceding remark.)

\smallskip
{\noindent\bf Relations:}
\begin{itemize}
\item []
\begin{enumerate}
\item $R\cup R'$.
\item $[z,  u^t]= [z,h]= [z,  u'{}^{t'}]=   [z,h']=1 $.
\item
$z^2 =1 $.
\item
$[z,t]^3=1 $.

\item  $([z,t][z,u ])^2= ([z,t][z,u^s])^2=1$. 
\item    $[z,u^{-1}]=[z,l'{}^{-1}]$      ~
    and  ~ $(zu)^2 = z^t k'{} z^t.$

\end{enumerate}
\end{itemize}
\smallskip

\noindent
Relations (6)  are especially important since they link the two copies
  for $\PSL(2,p)$.  Note that the  
first relation can be rewritten
$z^{u^{-1}}=z^{l'{}^{-1} } $.

Let $G$ be the group defined by this short, bounded presentation.
By the preceding remark, there is a surjection  $
\pi\colon
G\to
S_{p+2}$. (Note that this uses the permutation $\sigma$ 
that is
not actually
in the presentation!  The first relation in (6) is based on the fact
 that both $u^{-1}$ and $(l^\sigma)^{-1}$ conjugate 
$(\star,0)$ to $(\star,1)$
since $l^\sigma(0 ) =\sigma^{-1}(l(0))=1$.) 

We claim that  {\em $X$ and $z$ satisfy the relations  in} Lemma~\ref{long
Sp+2}. This is clear for Lemma~\ref{long Sp+2}(1)-(5).
Moreover, the present relation (6) implies Lemma~\ref{long Sp+2}(6):
$$
(zu)^{p+1} = (z^t k'{} z^t)^{(p+1)/2} = z^t {k'{}}^{(p+1)/2} z^t =1
$$
since $k'{}{}^{(p+1)/2}=1$ follows from the relations $R'$.
Thus, by    Lemma~\ref{It's a group}  we can identify
$G^*:= \langle z, X \rangle $  with  $S_{p+2}$
 in such a way that $z$ is a transposition  and $|z^{\langle z, X\rangle
}|=p+1$. Then $G^*=\langle z^ {\langle z, X  \rangle}  \rangle$.

We  can also identify
$\langle X' \rangle $ with $\PSL(2,p)$.
In place of \eqn{double coset} we  will use
$\langle X'\rangle=B' \langle k'\rangle \cup B' l'{}^{-1} \langle
k'\rangle,$   where $B'=\langle   u'{}^{t'},h' \rangle  $. Since $z$ commutes
with   $B'$ by (2), it  follows that
$$z^{\langle X'\rangle} = z^{\langle k'\rangle}\cup
(z^{l'{}^{-1}})^{\langle k'\rangle}=z^{\langle k'\rangle}\cup (z^{u^{-1}
})^{\langle k'\rangle}\subseteq z^{\langle z, X\rangle}$$
by both parts of (6).  
Then 
$$
p+1=|\pi( z^{\langle z,X'\rangle})|\le | z^{\langle
z,X'\rangle} |\le |  z^{\langle z,X\rangle} |=p+1,
$$
so that   
$z^{\langle z,X'\rangle} =z^{\langle z,X\rangle}$ 
and hence
$G^*=\langle z^{ \langle
z,X\rangle } \rangle   =\langle z^{\langle  z,X ,X'\rangle } \rangle  
=\langle z^G \rangle  \unlhd  G$.

Clearly, $G/G^*$ is a homomorphic  image of 
$\langle X'\rangle\cong
\PSL(2,p)$  in which $k'{}$ is mapped  to $1$,
   since  $k'{} \in \langle X,z \rangle= G^*$  by  (6).
    Therefore $G/G^*=1$   since
    $  \PSL(2,p) $ is simple.

Replace $ k'$  and  $l'$   by the words
representing them in order to obtain 
a presentation with 9 generators and 26 relations.
\qedd 

\smallskip
\smallskip


We need further  properties of  the preceding generators in order to
handle general symmetric groups:

\begin{Lemma}
\label{short words elements}
The following elements of $S_{p+2}$ can be written as words of length
   $O(\log p)$ in $X\cup X'$ whenever $1\le i\le p-1 $$:$
$$
\begin{array}{lllll}
\pi_{i} :=(\star    ,0,1,\dots,i)
\quad &
& \lambda_{i} :=(i+1,\dots,p-1,\infty) \vspace{2pt}
\\
\theta_i:=(i,i+1)
 &&
\tau:=(\star, \infty) \vspace{2pt}
\\\theta :  =(p-1, \infty) .
\end{array}
$$

\end{Lemma}
\Proof
In each case we will express the permutation in terms of
elements that can be moved into one of our copies of $\PSL(2,p)$
 and hence  can be
written as short words by Lemma~\ref{PSL2 length}.
We start with the elements
$$
u=(0,1,\dots,p-1),    k'{}\!=(z^t(zu)z^t)^2=(0,1,\dots,p-1,\infty)^2
,\,    z=(\star,0), \,   t = (\infty ,0) \cdots
$$
in $ X\cup X'.$  Then
$$
\begin{array}{llllccc}
\tau  = z^t   &\theta \,\, =  (0,\infty) ^u =  (z^u)^\tau
&
   \qquad(\star,1)=uzu^{-1}
\vspace{4pt}\\
   (0, 1) = (uzu^{-1})^z  \quad&  \theta_i =u^i(uzu^{-1})^z u^{-i} .
\end{array}
$$
If $1\le 2j-1\le p-1$, then
$$
(\infty,0,1,\dots,2j-1) =
(0,1,\dots,p-1)^{2j}(0,1,\dots,p-1,\infty)^{-2j}  =  u^{2j} k'{}^{-j}
,$$
 so that
$$
\begin{array}{llllccc}
\pi_{2j-1}  = z^t (\infty ,0,1,\dots,2j-1)z^t  = z^t
u^{2j} k'{}^{-j}z^t
\\ 
\pi_{2 j} = \pi_{2j-1} \theta_{2j-1}
\vspace{3pt}\\
\lambda_{p-1-2j} =(0,1,\dots,p-1,\infty)^{-2j}(0,1,\dots,p-1)^{2j}
=  k'{}^{-j}u^{2j}
\\
\lambda_{ (p-1-2j )-1}=\theta_{p-1-2j}\lambda_{p-1-2j}.
\end{array}
$$
By Lemma~\ref{PSL2 length}, the elements  $u^{i}, u^{2j}\in \langle
X\rangle $  and
$k'{}^{j} 
\in \langle X'\rangle $ can be expressed as short words in our
generators.~\qedd

\smallskip

\smallskip

Note that the same argument shows that all cycles $(i,i+1,\dots , j) $, and
all permutations with bounded support, can be written as short words.
 


   \subsection{{\boldmath
$S_{n}$} }
\addcontentsline{toc}{subsection}
{\protect\tocsubsection{}{\thesubsection}{$S_{n}$}}
\addtocontents{toc}{\SkipTocEntry}

 \label{S{n}} ~
 Finally, we  consider all symmetric groups:

\begin{Theorem}
\label{All symmetric groups}
    $S_n$  has a   presentation of length $O(\log n)$   with $18 $
generators and $58$ relations. 

\end{Theorem}
\Proof  
By Bertrand's Postulate, if $n \ge 6$  then there is a prime
$p$ such that\break
 $(n-2)/2   < p  < n -2$.
 Then  $i : = 2(p+2) -n -2 \ge1$.
 We  will use  two copies $\langle X\mid R\rangle$ and $\langle \tilde
X\mid
\tilde R\rangle$ of the   presentation of
$S_{p+2}$  in Theorem~\ref{Sp+2 Theorem} (where  $X\cap \tilde
X=\emptyset$),  and identify them along a
subgroup $S_{i+2}$.

The   map $x\mapsto \tilde  x $,
$x\in X$,   extends  
to an isomorphism ``tilde'':$\langle X\mid R 
\rangle \to  \langle \tilde  X \mid  \tilde  R   \rangle$.
Then the  elements
$ \lambda _i,\theta_i, \tilde \lambda_i,  \tilde \tau  $ 
are  as in  the preceding lemma. 
We will show that  
 $S_{n}$ has the following presentation:%
\smallskip

{\noindent \bf Generators:} $X, \tilde X $.
\smallskip


{\noindent\bf Relations:}
\begin{itemize}
\item []
\begin{enumerate}
\item $R\cup \tilde R  $.
\item $z = \tilde z$, $\pi_i=\tilde \pi_i$.
%
%
\item $[\lambda _i, \tilde \lambda_i]=[\lambda_i, \tilde \tau ]=
[\theta_i  ,\tilde \lambda_i]=[  \theta_i ,\tilde\tau]=1$.
\end{enumerate}
\end{itemize}
\smallskip

Let $G$ be the group defined by this bounded presentation of length
$O(\log n)$.    By Lemma~\ref{It's a subgroup},  $G$  is generated by two
copies  $\langle  X\rangle $ and $\langle \tilde X\rangle $  of  $S_{p+2} $.  These are linked as
follows:
\begin{equation}
\label{3 rows}
\begin{array}{ccccccc}
                        &  \hspace{-4pt}  i+1,  &   i+2,  & \dots  & p-1,&
\infty
\\
\star   , 0, 1, 2,\dots, i-1, i, \\
                        &  \hspace{-4pt} \widetilde{i+1}, &  \widetilde{i+2},
&\dots  &
\widetilde{p-1}, & \widetilde{{^{^{^{}}}}\infty ^{}  }
\end{array}
\end{equation}
Namely, the first copy of $S_{p+2}$ acts on the first and second rows,
and the second copy on the second and  third   rows;
the relations (2) identify the intersection, which is $S_{i+2}$ acting on
the second row.  Thus, we will abuse notation and
identify
$s =\tilde s $ for
$s\in\{\star,0,\dots,i\} $.
Note that $\Sym\{ i,  i+1, i+2,  \dots , p-1,  \infty\}=
\langle  \lambda_i,\theta_i \rangle   \cong S_{p-i+1} $ commutes with
$\Sym {
   \{
    \widetilde{i+1} ,  \dots , \widetilde{p-1},
\widetilde{{^{^{^{}}}}\infty}, \star
\}  } =\langle
\tilde\lambda_i, \tilde\tau\rangle  \cong S_{p-i+1} $,  by (3),
since  $i\ge1$.

In order to prove that this defines $S_n$ we will
use the Coxeter presentation \eqn{Coxeter} with the following
 ordering  of  our $n$  points:
\begin{equation}
\label{arrangement}
{\infty},   p-1,\ldots ,i+1,i=\tilde i,i-1,\ldots ,2,1,0,\star ,
\widetilde{{^{^{^{}}}}\infty ^{} },\widetilde{p-1} ,
\dots , \widetilde{i+1}.
\end{equation}
The first copy of  $S_{p+2}$ acts from $\infty$ to
$\star$, and the second copy acts  from $i = \tilde i $ to
$\widetilde{i+1}$.

For any adjacent points  $a,b$ in
\eqn{arrangement}, there is a well-defined element $g_{(a,b)}$ of
$G $ that acts as a transposition in at least one  of our two symmetric
groups
$S_{p+2}$.   The elements  $g_{(a,b)}$ generate $G$
since they generate  both symmetric groups.

The natural   surjection $G\to S_n$ maps
   $g_{(a,b)}$ to the transposition $(a,b)$. In order to prove that
$G\cong S_n$, it suffices to prove that the   elements $g_{(a,b)}$ satisfy
the Coxeter relations \eqn{Coxeter}. 
\vspace{-5pt}
\begin{table}[h]
\label{Sn table}
\begin{center}
\caption{Pairs of pairs of points}
            $ \begin{array}{|c|c|c|c|c|c|} \hline
   a,b   &  c,d &$ reason $\\
\hline 
   \infty,\dots,i & \infty,\dots, \star  &
S_{p+2}
\\
   \infty,\dots,i &\star,\dots,  \widetilde{i+1}  &
   [ \langle  \lambda_i,\theta_i \rangle , 
\langle \tilde \lambda_i,\tilde\tau\rangle ]=1
$ by (3)$
\\
i, i-1, \dots , \widetilde{i+1}   & i,i-1, \dots, \widetilde{i+1}  &
S_{p+2}
\\
       \hline
             \end{array}$
\end{center}
\end{table}

Consider four points $a,b,c,d$ occurring in this order in
\eqn{arrangement}, with $a,b $  and $c,d$ two distinct
 pairs of  adjacent  points.  (We allow the possibility $b=c$.)
We  tabulate the various cases in Table~2, including the reason that
the required relation in  \eqn{Coxeter} holds.

Thus, we obtain a  short  presentation for $S_n$ with 
 $2*9$
generators and $2*2 6 + 2+ 4$   relations.
\qedd


\begin{Corollary}
\label{All symmetric groups: plength}
    $S_n$ has a
bounded presentation   with bounded  plength.

\end{Corollary}
\Proof
Use the preceding presentation {\em without}  converting the various
powers of $u$ and $k'$ into short words
   in the proof of Lemma~\ref{short words elements}.
(Of course, we should also use Lemma~\ref{long Sp+2} in place of
Theorem~\ref{Sp+2 Theorem}, since the latter is not  relevant here.)
~\qedd

\vspace{2pt}
\Remark
\rm
Neither $n$ nor $\log n$ appears explicitly in this presentation.  However,
$n$ is encoded in the presentation (in binary):    $p$ certainly is, and
$i$ is  encoded through the word for $\lambda_i$ 
in Lemma~\ref{short words elements}.
In other words, $n$ can be reconstructed from  $p$ and $i$, which are
visible in the presentation.

\Remark
\label{parities}
\label{needed permutations} \rm
   For use   in Section~\ref{SL(n,q)}, for $n\ge 8$ we calculate  a number of
additional elements of $S_n$ that can be written  as short words in
   \rm  $X\cup  \tilde X $.

 (i) In  Section~\ref{SL(n,q)}  we will use the following
{\em renumbering  of our points}:
$$
\star, 0,1,\dots,i-1  , \widetilde{i+1} ,\dots, \widetilde{p-1} ,
\widetilde{{^{^{^{}}}}\infty  } ,i, \dots,p-1,\infty ~ = ~
{\bf2},{\bf3}, {\bf4},{\bf5},
\dots,{\bf n} ,
$$
where $ {\bf n }: =n+1$:  we will be considering the stabilizer of $\bf 1$
in the symmetric group $S_{\bf n} = \Sym\{{\bf 1},\dots,\bf n\}$.  Then
the permutations
$$
\begin{array}{lllll}
     z\hspace{1pt}=({\bf2},{\bf3})   &
   \sigma :=\pi_i \tilde \lambda_{i-1} \lambda_{i-1}=({\bf2},{\bf3},
{\bf4},{\bf5},  \dots,{\bf n})
 \vspace{2pt}
\\

z^ {\sigma^{-4}}=({\bf 6},{\bf 7}) & z
z^{\sigma ^{-1}}z^{\sigma^{-2}} = ({\bf 2},{\bf 3},{\bf4},{\bf 5})
 \vspace{2pt}
\\
 (z z^{\sigma^{-1} }z^{\sigma^{-2} } z^{\sigma^{-3} } )^{-1} \sigma =
({\bf6},{\bf7},  \dots,{\bf n})\quad
&
 (z z^{\sigma^{-1} }z^{\sigma^{-2} } z^{\sigma^{-3} }z^{\sigma^{-4} } )^{-1}
\sigma = ({\bf7},  \dots,{\bf n})
\end{array}
$$
{\em  can be written as words of length $ O(\log {\bf n})$ in our generators.}

(ii) We will also need information concerning parities  of some of our
generators:   the stabilizer of the point $\star =\bf  2  $ is
generated by a set $X \cup  \tilde X$  of even permutations
(the generators of our two copies of  $\PSL(2,p) $) together with
the odd permutation   $ y:=z^{\sigma^{-1}}  = ({\bf 3},{\bf 4})$ that we
{\em  will now include among our $18+1$ generators.} 

(iii)
 Since
$yy^{\sigma^{-1}} = ({\bf 3},{\bf 4},{\bf 5})$
 we see that, in the alternating
group  $A_{\bf n-1}$  on $\{{\bf2},{\bf3}, {\bf4},{\bf5},
\dots,{\bf n} \}$,
$\langle X, \tilde X , yy^{\sigma^{-1}}\rangle $ is the stabilizer of
{\bf 2},   and
 $\langle  u^t,h, \tilde X , (yy^{\sigma^{-1}})^{\sigma^{-1}}  \rangle
\break
 = \langle B,\tilde X,
(yy^{\sigma^{-1}})^{\sigma^{-1}}
\rangle $ is the stabilizer  of both  {\bf 2} and {\bf 3}.

   \begin{Corollary}
\label{extensions of alternating}
Every perfect central extension of $A_n$ has a presentation
with  a bounded number of relations and length $O(\log n)$.
\end{Corollary}
\Proof
For $A_n$ this follows from Theorem~\ref{All symmetric groups} and
Lemma~\ref{Schreier}(i).
When  ${n>6}$ the center of the universal perfect central extension of
$ A_n$ has  order 2  \cite[Theorem 5.2.3]{GLS}, so that
Proposition~\ref{Alex'} can be applied.~\qedd

\section{Rank 1 groups }
\label{Rank 1 groups}
Simple groups of Lie type are built out of two pieces:     rank 1
subgroups such as $\SL(2,q)$, and Weyl groups such as $S_n$.
 Therefore, as stated earlier,  we have placed  special emphasis on
 these types of groups.


   \subsection{The BCRW-trick}
   \label{BCRW trick }
The following lemma is based on an idea in \cite{CRW1}, which \cite{KoLu}
called {\em the CRW-trick}.
A year after first using this  we learned 
 that a similar idea had
been used earlier (for a different goal) in \cite{Bau}.  We therefore
rename  this {\em the BCRW-trick}.
 This makes remarkably effective use of a simple observation:
if $u,v,w$ are group elements such that $w=uv$, then any element commuting with
two of them also commutes with the third.  This  observation  is used
to  show  that, under suitable conditions,  a subgroup can be proven  to be  {\em
abelian} using a very   small number of relations (see
Section~\ref{Some combinatorics behind the BCRW Trick} for a different
viewpoint  and a generalization):
   \begin{lemma}
\label{BCRW1 trick}
Let $G$ be a group generated by elements  $u, a,  b$   satisfying the
following relations$:$
\smallskip
\begin{itemize}
\item []
\begin{enumerate}
\item $[ a, b] =1;$
\item $[u,u^{  a}] =1$ and $[u,u^{ b}] =1;$ and
\item each of $u,$ $u^{ a},$  $u^{ b}$ can be expressed as a word in the
other two.
\end{enumerate}
\end{itemize}
\smallskip
Then $\langle u^{\langle a,b\rangle }  \rangle$ is abelian.

   \end{lemma}

It is clear that (3) and the first part of (2) imply the second part.
We have stated the result in the present manner in order to emphasize
the symmetry between $a$ and $b$.
\smallskip

\Proof
Define $u_{s,t} := u^{ a^s b^t}$ for all integers $s,t$. It suffices to
prove that
\begin{equation}\label{BCRW1inductive goal}
    \mbox{\em
     $u =u_{0,0}$ and $u_{s,t}$ commute for all
integers
$s,t$.   \qquad\qquad\qquad\qquad \quad\quad\quad
   }
\end{equation}

First note that, by  (1) and (3),
\begin{equation}\label{(A)}
\hspace{-25pt}
\begin{array}{lll}
   \mbox {\em
If $u$ commutes with two members of a
triple
$u_{i,j},$ $u_{i+1,j},$ $u_{i,j+1}$  }
\\
  \mbox {\em then  it commutes with the
third}.
\end{array}
\end{equation}

Also, {\em if $u$ commutes with $u_{i,j}$ then $u^{ b^{-j}  a^{-i}  } =
u_{-i,-j}$ commutes with $u $} (using (1)).
It follows that
\begin{equation}\label{(B)}
\hspace{-25pt}
\begin{array}{lll}
\mbox {\em If $u$ commutes with two members of a
triple
$u_{i,j},$ $u_{i-1,j},$ $u_{i,j-1}$  }
\\
\mbox {\em  then it commutes with the
third.}
\end{array}
\end{equation}
   For,  $u$ commutes with
two members of the  triple $u_{-i,-j},$ $u_{-i+1,-j},$ $u_{-i,-j+1}$,
   so by \eqn{(A)} it commutes with the third and hence with the third
member of the triple in  \eqn{(B)}.

We will prove   \eqn{BCRW1inductive goal}    by induction on $n:=|s|+|t|
$.

When  $n=1$ our hypothesis states that
$u=u_{0,0} $  commutes with  $u_{1,0}$  and $u_{0,1}$; these and
$u_{-1,0}$  and $u_{0,-1}$  are the only possible elements having $n=1$,
and we have already seen that $u$ commutes with the preceding two
elements.

Assume that \eqn{BCRW1inductive goal} holds for some $n\ge 1$ and
   consider  nonnegative integers  $s,t$ such that $s+t=n+1$.

\smallskip
   Case: $0= s\le  t= n+1$. By induction, $u$ commutes with
the last   two members of the triple 
$u_{1,n},$ $u_{0,n},$ $u_{1,n-1}$,  hence   by  \eqn{(B)} with the first.
Then, by induction and \eqn{(A)}, $u$ commutes with the
first  two members of the triple
$u_{0,n},$ $u_{1,n},$ $u_{0,n+1}$,  hence  with the third.
\smallskip

   Case: $1\le s\le  t$.   By induction and  \eqn{(B)}, $u$ commutes with
the last  two members of the triple 
$u_{s,t}$,   $u_{s-1,t}$,   $u_{s,t-1}$, hence with the first.
\smallskip

By the symmetry of our hypotheses, this completes the proof when
$s,t\ge0,$ and hence also when
$s,t\le0$.   By symmetry, it remains to   consider the case  $s<0 <t$,
where we  use \eqn{(A)} and the triple
   $u_{s  ,t-1}, u_{s+1 ,t-1},  u_{s ,t}$.
\qedd
\smallskip\smallskip

We will need various  nonabelian versions of the preceding lemma.  The
most straightforward  one is as follows:


   \begin{lemma}
\label{BCRW1 trick nonabelian}
Let  $U_0$ and $W_0$ be subgroups of  a group $G,$ and let $u, w, a,  b$
be elements of $G$ satisfying the following conditions$:$ 
\begin{enumerate}
\item $[ a, b]=1;$
\item each of $u ,  u^{ a} ,  u^{ b}$ can be expressed
as a word in the other two together with  elements of $U_0;$
\item each of $w ,  w^{ a} ,  w^{ b}$  can be expressed
as a word in the other two together with  elements of  $W_0 ;$
\item $[u,w]=[u^a,w]=1;$
\item $U_0^a=U_0^b=U_0$ and $W_0^a=W_0^b=W_0;$ and
\item $[U_0,w]=[u,W_0]=1$.  
\end{enumerate}
Then  $ [\langle u^{\langle a,b\rangle }\rangle,\langle w^{\langle
a,b\rangle }
\rangle]=1 $.
\end{lemma}

This time we have not preserved the symmetry between $a$ and $b$; but
note that $[u^b,w]=1$ by (2), (4)  and  (6).

\smallskip
\Proof
Define $u_{s,t} := u^{ a^s b^t}$ and $w_{s,t} := w^{ a^s b^t}$ for all
integers $s,t$. It suffices to  prove that
\begin{equation}
\label{inductive goal-noncommutative}
    \mbox{\em
     $u_{s,t}$ and $w$ commute  for all integers $s,t$.
\qquad\qquad\qquad\qquad \quad\quad\quad
    }
\end{equation}

The proof parallels that of Lemma~\ref{BCRW1 trick}.
This time,
\begin{equation}
\label{Aw}
\hspace{-6pt}
\hspace{-25pt}
\begin{array}{lll}
\mbox{\em If  $u\!$
commutes with two members of a triple
$w_{i,j}, w_{i+1,j}, w_{i,j+1}\!$}
\\
  \mbox{\em then it
commutes with the third.}
\\
\mbox{\em If $w $
commutes with two members of a triple
$u_{i,j}, u_{i+1,j}, u_{i,j+1}\!$}
\\
  \mbox{\em then it  commutes with the
third.}%
\end{array}
\end{equation}

 For the first of these, by (1),  (3) and (5)  each of the three
indicated  elements is a word in the others together with elements of
$W_0^a=W_0^b=W_0$; now use (6).  Similarly, the second assertion follows
from  (1),  (2), (5) and (6). 

Even more trivially, in view of (1), {\em if $u $ commutes with
$w_{i,j}$  then
    $u^{ b^{-j}  a^{-i}  } =  u_{-i,-j}$ commutes with $w $};
    and {\em if
     $w$ commutes with $u_{i,j}$   then
$w^{ b^{-j}  a^{-i}  } =  w_{-i,-j}$ commutes with $u $}.
It follows that
\begin{equation}\label{Bw}
\hspace{-25pt}
\begin{array}{lll}
\mbox {\em If $w$ commutes with two members of a
triple
$u_{i,j},$ $u_{i-1,j},$ $u_{i,j-1}$  }
\\
\mbox {\em  then it commutes with the
third.}
\end{array}
\end{equation}
Namely,   $u$ commutes with
two members of the  triple $w_{-i,-j},$ $w_{-i+1,-j},$ $w_{-i,-j+1}$,
hence, by \eqn{Aw}, $u$  commutes with the third, and hence  $w$
commutes  with the third member of the  triple in  \eqn{Bw}.

As before we  prove   \eqn{inductive goal-noncommutative}    by
induction on $n:=|s|+|t|$.

When  $n=0$,    $u=u_{0,0}$ and $w=w_{0,0}$ commute by (4).

Consider the case $n=1$. By   (4),  $w$ commutes with the
first two members of the triple $u_{0,0},$ $u_{ 1,0 },$ $u_{0,1}$, and hence
by
\eqn{Aw}  with
$u_{0,1}$.  Then $w$    commutes with the first and
  last     members of the triple
$u_{0,1},$ $u_{-1,1},$ $u_{0,0}$, and  hence  by \eqn{Bw} with $u_{-1,1}$.
Then $w$ also commutes with the   last two members of the triple
$u_{-1,0},$ $u_{0,0},$ $u_{-1,1}$, and  hence  by \eqn{Aw}  with
$u_{-1,0}$.  Finally, $w$ commutes with the first two members of the triple
$u_{0,0},$ $u_{-1 ,0},$ $u_{0,-1}$, and  hence  by \eqn{Aw}  with
$u_{ 0,-1}$.

Assume that \eqn{inductive goal-noncommutative} holds for some $n\ge 1$,
and  consider  nonnegative integers  $s,t$ such that $s+t=n+1$.%
\smallskip

   Case: $0= s\le  t= n+1$. By induction, $w$ commutes with
the last   two members of the triple 
$u_{1,n},$ $u_{0,n},$ $u_{1,n-1}$, and  hence  by \eqn{Bw}  with the
first.   Then, by induction and \eqn{Aw}, $w$ commutes with the
first  two members of the triple
$u_{0,n},$ $u_{1,n},$ $u_{0,n+1}$,  hence  with the third.

\smallskip

   Case: $1\le s\le  t$.   By induction and  \eqn{Bw}, $w$ commutes with
the last  two members of the triple 
$u_{s,t}$,   $u_{s-1,t}$,   $u_{s,t-1}$, hence with the first.
\smallskip

By the symmetry of our hypotheses, this completes the proof when
$s,t\ge0.$  When
$s,t\le0$, proceed  as in the above two cases but negate all subscripts.
By symmetry, it remains to   consider the case
$s<0 <t$,  where we  use  \eqn{Aw} and the triple
   $u_{s  ,t-1}, u_{s+1 ,t-1},  u_{s ,t}$ to prove that  $w$ commutes with
$ u_{s ,t}$.~\qedd

\smallskip
In the proof of Proposition~\ref{SU3 Borel}  we will use this lemma for 
more than one choice of
$U_0$ and
$W_0$.  The simplest use has $W_0=1$  and  $w$ in the abelian group $U_0$.

\subsection{Some combinatorics behind the BCRW-trick}
\label{Some combinatorics behind the BCRW Trick}
The proof of  Lemma~\ref{BCRW1 trick} was concerned primarily with triples of
subscripts. In  order to handle Suzuki and Ree groups, we will need to
generalize this subscriptology.    For these
purposes, and in order both  to see the connection with the original BCRW
Trick  and to {\em use} the original one, we state the desired situation in
some generality.

\begin{definition}
\label{irreducible}
\label{closed set}\rm
Let $\mathcal{M}$ be a finite collection of finite sets of vectors in
$\Z^k$.

(1)
A set $S\subseteq \Z^k$ is called  {\em closed}   provided that,
for any $v\in \Z^k$ and any $M \in \mathcal{M}$, if $S$
contains all but perhaps one  of the vectors in the translate $v+M=\{v+m
\,|
\,m\in M\}$ then $S$  contains the entire translate.
\smallskip

(2)
Fix   a length function $|v|\ge 0$ for vectors
$v\in \Z^k$ such that each ball about 0 contains only finitely many
vectors. 

  A vector
$v\in
\Z^k$ is called   {\em reducible} if there exist
$u\in \Z^k$ and
$M\in \mathcal{M}$ such that $v$ is the {\em unique} longest vector in
$ u+M$. In this case we will say that  $v$ {\em is reducible  via  $u+M$.}

A vector is called  {\em irreducible} if
it is not reducible.
\smallskip

(3)
The collection $\mathcal{M}$ is called of {\em finite type}  (with respect
to the length function  $|\cdot |$)   if    there are only finitely many
irreducible vectors in $\Z^k$.
\smallskip

(4)
If a collection is denoted  
$\mathcal{M}_{\mathcal{X}}$, then $\mathcal{M}_ {-\mathcal{X}}$
denotes $\{ -M\mid M\in  \mathcal{M}_{\mathcal{X}}\}$.
\end{definition}

\begin{exam}
\label{pm BCRW}\rm
In order to relate this to Lemmas~\ref{BCRW1 trick} and~\ref{BCRW1 trick nonabelian}, we use
$\mathcal{M}_{\rm BCRW}=\big\{\{(0,0),(1,0),(0,1)\}\big\}$, consisting of
a  single triple
in $ \Z^2$, together with
$$\mathcal{M}_{\pm \rm BCRW}  \hspace{-1pt} = \hspace{-1pt} 
\mathcal{M}_{\rm BCRW} \cup \mathcal{M}_{- \rm BCRW} \hspace{-1pt}  =
\hspace{-1pt} 
\big\{\{(0,0),(1,0),(0,1)\},
\{(0,0),(-1,0),(0,-1)\}\big\}.
$$
The two collections  $\mathcal{M}_{\pm \rm BCRW}$   arise  naturally in
Lemma~\ref{BCRW1 trick}, as follows.

 Let $u=v$ and the group $H:=\Z^2$
acts  on
$U =\langle  u^{\langle a,b\rangle }\rangle $ via
$u^{(i,j)} =u^{ a^i b^j}.$
The  set $\mathcal{M}_{\rm BCRW}=\{ M\}$   corresponds to the relation
$u^b=u u^a$, and hence to the fact that each   element in
$\{u,u^a,u^b\} = u^M=\{ u^m \mid m\in M \}$  lies in the subgroup
generated  by the remaining elements (as in  Lemma~\ref{BCRW1 trick}(3)).
If $v\in H$ then $(u^M)^v=u^{v+M}$.  Hence,
the translates of
$\mathcal{M}_{\rm BCRW}$ are exactly the triples of subscripts appearing in
\eqn{(A)}, while   the translates of
$\mathcal{M}_{\pm\rm BCRW}$ are exactly the triples of subscripts appearing in
\eqn{(A)}  and \eqn{(B)}.

Let  $S:=\{h\in H \mid [u^h,u]=1 \}$.
Since $u^{v+(0,1)}=u^{ v} u^{v+(1,0)}$ 
for any $v\in H$, if $u$ commutes with two of the elements 
$u^{v+(0,1)}$, $u^{ v}$, $ u^{v+(1,0)}$
then it commutes with the third, so that    $S$
is closed under
$\mathcal{M}_{\rm BCRW}$. Moreover,  if $[u^h,u] =1 $ then  
$[u^{-h},u]=1$, so that   $S = -S$, and hence  $S$ is also
closed under
$\mathcal{M}_{-\rm BCRW}$.

 We used the length function  $|(x,y)|=|x|+|y|$ in
the  proof of Lemma~\ref{BCRW1 trick}.   Moreover, the way we used 
\eqn{(A)}  and \eqn{(B)} emphasized the fact
that $\Z^2$ is the smallest $\mathcal{M} _{\pm\rm BCRW}$-closed subset of $\Z^2$
containing the vectors
$ (0,0),(\pm1,0),(0,\pm1)$
(a notion generalized below in  Lemma~\ref{irreducible = all}).

The notion of ``reduced" was involved in our induction:  commuting with
two elements of a triple of group elements implies commuting with the
third; in the induction,  the third vector was ``longer'' than the other
two. Our  hypotheses gave us the above  initial set of 5 irreducible
elements of
$\Z^2$,  from which we   ``generated'' all of $\Z^2$ using the notion of
closed.

The situation in  Lemma~\ref{BCRW1 trick nonabelian}
was very similar.  There we implicitly used $S:=\{h\in H \mid [u^h,w]=1
\}$, and our hypotheses were designed to imply that $S$ is 
closed under
$\mathcal{M}_{\pm\rm BCRW}$.  We will formalize this somewhat more in
Lemma~\ref{commutator and systems}.
\end{exam}

\begin{exam}\rm
\label{Original BCRW triples}
We now reconsider the set $\mathcal{M} _{\rm BCRW}$ in the preceding example,
but this time using the Euclidean length $|v|^2=x^2+y^2$, $v=(x,y)$. We
wish to determine when  $v$ is irreducible.
   Note that

if $y < 0< x$  then $(x,y)$ reduces via $(x-1,y) ,(x,y) ,(x-1,y+1)$,

if $x < 0< y$ then $(x,y)$ reduces via $(x,y-1) ,(x+1,y-1), (x,y)$,

if $x,y<0$ then $(x,y)$ reduces via $(x,y),(x+1,y),(x,y+1)$,

if $0\le y< x-1$ then $(x,y)$ reduces via $(x-1,y),(x,y),(x-1,y+1)$, and

if $0\le x< y-1$ then $(x,y)$ reduces via $(x,y-1),(x+1,y-1),(x,y)$.

\noindent
(Compare the proof of Lemma~\ref{BCRW1 trick}.)
The only vectors $(x,y)$  not covered by these observations are
$(0,y)$ for $y \leq0$,
$(x,0)$ for $x \leq0$, and
$(x,y)$ for $ x,y \geq 0$ and $|x-y| \leq 1$,
   all of which are irreducible.

For Lemma~\ref{BCRW1 trick} we needed to have a vector $(x,y)$ such that both
$(x,y)$ and $(-x,-y)$ have the above form. It is clear that for this to
happen one of $x,y$  has to be $0$ and the other one $0,1$ or $-1$.  This
proves that \emph{the only irreducible vectors for  $\mathcal{M}_{\rm
BCRW}$ are the $5$ points inside a ball of radius $1$ around $(0,0)$.}

\end{exam}

The next results deal with the ways these  notions
will be used later.
First we  
formalize the use in  a group-theoretic setting  already hinted
at in Example~\ref{pm BCRW}:

\begin{lemma}
\label{commutator and systems}
Let $\mathcal{M_U}$ and $\mathcal{M_W}$ be  finite collections of
finite sets of elements of a subgroup 
$H  \cong \Z^k$  of a group $G$. Let
$U_0,W_0\le G$ and
$u, w\in G$ satisfy  the following conditions$:$
\begin{enumerate}
\item If $M \in \mathcal{M_U}$ and $h\in M $ then
$u^h$ can be expressed as a word in $U_0$ and $\left\{u^l \mid l\in M
\setminus \{h\} \right\}$$;$
\item If $M \in \mathcal{M_W}$ and $h\in M $ then
$w^h$ can be expressed as a word in $W_0$ and $\left\{w^l \mid l\in M
\setminus \{h\} \right\}$$;$   and
\item $[U_0^H ,w]=[u,W_0^H ]= 1$.

\end{enumerate}
Then $S:=\left\{ h \in H \mid [u^h,w]=1 \right\}$
is closed under $\mathcal{M_U}$ and $\mathcal{M_{-W}}$.
\end{lemma}

Condition (1) essentially says that, for each $M\in \mathcal{M_U}$, there
is a relation among the elements $\{ u^h \mid h \in M\}$  that  can be
solved for each of these elements in terms of the others. Of course, we
will want $S$ to be nonempty, which will always be forced by having
$[u,w]=1$.

\smallskip
\smallskip

\Proof 
We will write $H$ additively.
First  note that  $S$ is closed under $\mathcal{M_U}$.  
For, if  $h\in M\in
\mathcal{M_U}$ and $v \in H$, by (2) we can express $u^{v+h}$  as a word
in
$U_0^v$ and $\left\{u^{v+l} \mid l\in M \setminus \{h\} \right\}$.
Therefore, by (3), if $w$ commutes with
$u^{v+l}$   for each $l\in M \setminus \{h\}$  then
it commutes with
$u^{v+h}$.

The preceding
argument shows that  $\bar S:=\left\{ h \in H \mid [u,w^h]=1 \right\}$ is
closed under $\mathcal{M_W}$, so that     $S = -\bar S $ is also closed
under
$\mathcal{M_{-W}}$  
(see Definition~\ref{irreducible}(4)).~\qedd
\smallskip
\smallskip

Next we observe how the notions ``closed'' and ``finite type'' 
(Definition~\ref{irreducible}(1),(3)) 
will be used  (compare the proofs of 
Lemmas~\ref{BCRW1 trick} and \ref{BCRW1 trick nonabelian}):

\begin {lemma}
\label{irreducible = all}
Fix a length function as in \emph{Definition~\ref{irreducible}(2)}. 
If $S$ is a closed subset of
  $\Z^k$   containing the set $S_0$ of all irreducible
vectors for $\mathcal{M},$ then $S=\Z^k$.
\end {lemma}
\Proof
  List all vectors in $\Z^k$ by non-decreasing length: $ v_1,v_2,\dots$
(recall the finiteness assumption in Definition~\ref{Some combinatorics
behind the BCRW Trick}(2)).
Using induction on $k$ we will show that $v_k\in S$. If $v_k$ is 
irreducible  it is in $S_0$ and therefore in $S$. If $v_k$ is reducible 
then there is a translate
$v+M$,  $v\in \Z^k$,
$ M\in \mathcal{M}$, in which  $v_k$ is the longest vector.
By   induction the other vectors in $v+M$
are in $S$. Since $S$ is closed, also  $v_k\in S$.
\qedd

\begin{corollary}
If the set $\mathcal{M}$ is of finite type then there is a finite set
$S_0\subseteq \Z^k$ such that$,$ if a subset  $S\subseteq \Z^k$ is
$\mathcal{M}$-closed and contains $S_0,$ then $S=\Z^k$.
\end{corollary}

Although  the length function does not appear in the statement of the
corollary  it was needed in the proof. Moreover: 
\smallskip

\begin{remark}
If we only consider lengths arising from  positive definite bilinear forms on
$\Z^k$, then the property that $\mathcal{M}$ is of finite type does not
depend on the specific choice of length.
 Of course, the precise  number of irreducible vectors does depend on
the  choice of  length function.

\end{remark}

\subsection {Central extensions of Borel subgroups}
\label {Central extensions of Borel subgroups}
In this section we apply the\break results in the preceding sections  to
go one step closer to simple groups of rank~1:
in Propositions~\ref{SL2(q) Borel short bounded presentation}, \ref{SU3 Borel}
and \ref{Sz Borel} we will provide short bounded presentations for {\em
infinite} central extensions of their Borel subgroups. (Recall from
Remark~\ref{unbounded cyclic} that Borel subgroups do not have such
presentations.)

\addtocounter{subsubsection}{-1}
\addtocontents{toc}{\SkipTocEntry}
\subsubsection{Polynomial notation}
\label{Polynomial Notation}

 \label{power notation}
 \label{FIELD notation}
All of our presentations for groups of Lie type involve dealing with
polynomials over finite fields. 
\smallskip

1. 
  If $\gamma$  lies in an extension field of $\F_p$, then $m_\gamma (x)$ will denote
its minimal polynomial over $\F_p$. 
If $\delta \in \F_p[\gamma]$, let $f_{\delta;\gamma}(x)\in\F_p[x]$ with 
$f _{\delta;\gamma}(\gamma)=\delta$
and $\deg  f _{\delta;\gamma} <\deg m_{\g}$.
\smallskip

2. 
For any polynomial $g(x) = \sum_0^e g_i x^i\in
\Z [x]$  and 
any two elements $u,h $ in a group $G$,
define    ``powers'' as follows, multiplying 
in the stated order:
\begin{equation}
\label{how polynomials act}
[[u^{g (x)  }]]_h =  (u^{g_0})  (u^{g_1})^{h^{ 1}}
\cdots (u^{g_e})^{h^{ e}} \! ,
\end{equation}
so that
\begin{equation}
\label{Horner for polynomials}
[[u^{g (x)  }]]_h =
u^{g_0} h^{-1}u^{g_1} h^{-1} u^{g_2} \cdots h^{-1} u^{g_e}
h^{e}.
\end{equation}
If  $0\le g_i<p$  then  
this is a word of length $O(pe)$ in  $u$, $h$;
 \eqn{collapse}
explains the collapse of the exponents here  (cf.
\cite[p.~512]{BKL}).

Moreover, if  
 $0\le g_i<p$ and  some $h_2\in G$ satisfies
$u^{h_2}=u^4$ or $u^{h_2}=u^2$, then
we can apply (a version of)
\eqn{base 4} for each $u^{g_i}$, and hence we can  write
$[[u^{g (x)  }]]_h$ as a word of length $O(e \log p)$ in $u$,  $h_2$ and
$h$.
Note that  
$[[u^f(x)]]_h$ then depends only on the image of $f(x)$ in $\F_p[x]$.

Care is needed  with \eqn{how polynomials act}: we cannot use it  to
calculate
$[[u^{g (x) +k(x) } ]]_h$ or $[[u^{g (x)  k(x) } ]]_h$
in the obvious manner
unless we know that all of the
    conjugates $u^{h^i}$  commute with one another.  If this
occurs and $u^p=1$, then 
$[[u^{g (x) +k(x) } ]]_h=[[u^{g (x)}]]_h [[u^{ k(x) } ]]_h$, and
     $[[u^{f (x) g(x) } ]]_h =1$ if $[[u^{f (x) }]]_h =1$.

We will usually use these
formulas when $h$ is such  that the conjugates $u^{h^i}$
generate an isomorphic image of the additive group of $\F_q$ and the
action of 
$h$ corresponds to   multiplication by some $\gamma$  on
$\F_q=\F_p[\gamma]$. In this case,   we have
$[[u^{m_\gamma(x)} ]]_h=1$,
and an isomorphism is given by
$$
 \mbox{$ \delta \mapsto [[u^{f_{\delta;\gamma}(x)}]]_h ~$   for  $~\delta
\in
\F_q  $}.
$$
(Relations  used in   Propositions~\ref{nonabelian Borel short
bounded presentation}, \ref{Sz Borel}  and \ref{Ree Borel}    contain
variations on this situation  in   nonabelian settings.)

\subsubsection{Special linear groups}
\label{Abelian p-groups}
Let $\zeta $ be  a generator of   $\F_q^*$,
where $q=p^e $ for a prime $p$.
Our first result (Proposition~\ref{SL2(q) Borel short bounded
presentation}) requires additional data that exists only when $q \ne
2,3,5,9$:  we need  $a, b\in \F_q ^*$   such that
\begin{equation}
\label{a and b}
\mbox{ $b^2 =   a^2 + 1 ~$ and
$~\F_q=\F_p[a^2]$.}
\end{equation}
Such elements   exist for the stated $q$: if $\zeta ^2+1$ and $\zeta ^2 -1$
are nonsquares then $\zeta ^4 -1$ is a nonzero square $a^2$;  and   this
generates
$\F_q$  for the stated $q$ since $\zeta ^4$   is in no proper subfield.
We emphasize that we will be working  with $a^2$ rather than
$\z$; the latter will be needed only to obtain an element of the desired
order $q-1$.
  See  Lemma~\ref{AGL(1,q) short bounded presentation}
for a variation on the present theme that does not require the equation
$b^2=a^2+1$.

 We base our presentation on   the following
matrices:
\begin{equation}
\label{SL2 generators}
   u=\left(\begin{array}{cc} 1 & 1 \\ 0 & 1\end{array}\right) \!,   ~~
   h_\star= \left(\begin{array}{cc} \star^{-1} & 0 \\
\!\!\!\! 0&\star\end{array}\right)
~~\mbox{ for $\star\in\{2,a ,b,\z\}$.}
\end{equation}
(If $p=2$ discard $h_2$.)
Then $B:=\langle u,h_\zeta \rangle $ \emph{is a Borel subgroup  of}
$\SL(2,q)$, of order $q(q-1)$.

With this preparation,
using  the polynomial
notation in Section~4.3.0   we can consider the following
presentation.

\smallskip
{\noindent \bf Generators:} $u,  h_2, h_a, h_b, h_\zeta $.
(If $p=2$ discard $h_2$.)

{\noindent\bf Relations:}
\begin{itemize}
\item []
\begin{enumerate}
\item $[h_\star, h_{\bullet}] =1$ for $\star,\bullet \in \{2,a,b,\zeta \}$.
\item
    $u^{h_b} = u  u^{h_a} =   u^{h_a}u$.
\item $u^{h_2}=u^4$.
\item $u^p=1$.
\item $[[u^{m_{a^2}(x)} ]]_ {h_a } =1$.
\vspace{2pt}
\item $u^{h_\zeta }= [[u^{f_ {\zeta^2 ; a^2   } (x)} ]]_ {h_a }$\!.

\end{enumerate}
\end{itemize}
Here and in the remainder of this section
{\em we view each of our relations as
having length} $O(\log q)$ by using  \eqn{base 4} and  \eqn{Horner for
polynomials}.

   \begin{Proposition}
\label{SL2(q) Borel short bounded presentation}
If $ q\ne2, 3, 5, 9, $  then this produces a   bounded   presentation
    of length $O(\log q)$ for an infinite  central extension of
a Borel subgroup  $B$ of $\SL(2,q)$.
   \end{Proposition}

\Proof
Let $\hat B$ be the group defined by this   presentation.
Using the matrices written above it is easy to check  that the natural
map
$\pi\colon \hat B\to B$  sending $u,  h_\z, h_a, h_b,h_2 $ to
``themselves'' is a homomorphism.  In particular, (2) arises from the
fact  that
$h_b$ acts on
$\bigl( \begin{smallmatrix}
1 &  \star \\ 0& 1
\end{smallmatrix} \bigr)$ as multiplication of $ \star $  by $b^2=a^2+1$.

By (1), (2)  and
Lemma~\ref{BCRW1 trick} (with $h_a$ and $h_b$  in place of  $a$ and $b$),
$U:=\langle u^{\langle h_a
\rangle}
\rangle =
\langle u^{h_a^ih_b^j}\mid i,j\in \Z \rangle $
is abelian, and hence is  an elementary abelian $p$-group by (4).
In particular,  the definition of $[[u^{g(x) } ]]_{ h_ a}$ is     independent
of the  order of the  terms  for $g(x)\in \F_p[x]$.

Clearly $h_a$ acts as a linear transformation on $U$.  By (5), its
minimal polynomial  divides  the irreducible polynomial $m_{a^2}(x)  $.
Since the images of $u$ under powers of this linear transformation span
$U$, we have  $|U|\le q$, and then we obtain equality  using
$\pi$.

Now  
$\delta\mapsto [[ u^{f_{\delta;a^2}(x)} ]]_{h_ a}$ is an
isomorphism
$\F_q^+\to U$. Since $f_{\delta a^2;  a^2} (x)\equiv xf_{\delta ;
a^2}(x)$ (mod $m_{a^2}(x)) $  by the definition of $f_{\delta ;
a^2}(x)$ in Section~4.3.0, this   map  extends to  an
$\F_p[x]$-module isomorphism with $x$ acting on  $\F_q$ as multiplication
by $a^2$ and on $U$  as conjugation by
$h_a$.  (Note that there are choices for this isomorphism:
in place of $a^2$ we could use
any conjugate of it  by a field automorphism.)

By (6), the linear transformation induced by  $h_\z$ on $U$ has order
$(q-1)/(2,q-1)$.
   Thus, $h_\z$ acts on $U$ exactly as occurs in $B$.
All $h_\star$ do as well.  It follows from  (1) that $\hat B$ is a central
extension
 of
$B$, as required.

In view of \eqn{how polynomials act},
this presentation involves integer
exponents that can be much larger than
$\log p$.  However,   using  (3) these can be shortened as in \eqn{base
4}. Similarly, exponents that are polynomials are dealt with   in
\eqn{Horner for polynomials} using short words. Thus,  this
presentation  produces one that is short and bounded.

Finally, there is a semidirect product $U\semi \Z$ that satisfies the
above relations:  we are dealing with an infinite group.~\qedd
\smallskip

Note that
$$
Z(\hat B)=\langle h_\z^d, h_2 h_\z^{-n(2)},
h_ a h_\z^{-n( a)},h_b h_\z^{-n(b)}  \rangle ,
$$
where $d=|\z^2| =(q-1)/(2,q-1)$ and $\z^ { n(\star) } =\star$ for
$\star\in\{2,a,b  \}$.
Almost all of $Z(\hat B)$  will disappear when we deal with
simple groups.

As we saw above,  we can choose   $a=\z$ or $b=\z$ or
$ b= \z^2$; correspondingly  $h_ a=h_ \z$ or $h_ b=h_ \z$ or
$h_ b=h_ \z^2$.
This  allows us to obtain a shorter presentation by removing  one of
the generators
$h_a$ or
$h_b$,  and hence also three relations in (1).

The following is a simple   variation  on  this proposition.

\begin{lemma}
\label{AGL(1,q) short bounded presentation}
If $ q\ne2,  $  then the following produces a  bounded
presentation
    of length $O(\log q)$ for an infinite central extension of
the $1$-dimensional affine group $\AGL(1,q)=\{ x\mapsto \gamma x+\delta \mid
\gamma\in \F_q^*, \delta \in \F_q \}$$:$
\rm

\smallskip
{\noindent \bf Generators:} $u,  h'_2, h'_\zeta , h'_{\zeta +1}$.
(If $p=2$ discard $h_2'$.)
\smallskip

{\noindent\bf Relations:}
\begin{itemize}
\item []
\begin{enumerate}
\item  $[h'_\star, h'_{\bullet}] =1$ for $\star,\bullet \in \{2,\zeta
,\zeta +1\}$.
\item
    $u^{h _{\zeta+1}'} = u  u^{h_\z'{} } =   u^{h'_\zeta }u$.
\item $u^{h'_2}=u^2$.
\item $u^p=1$.
\item
$ [[  u^{m_\z(x)} ]  ] _{h'_\z } =1$  (cf. \eqn{how polynomials act}).

\end{enumerate}
\end{itemize}
   \end{lemma}

\Proof
This is similar to the preceding proof, using the same
$u$ but based on the new  diagonal matrices
   $h' _\star = \bigl( \begin{smallmatrix}
1 & 0 \\ 0& \star
\end{smallmatrix} \bigr)
   \mbox{ for $\star\in\{2,\z, \zeta+1\}$.}$~\qedd


\subsubsection{Unitary  groups}
\label{Nonabelian $p$-groups: unitary groups.}
\noindent   
We  now turn to  nonabelian versions of the same theme:  the unitary
groups
$\SU(3,q)$. This time let $F =\F_{q^2} $,  where $q=p^e$ 
and $F^*=\langle \z \rangle$
as usual.  We
will use  the following matrices that preserve the  form
$x_1\bar x_3+x_2\bar x_2+x_3\bar x_1
$:
\begin{equation}
\label{unitary matrices}
u(\alpha ,\beta)=\begin{pmatrix}
        1& \alpha & \beta \\
          0 & 1 & \!\!\!\! -\bar \alpha \\
   0 & 0 &1\\
     \end{pmatrix}\!\!, ~
h_2= \begin{pmatrix}
       1/ 2& 0 & 0\\
          0 & 1 &0 \\
          0 & 0 & 2\\
     \end{pmatrix} \!\! ,~
h_\delta= \begin{pmatrix}
        \bar \delta^{-1} & 0 & 0\\
          0 &  \bar \delta/\delta &0 \\
          0 & 0 &  \delta\\
     \end{pmatrix} \! \! , \,
\end{equation}
for any field elements such that
 $\beta+\bar \beta = - \alpha \bar \alpha
$ and $\delta\ne 0$.
Let  $U$ be the set of all matrices $u(\alpha ,\beta)$; its center has
$\alpha=0$.  Moreover,
\begin{equation}
\label{SU3 h action}
\mbox{$u(\alpha ,\beta)^{h_\z}=u(\z^{2q-1}\alpha ,\z^{ q+1} \beta )    $,
and  $u(\alpha ,\beta)^{h_ 2}=u( 2\alpha , { 4} \beta ) $ if $p\ne2$.  }
\end{equation}
Then  $B=U\rtimes \langle h_\z\rangle $ is a Borel subgroup of $\SU(3,q)$.
Its structure is as follows:
$|U|=q^3$, $Z(U)=U'=\Phi(U)$ is elementary abelian of order $q$, and
$U/Z(U)$ is elementary abelian of order $q^2$;  the
action of $h_\z$ on $U$ is given above, and it has order $(q^2-1)/d$
in that action, where
 $d=( 2q-1, q^2-1)=(3, q+1) $.

We  start with a property of $F$ somewhat analogous to \eqn{a and b}:

\begin{lemma}
\label{a,b}
\label{SU3 Borel field lemma}
If $q\not=3,5$ and   $q$ is odd$,$
then there are elements   $a,b\in F$ satisfying 
\begin{equation}
\label{a,b equation}
 a^{2q-1}+ b^{2q-1} =1~ \mbox{ and } ~
 a^{q+1}+ b^{q+1} =1 ,
\end{equation}
 such that
$F=\F_p[a^{2 q-1 }]$
and
$\F_q=\F_p[a^{q+1}]$.

If $q=2^e$ with $e\ge 2,$  then there are elements  $a,b\in F$
satisfying {\rm\eqn{a,b equation}} such that\/
$\F_q=\F_p[a^{q+1}]=\F_p[a^{2 q-1 }]$.
\end{lemma}

We give a proof  in  Appendix~\ref{Field lemma}.

\begin{Proposition}
\label{nonabelian Borel short bounded presentation}
\label{SU3 Borel}
There is an infinite  central extension of the above Borel subgroup $B\!$
{of\/ $\SU(3,q)$}  having  a   bounded   presentation
    of length $O(\log q)$.
   \end{Proposition}
\Proof 
Until the end of this proof we assume that   $q\ne 2, 3, 5$.
We use the abbreviations $\gamma'=\gamma^{q+1}$, $\gamma''=\gamma^{2q-1}$
for $\gamma\in F$.
  Section~4.3.0   contains the polynomial notation used
below.

Fix elements $a,b\in F$ as in the preceding lemma.
We base our presentation on
arbitrary elements  $u=u(\alpha_0 ,\beta_0)\notin Z(U)$ and $1\ne w\in
Z(U)$. 

We will use the following presentation with  (at most) 6 generators and 
20 relations.
\smallskip

 {\noindent \bf Generators:} $u,  w, h_2,   h_{a} , h_{b},
h_\z$.
(If $p=2$ discard $h_2$.)
\smallskip

{\noindent\bf Relations:}
\begin{itemize}
\item []
\begin{enumerate}
\item $[h_\star, h_{\bullet}] =1$ for $\star,\bullet \in \{2, a, b, \z
\}$.
\item $w= w^{h_a} w^{h_b} =w^{h_b} w^{h_a}$\!.
\item $w^{h_2}=w^{4}$.
\item $w^p=1$.
\item $[[w^{m_{a'}(x) } ]]_{h_a }=1$.
\item $[[w^{f_{\z';a'}(x) } ]]_{h_a }=w^{h_\z}$.
\item $u = u^{h_a}u^{h_b}w_1 = u^{h_b}u^{h_a} w_2 $.
\item
$[u,w] = [u^{h_a},w] =
1$.
\item $u^{h_2}=u^{2}w_3$.
\item $u^p=w_4$.
\item $[[u^{m_{a''}(x) } ]]_{h_a}=w_5$.
\item[(12$'$)]
$[[u^{f_{\z'';a''}(x) } ]]_{h_a }=u^{h_\z}w_6$  if $q$ is odd.
\item[(12$''$)]
$
([[u^{f_{\alpha ;a''}(x) } ]]_{h_a })^{h_\z}
[[u^{f_{\beta;a''}(x) } ]]_{h_a }
=u^{h_\z^2}w_6$
if $q$ is even and $\z''$ satisfies $\z''{}^2 = \alpha \z''{} +\beta $ for
$\alpha ,\beta \in \F_q$. 
\item[(13)]
$[u,u^{ h_\z}]= w _7 $  and  $[u^{ h_a},u^{ h_\z}]= w_8$  if $q$ is even.
\end{enumerate} 
\noindent
Here   $  w_1,\dots,w_8$ are  elements of
$W:=\langle w^{\langle  h_{a}\rangle}\rangle$;  they are  determined by
$\z$  and  our chosen matrices
$u$,  $w$   and  $h_\star$.
\end{itemize}

There is a surjection   $\pi\colon\hat B \to B$ from  the group $ \hat
B$ defined by this presentation onto $B$;  this is the significance of
the specific elements $w_i$.  Note that  (12$''$)
uses the fact that $\z'' $  is not in   $  \F_q$,  and hence
has a quadratic minimal polynomial over $\F_q$.

By (1), (2) and Lemma~\ref{BCRW1 trick},
$W $  is  abelian (using  $h_a,h_b$ in place of the elements $a,b$
in that lemma)  and
hence elementary abelian by (4).

Using $\pi$ we see that $|W|\ge q$.
Now (2), (5) and (6) let us   identify $W$ with the additive group of
$\F_q$ in such a way that each $h_\star$ acts as multiplication by 
a  conjugate of
$\star'=\star^{q+1}$.
 
By (1), (7) and (8), we can apply Lemma~\ref{BCRW1 trick nonabelian} 
to deduce that 
$[u,W]=1$ (using    $U_0=W $,  $W_0=1$  and $h_a, h_b$ in place of
the elements $a, b$ in that lemma).

If  $U:=\langle u^{\langle h_a\rangle },W \rangle $, then $[U,W]=1$  
since $W^{\langle h_ a, h_ b\rangle }=W $  by (2).

In view of (7),
  another application of   Lemma~\ref{BCRW1 trick} shows that
$U/W$ is abelian,
and by (10) it is elementary abelian.

Using $\pi$ we see that $|U/W|\ge |\F_p[a^{2q-1}]|$.
Now  (11) lets us identify $U/W$ with
$\F_p[a^{2q-1}]$ in such a way that each $h_\star$ acts as multiplication by
$\star''=\star^{2q-1}$.

  \smallskip
   {\noindent\bf Case $q$ odd.} 
Since $F=\F_q[\z'']$,    $(12')$ implies that
$h_\z$  acts irreducibly on $U/W$, and hence  $ |U/W| \le q^2$.
Thus, $|U|=q^3$.
Since $W=\Phi(U)$, the action of $h_\z$ on $U/W$ forces its action on $W$.
Since $U\unlhd \hat B$,
 by  (1)  it follows that $\hat B$ is a central extension  of $B$,
as required.

  \smallskip
   {\noindent\bf Case $q$ even.}
We   use Lemma~\ref{BCRW1 trick nonabelian} 
in order to bound $|\langle U,U^{h_\z}\rangle/W |$:
in that lemma we use $G=\langle h_a,h_b,U,U^{h_\z}\rangle/W$,    also
$uW$ and $u^{h_\z}W$ in place of  of $u$ and $w$, and finally $U_0 = W_0
=1$.  The  only condition that needs to be
checked is Lemma~\ref{BCRW1 trick nonabelian}(4), and this is just
$(13).$
Thus, $[u^{\langle h_a,h_b \rangle} , u^{h_\z\langle h_a,h_b \rangle}   ]
\le W$.  It follows from (1) that  
$ [U, u^{h_\z} ]\le W $ since  $U=\langle u^{\langle
h_a\rangle },W \rangle $.  

Then  $|\langle U,U^{h_\z}\rangle |\le q^3$  by (1), and hence equality
holds using
$\pi$.  By  (12$''$),
$ \langle U,U^{h_\z}\rangle ^{h_\z} = \langle U,U^{h_\z}\rangle $.
Now we can proceed as in the preceding  case, using $\langle U,U^{h_\z}\rangle $
in place of $U$.

\smallskip
Once again  $\hat B$  is infinite.   There are  at most  20 relations:
20   if $q$ is odd, and only 17 if $q$ is even since we delete 3 from (1)
together  with (3) and (9).

We still need to show that the  length is $O(\log  q )$.
By (3),
\eqn{base 4} and  \eqn{Horner for polynomials}, each element of $W$ has
length $O(\log q)$ in our generators.
Similarly,  by
(9),  (10),    \eqn{base 4} and  \eqn{Horner for polynomials},
every element in  $\langle U,U^{h_\z}\rangle/W$ 
 has length $O(\log q)$, and hence  
the  relations (11),  (12$'$) and (12$''$) can be viewed as  short.

This completes the proof of the proposition, but for future reference we
include   presentations for the omitted cases  $q = 2, 3,  5$.
When $q=2$, $B\cong  Q_8\times C_3$. For the remaining   $q $  we use
the     irreducibility of   $h_\z$  on $Q/Z(Q)$.
\smallskip


 {\noindent \bf Generators:} $u,  w, h_\z  $.
\smallskip

{\noindent\bf Relations:}
\begin{itemize}
\item []
\begin{enumerate}

\item $w^p=1$.
\item $[u,u^{h_\z} ]=w$.
\item $u^p=1$.
\item   $[u, w^{ h_{\z}^i}]=1$,  $i=0,1$.

\item $[[w^{m_{  \z ' }(x) } ]]_{h_ \z }=1$.
\item $[[u^{m_ { \z''} (x) } ]]_ {h_ \z }=w_1$.

 \end{enumerate} 
\noindent
Here,  $  w_1 \in W:=\langle w^{\langle  h_{ \z
}\rangle}\rangle$.
\end{itemize}
\smallskip
In the group  with this presentation,    $W $ is
cyclic of order   $q=p,$ and  by (5) it  is
normalized by $h_\z$.   By (4), $W$ commutes with    $U:=\langle
u^{\langle h_\z\rangle }\rangle $.  By (2), (3) and (6), $\langle U,W\rangle /W$ is elementary abelian,
and we can proceed as before.~\qedd


\subsubsection{Suzuki groups}
\label{Nonabelian $p$-groups: Suzuki groups.}


We start with analogues  of \eqn{unitary matrices} and \eqn{SU3 h action},
based on \cite[p.~133]{Suz}.  This time,
  $F=\F_q$, $q=2^{2k+1}$  and $F^*=\langle \z\rangle $. Let   $\theta\colon x \mapsto
x^{2^{k+1}}$  be the field automorphism such that
$\theta^2=2$.
Then  $B=U\rtimes \langle h_\z\rangle $, where
$U$ consists of all $(\alpha ,\beta)\in F^2$ with multiplication rule
\begin{equation}
\label{Suzuki rule}
(\alpha ,\beta) (\gamma,\delta)= (\alpha+\gamma,\beta+\delta+\alpha
\gamma^\theta),
 \end{equation}
and $W:=Z(U) =U'$ is all $(0,\beta)$.
Moreover, if  $\e\in \F^*$  let $h_\e$ denote  the automorphism of $U$
defined by
\begin{equation}
\label{Suzuki conjugation}
(\alpha ,\beta)^{h_\e} =( \e \alpha  , \e^{ \theta +1}  \beta).
\end{equation}
Then
$\langle  h_\z\rangle $ is transitive on the nontrivial elements of both
$Z(U) $ and
$U/Z(U)$.

If $x\in U$ and $v=(i,j,k,l )\in \Z^4$ 
   we write
\begin{equation}
\label{$Z^4$ action}
x^v=x^{(i,j,k,l )} = x^{h} \mbox{~where~} h={
h_{\z} ^i
h_{\z^\theta}^j
h_{\z+1} ^k
h_{\z^\theta+1} ^l
}.
\end{equation}

By \eqn{Suzuki conjugation},  
the following  trivial  identities in the field  $F$
$$  \hspace{6pt}
(\z) + 1 = (\z +1)
\qquad
(\z^\theta) + 1 = (\z^\theta +1)
$$
  translate into    relations modulo $W$:  
$$
u^{(0,0,0,0)} u^{(1,0,0,0)} \equiv  u^{(0,0,1,0)}
\qquad
u^{(0,0,0,0)} u^{(0,1,0,0)}  \equiv  u^{(0,0,0,1)}
$$
for any  $u\in U$. These relations correspond  to the following
collection $\mathcal{M_U} $ of subsets of $\Z^4$: 
$$
 \big\{ ^{^{^{}}}\!
\{(0,0,0,0),(1,0,0,0),(0,0,1,0)\},
\{(0,0,0,0),(0,1,0,0),(0,0,0,1)\}\big\}.
$$

Since $\theta^2=2$ we have  
$\lambda^{(\theta-1)(\theta  + 1) } =\lambda$
and 
$\lambda^{(2-\theta )(\theta  + 1) } =\lambda^\theta$ for any
$\lambda\in F$.  Then 
$$
(\z^{-1})^{\theta+1}(\z^\theta)^{\theta+1} + 1 =
((\z+1)^{-1})^{\theta+1}(\z^\theta+1)^{\theta+1}
$$
$$
(\z^2)^{\theta+1}((\z^\theta)^{-1})^{\theta+1} + 1 =
((\z+1)^2)^{\theta+1}((\z^\theta+1)^{-1})^{\theta+1},
$$
which, by \eqn{Suzuki conjugation},  imply that
$$
w^{(0,0,0,0)} w^{(-1,1,0,0)} = w^{(0,0,-1,1)}
~\mbox{ \rm and }~
w^{(0,0,0,0)} w^{(2,-1,0,0)} = w^{(0,0,2,-1)}
$$
for any  $w\in W$. These correspond  to  another collection
$\mathcal{M_W} $ of subsets of  $\Z^4$: 
$$
 \big\{ ^{^{^{}}}\!
\{(0,0,0,0),(-1,1,0,0),(0,0,-1,1)\}\!, \,
\{(0,0,0,0),(2,-1,0,0),(0,0,2,-1)\}\big\}.
$$
 
Our presentation of a Borel subgroup of $\Sz(q)$   
 uses Lemma~\ref{commutator and systems} to show that $U$ is
nilpotent of class 2.  This  relies heavily on  the observation that the 
collections $ \mathcal{M_U}  \cup \mathcal{M_{-U}},$  $ \mathcal{M_W} 
\cup \mathcal{M_{-W}}
$ and 
$\mathcal{M_U}  \cup
\mathcal{M_{-W}}$ are of finite type. This is almost trivial for the
first two  collections because they contain copies of $\M_{\pm \BCRW} $
up to a linear transformation.
The next lemma, whose proof is in Appendix~\ref{Suzuki triples},  gives
the same for the last collection.

\begin{Lemma}
\label{A Sz}
Let $\mathcal{M}_\Sz:=\mathcal{M_U} \cup \mathcal{M_{-W}}$
 consist of the  following
four triples of  vectors of $\Z^4$$:$
{\footnotesize
$$
\begin{array} {llll}
({\rm a})=\{(0,0,0,0),(1,0,0,0),(0,0,1,0)\}&&
({\rm b})=\{(0,0,0,0),(0,1,0,0),(0,0,0,1)\}
\\
({\rm c})=
\{(0,0,0,0),(1,-1,0,0),(0,0,1,-1)\}&&
({\rm d})=\{(0,0,0,0),(-2,1,0,0),(0,0,-2,1)\}.
\end{array}
$$
}
\vspace{-6pt}

\noindent
Then $\mathcal{M}_\Sz$ is of finite type and
there are exactly $139$ irreducible vectors for $\mathcal{M}_\Sz$,
with respect to the  length function
$|(x_1,x_2,x_3,x_4)|^2=x_1^2+2x_2^2+x_3^2+2x_4^4$  for
$ (x_1,x_2,x_3,x_4)\in \Z^4$.

Moreover$,$ there is a  set $S_0\subset \Z^4$  consisting of 
$16$ vectors such that$,$ if
$S\subseteq \Z^4$ contains $S_0$ and is
$\mathcal{M}_{\Sz}$-closed$,$ then    $S=\Z^4$.

\end{Lemma}

The proof of the first part of the lemma involves tedious computations
requiring  only high school  algebra  and bookkeeping.   More effort 
 is needed for the refinement in the second part (see 
Appendix~\ref{Suzuki triples}).


\begin{Proposition}
\label{Sz Borel}
There is an infinite  central extension of the above Borel
subgroup $B$ of
$\Sz(q)$  having  a   bounded   presentation
    of length $O(\log q)$.
   \end{Proposition}

\Proof
We base our presentation on an arbitrary  choice of $u\in U\setminus Z(U)$
and on $w=u^2 \in Z(U)$.
Let  $S_0$ denote the set of  vectors  appearing in Lemma~\ref{A Sz}. 
We will use the following 5 generators and $17+|S_0|=33$ relations:

\smallskip

{\bf Generators:} $u, h_\star $ for $\star \in \{\zeta, \zeta^\theta,
\zeta+1,\zeta^\theta+1\}$.

\smallskip
{\bf Relations:}
\begin{itemize}
\item []
\begin{enumerate}
\item $[h_\star, h_{\bullet}] =1$ for $\star,\bullet \in  \{\zeta,
\zeta^\theta,
\zeta+1,\zeta^\theta+1\}$.
\item $w^2=1$, where  $w:=u^2$.
\item $ w^{(0,0,-1,1)} =w w^{(-1,1,0,0)} $
in the notation of \eqn{$Z^4$ action}. 
%
\item $ w^{(0,0,2,-1)} =w w^{(2,-1,0,0)} $.
\item $[w,w^{(-1,1,0,0)} ]=1$. 
\item $[[w^{m_\z(x)}]]_{a}=1$
where $a:=h_{\z}^{-1} h_{\z^\theta}$. 
%
\item $w^{(1,0,0,0)} =[[w^{x^{2^{k+1}+1}}]]_{a}$,
where   $x^{2^{k}}$ is  reduced mod $ m_{\z }(x)$ in order 
to obtain  a short relation.
 
 \vspace{2pt} 
%
%
%

\item $u^{(1,0,0,0)} = u u^{(0,0,1,0)} w_1$. 
\item $u^{(0,1,0,0)} = u u^{(0,0,0,1)} w_2$. 
\item $[u^{(i,j,k,l)},w]=1$
for every $(i,j,k,l) \in S_0.$

\item $[u,u^{(1,0,0,0)}]=w_3$.
\item $[[u^{m_\z(x)}]]_{h_\zeta} = w_4$.
\item $u^{(0,1,0,0)} =[[u^{x^{2^{k+1}}}]]_{h_\zeta}
w_5$.
%
\end{enumerate} 
\noindent
Here, $w_1,\dots$ are suitable elements of  $W:=\langle  w^{\langle
a\rangle } \rangle $  determined by our choice of $u$ and $\z$.
\end{itemize}

\smallskip

First of all,   these relations are satisfied in  the Borel subgroup $B$: 
if 
$w=(0,\beta)$ then
$(0,\beta)^{h_{\zeta^\theta}
h_{\zeta}^{-1}} =(0,\z\beta)$  by \eqn{Suzuki conjugation}, which
yields (6) and (7); \eqn{Suzuki conjugation} also yields  
  (12)  and
(13).  We have already seen  (3), (4), (8) and (9), while (1), (2), (5),
(10)  and (11) are clear.%

By (1), (3),  (5) and Lemma~\ref{BCRW1 trick} (with $a$ as above and
$b:=h_{\z+1}^{-1}h_{\z^\theta+1}$),
$W$ is abelian. By (2), (3)  
  and  (6), $ W$ can be identified with $F$ in
such a way that
$a$ and
$b$ 
 act as multiplication by
$\zeta$ and $\zeta+1$, respectively.
By (7), $ h_\z$ acts as
multiplication by  $\z^{\theta +1}$, 
which together with (3) and (4)  implies that $h_\star$  acts on $W$ as
$\star^{\theta +1}$ for   $\star \in \{\zeta, \zeta^\theta,
\zeta+1,\zeta^\theta+1\}$.

Once again we use the notation 
\eqn{$Z^4$ action}. 
We  proceed essentially as in Section~\ref{BCRW trick }.   Let
$$
S:=\{ v\in \Z^4~\big|~  [ u^v,w]=1\}\!.
$$
Let $\mathcal{M}_\Sz$ be as in Lemma~\ref{A Sz}. 

By Lemma~\ref{commutator and systems} (with $U_0=W$ and  $W_0=1$), 
 combined with (3), (4), (7), (8)
and the fact that
$[w,W]=1$, the set   $S$ is
$\mathcal{M}_\Sz$-closed (Definition~\ref{closed set}(1)). 
Since $S$ contains $S_0$ by (10), we have $S=\Z^4$ by
Lemma~\ref{A Sz}, so that $[u^{(i,j,k,l)},w]    
=1$ for all
$(i,j,k,l)\in
\Z^4$. Thus, each $[u,w^{(i,j,k,l)}]=1$ and hence $[u,W]=1$.

It follows that  $W=W^{h_\z}$ commutes with $U:=\langle u^{\langle h_\z\rangle }
\rangle$. By (1), (8), (9) and (13), $U$ is invariant under each
$h_\star$, and hence   contains 
$\langle (u^2)^{\langle a\rangle
}\rangle =\langle w^{\langle a\rangle }\rangle =W$. 
  Thus, $W\le Z(U)$.

By (1), (8) and (11), we can apply Lemma~\ref{BCRW1 trick} to
$U/W.$ Then this group is elementary abelian  since $w=u^2$.
%
By   (12),  $|U|\le q^2$.
As in Proposition~\ref{nonabelian Borel short bounded presentation},
it follows that  the presented group $\hat B$ is a  central extension of
$B$.  Once again  $\hat B$  is
infinite.~\qedd


\subsubsection{Ree  groups}
\label{Nonabelian $p$-groups: Ree groups.}
\label{Ree  groups}
We cannot prove Theorem~\ref{A} for the Ree groups $^2G_2(3^{2k+1})$, but
we can at least handle their Borel subgroups.  Let $F=\F_{3^{2k+1}}$, 
$\theta\colon x\mapsto x^{3^{k+1}}$ and
$F^*=\langle \z\rangle $.
 By \cite{KLM},  $B=U\rtimes \langle h_\z\rangle $ where
   $U  = \{x(a,b,c)\mid a,b,c \in F \} $ and 
\begin {align}
\label{Ree multiplication}
  x (a_1,b_1,c_1) x (a_2,b_2,c_2) &=
x (a_1 + a_2,b_1 + b_2 +
a_1 a_2^{\theta},  
\\ \nonumber
&\quad \quad \quad \quad \quad \quad  c_1 + c_2 -   
a_1^2 a_2^{\theta } + a_2 b_1 + a_1 a_2^{\theta + 1}) 
\\
\label{Ree conjugation}
   x(a,b,c)^{h_\z} &=  x(\z^{2-\theta}a,\z^{\theta-1}b,\z c)  .
\end {align}

\begin{Proposition}
\label{Ree Borel}
There is an infinite  central extension of a Borel
subgroup $B$ of
$^2G_2(q)$  having  a   bounded   presentation
      of length $O(\log q)$.
     \end{Proposition}

A sketch of a proof is given in Appendix~\ref{Ree triples}. 
These groups are harder to work with than previous ones.  In particular,
whereas we used $\Z^2$ for $\SL(2,q)$, and $\Z^4$ for Suzuki
groups, for these Ree groups we  use $\Z^6$.

\subsection {Presentations for rank 1 groups}
\label {Presentations for rank 1 groups}
\label {A presentation for PSL(2,q)}
\label{new PSL(2,q)}

We are now ready to consider the rank~1 groups $\SL(2,q), $  $
\PSL(2,q), $  $\PGL(2,q), $  $ \PSU(3,q) $,  $ \SU(3,q) $  and  $\Sz(q)$.
Except for headaches due to field extensions or nonabelian $p$-groups,
the idea in this section is the same as the one in
Section~\ref{PSL(2,p)}

In each case there is a
Borel subgroup $B=U\semi \langle h \rangle $, with   $U$ a
$p$-group.  There is an involution  $t$ (mod $Z(G)$ in the case
$\SL(2,q)$ with $q$ odd) such that $h^t=h^{-1} $
(or $h^{-q }$ in the unitary case).  These lead us to the {\em Steinberg
presentation}  for these groups \cite[Sec.~4]{St2}:
\begin{theorem}
\label{Steinberg presentation rank 1}
{\rm(Steinberg presentation)} 
Each of the above groups  has a presentation using
\begin{enumerate}
\item any presentation for $B,$
\item any presentation for  $ \langle h ,t \rangle ,$ and
\item $|U|-1$ relations of the
form
\begin{equation}
\label{Steinberg relation}
u_0^t=u_1h_0tu_2,
\end{equation}
with $u_0,u_1, u_2$ nontrivial
elements of
$U$ and
$h_0\in \langle h \rangle $ $($one such relation for each choice of
$u_0$$).$
\end{enumerate}
\end{theorem}
We have already
used a special case of this presentation in Section~\ref{PSL(2,p)}.
We now use this presentation   in order to prove  the following


\begin{Theorem}
\label{rank 1}
\label{PSL2q length}
\begin{itemize}
\item[(a)] 
All of the groups $\SL(2,q),$   $\PSL(2,q),\hspace{-1pt} $ $ \PGL(2,q),\hspace{-1pt} $  
$\PSU(3,q),$  
$\SU(3,q)   $ and
 $\Sz(q)$    have bounded presentations of length $O(\log q)$.
\item[(b)]
Each  element of each of the preceding groups  can be written as a word
of length
$O(\log q)$ in the generating set used in the above presentation.
\end{itemize}
\end{Theorem}
\Proof
 We will treat the groups separately.
In each case we will show that we only need some of the relations  
\eqn{Steinberg relation}
 in order to deduce all of those relations.


\subsubsection{Special linear groups~}
\label{Case $SL(2,q)$.}
  Here we  will be more explicit than in
the  unitary or    Suzuki cases.
 We may assume that $q>9$.
We will use the notation leading up
to Proposition~\ref{SL2(q) Borel short bounded presentation}.  
We will show that   $\SL(2,q)$ has the   following presentation using (at
most) 6 generators and 22 relations.

\smallskip

{\noindent
\bf Generators:}
$u, t, h_2, h_a, h_b, h_\zeta $.  (If $p=2$ discard $h_2$.)
\smallskip


{\noindent\bf Relations:}
\begin{itemize}
\item []
\begin{enumerate}
\item $[h_\star, h_{\bullet}] =1$ for $\star,\bullet \in \{2,a,b,\zeta \}$.

\item    $u^{h_b} = u  u^{h_a}=   u^{h_a}u$. 
\item $u^{h_2}=u^4$.
\item $u^p=1$.
\item $[[ u^{m_{a^2} (x) } ]]_{h_ a} =1$.
\item $u^{h_\zeta }= [[ u^{f_{\zeta ^2; a^2} (x)} ]]_{h_ a }$\!. 
\item $h_\star^t=h_\star^{-1}$ for $\star \in \{2,a,b,\zeta \}$.
\item $[t^2,u] =1$.
\item  $t= u  u^{t} u   $.
\item $h_\star  t = [[ u^{f_{\bullet; a^2}  (x)  } ]]_{h_ a} \cdot      [[
u^{f_ {\star; a^2}(x) }]]_{h_ a}^{\hspace{1pt} t}  \cdot [[u^{f_{\bullet;
a^2}   (x)} ]]_{h_ a}    $, where
$\bullet =  \star^{-1}$ for
$\star
\in \{2,a,b,\zeta \}$.
\end{enumerate}
\end{itemize}
\smallskip
 
This   certainly involves   more  relations  than
one might expect (\cite{CRW1} uses at most 13 relations). In
\cite{GKKL2}  we give a similar presentation (of bit-length
$O(\log q)$ but  {\em not} of length
$O(\log q)$)
 using only 9 
relations.

 In order to show that the group $G$ presented in this manner
is isomorphic to $\SL(2,q)$,  note that $\SL(2,q)$ is a homomorphic
image of $G$ such that 
$t$ is in \eqn{generators} and 
$\hat B:=\langle u,h_2,h_{a},h_{b},h_\zeta \rangle $  
(which we identify with a subgroup of $G$)  maps onto
a Borel subgroup  $B$ of $\SL(2,q)$.
By (1)-(6),  we can apply Proposition~\ref{SL2(q) Borel short
bounded presentation}:  $\hat B$ is isomorphic to a central extension
of $B$.

In particular,  $U:=\langle u^{\langle h_{a}\rangle }\rangle $ is a normal subgroup of
$\hat B$,  and  $h_{a}$ acts  irreducibly on
$U $. Then $U=[U,h_a]\le \hat B'\le G'$, so that
$ \langle U,U^t\rangle \le G'$.
By  (9), $t \in \langle U,U^t\rangle $, and then all $h_\star\in \langle
U,U^t\rangle $ by (10), so that
$G=\langle U,U^t\rangle $  and  $G$ is a perfect group.

     In  the identification of $U$ with $\F_q$
in Proposition~\ref{SL2(q) Borel short bounded presentation},
   $h_\zeta $ acts as multiplication by $\zeta ^2$.
By (7), since $h_\z^{(q-1)/(2,q-1)}$ centralizes $U$ it also
 centralizes  $U^t$ and hence $G$.
Similarly,    by (8) we have
$Z:=\langle h_\z^{(q-1)/(2,q-1)},t^2\rangle  \le Z(G)$.

We claim that  $G/Z $ satisfies the presentation for
$\PSL(2,q)$  in Theorem~\ref{Steinberg presentation rank 1}.  For,
$\langle h_\z,t\rangle  / Z $ is dihedral of order $2(q-1)/(2,q-1)$ by
(7). Also,  $U\langle h_\z,t^2\rangle  / Z$ is isomorphic  to a Borel
subgroup of
$\PSL(2,q)$.
Moreover,   $\langle h_\zeta \rangle $ acts on the
nontrivial elements of $U$ with
    at most   two orbits;  orbit-representatives are
   $u^1$ and  $[[ u^{f_ {\zeta; a^2} (x) }]]_{h_a}$. 
By    (9) and (10),  relation
      \eqn{Steinberg relation} holds for $u_0=u^1$ or 
$[[ u^{f_ {\zeta; a^2} (x)}]]_{h_a}$; 
   conjugating
      these relations  by $\langle h_\z\rangle $ produces all of the required  relations
      \eqn{Steinberg relation}.
This proves our claim.

The only perfect central extensions of $\PSL(2,q)$  are
$\PSL(2,q)$  and $\SL(2,q)$ (since $q> 9$).  Hence, $G\cong
\SL(2,q)$.

This presentation is clearly bounded.  It can be made short:
relations (5), (6) and (10) have length $O(\log q)$ by \eqn{Horner for
polynomials} and
\eqn{base 4}. 
(Note that the generators $h_a$ and $h_b$ are not needed here.)

We  still need to  verify \   Theorem~\ref{rank 1}(b).
Every element of $\SL(2,q)$ is in $UU^t UU^t U$.  Thus, we only need to
verify (b) for $U$,  and again this  follows from   \eqn{Horner for
polynomials} and \eqn{base 4}.

\smallskip

{\noindent \em Case $\PSL(2,q)$.~}
Replace the previous  relation (8)  by  $ t^2  =1$.

\smallskip

{\noindent \em Case $\PGL(2,q)$.}  This is proved as in the preceding
case, using  Lemma~\ref{AGL(1,q) short bounded presentation}
in place of
Proposition~\ref{SL2(q) Borel short bounded presentation}.

\subsubsection{Unitary groups}
\label{unitary case}
 We may assume that $q>5$.
We use the  matrices in \eqn{unitary matrices}, together with the
additional matrix
\begin{equation}
\label{SU3 matrix t}
t=\begin{pmatrix}
        0& \, \, 0 & 1 \\
          0 & \!\!  -1   & 0\\
   1 &\, \,  0 & 0\\
     \end{pmatrix}.
\end{equation}
%
In the notation of \eqn{unitary matrices}, 
if $\beta \ne0$ and $\beta+\bar \beta+\alpha \bar \alpha =0$ then
\eqn{Steinberg relation} becomes 
\begin{equation}
\label{Steinberg unitary relation}
u(\alpha ,\beta)^t=
u(-\alpha/\bar\beta,1/\beta)\,h_ \beta \,t\,
u(-\alpha/ \beta, 1/\beta).
\end{equation}

We will make fundamental use of a result of Hulpke and Seress \cite{HS}, 
who showed that the presentation in Theorem~\ref{Steinberg presentation
rank 1} can be shortened by using (1), (2) and a carefully chosen set of
at most 7 of the relations (3).

The following presentation uses (at most) 14    generators and  52
relations.  We will show that  it defines a group $G$ that is   
isomorphic to
$\SU(3,q)$.
\smallskip

{\noindent
\bf Generators:}
$u, w, t, h_2,   h_{a} , h_{b},
h_\z , \breve h_i$ for
$1\le i\le 7$.  (If $p=2$ discard $h_2$.)
\smallskip

{\noindent\bf Relations:}
\begin{itemize}
\item []
\begin{itemize}
\item [(0)] All relations (1)-(13) used in Proposition~\ref{nonabelian
Borel short bounded presentation}.
\item [(14)] $ t^2 = 1$.

\item [(15)] $h_\star  =  v_{\star 1} v_{\star 2} ^t v_{\star 3}  t $
for $\star\in \{1,2,a, b, \z  \}$, where $h_1 : =1$.

\item[(16)] $u^{h_\star^t} = [[u^{f_\star(x)}]]_{h_\z}$ for $\star\in
\{2,a,b,\z\}$. 

\item [(17)] $ u_i^t=   u_{i1} \breve  h_{i}   t   u_{i2} $ for
$1\le i\le 7$,   relations due to Hulpke and Seress~\cite{HS}.
 
\vspace{2pt}
\item[(18)] 
$[h_\z,  h_\z ^t]= 1$.

\item[(19)] 
$ [h_\z, \breve h_i]=   1$ for $1\leq
i\leq 7$.

\vspace{2pt}
\item[(20)] 
$u^{\breve h_i} = [[u^{g_i(x)}]]_{h_\z}$ and  
$u^{\breve h_i^t} = [[u^{ g_i^\bullet (x)}]]_{h_\z}$
 for $1\leq i\leq 7.$

\vspace{2pt}


\end{itemize}
\noindent
Here $u_{ i},u_{ i1} ,u_{ i2},   v_{ \star 1} ,v_{ \star 2}   ,v_{
\star 3}$ are specific elements of $U:=\langle u^{\langle h_\z \rangle}
\rangle$, and
 $f_\star (x)$,  $ g_i (x)$,  $ g_i^\bullet (x)$  are specific
polynomials in $F[x]$.
 
\end{itemize}


There is a surjection
$\pi\colon G\to\SU(3,q)$ taking each of the generators  to ``itself''.
This is clear for (0) and (14), while (15) comes from  \eqn{Steinberg
unitary relation}.     Relations (16)   and (18)-(20) encode  the fact
that,  in $\SU(3,q)$,  $h_\star$ and $\breve h_i $    are powers of
$ h_\z$,
while $h_\star^t=h_\star^{-q}$ and $\breve h_i^t=\breve h_i^{-q}$.  
 There are
at most  7 relations of the form (17) in the presentation for $\SU(3,q)$
given in
\cite{HS}.
 (All of
these relations are given far more explicitly in \cite[(11)~and~(12)~on
p.~103]{BGKLP} and \cite{HS}, using   \eqn{Steinberg unitary relation}.)

Moreover, $\pi$  maps  
$\hat B :  =\langle u, w, h_2,h_ a,h_b,h_\z \rangle $
(which we identify with a subgroup of $G$) 
onto a Borel subgroup $B$  of $\SU(3,q)$.
Proposition~\ref{nonabelian Borel short bounded
presentation} implies that $\hat B$ is   a
central extension of $B$.
In particular,  since $q>5$, $U$ is a
normal  subgroup of
$\hat B$ of order $q^3$   and  $h_{\z}$ acts  irreducibly on
its  Frattini quotient
 $U/Z(U) $.  It follows that $U=[U,h_\z]\le G'$.
By  (15) for $h_1$, we have  $t \in \langle U,U^t\rangle $,  and 
then all $h_\star, \breve h_i \in \langle U,U^t\rangle $  by (15) 
and (17), so that
$G=\langle U,U^t\rangle $  and  $G$ is a perfect group.

In particular, in  the identification of $U/Z(U)$ with $\F_{q^2}$
in Proposition~\ref{nonabelian Borel short bounded
presentation},
   $h_\zeta $ acts as multiplication by $\zeta ^{2q-1}$ (cf. \eqn{SU3 h
action}), and  hence has order $d:=(q^2-1)/(3,q+1)$ in  its action on 
$U$.

By (16)  and  (18), $h_\z^t$ also acts on $U,$ and  is irreducible on
$U/Z(U)$.  Thus,  
$U=\langle u^{\langle
h_\z^t\rangle}\rangle $.

Let $Z= Z(G)$. 
We claim that  $h\in \langle h_\z\rangle Z$
if $h$ is any    of the elements  $h_\star,
\breve h_i$. 
By Proposition~\ref{nonabelian Borel short bounded
presentation},   (19)  and (20),  
 $h$ normalizes $U$  and 
there is an integer $k$ such that 
$1\le k < d$ and $hh_\z^k$ centralizes $U$.  By 
(16),  (19) and (20), since  $h^t$ commutes with $h_\z ^t$ it also
normalizes 
$U=\langle u^{\langle
h_\z^t\rangle}\rangle $, and  there is also an
integer
$k'$ such that
$1\le k' < d$ and 
$(hh_\z^{k'})^t$ centralizes $U$.  
 Since $\pi(G)=\SU(3,q)$,
we have   $k=k' $, so that   $hh_\z^k$ centralizes 
$\langle U,U ^t\rangle  =G$, which proves our claim.

We can now use Theorem~\ref{Steinberg presentation rank 1} to show that
$G/ Z  \cong \PSU(3,q)$.
For, (16)  implies that
$\langle h_\z,t, Z\rangle  / Z $ has order $2(q^2-1)/(3,q+1)$, and hence
is isomorphic to the required subgroup of $\PSU(3,q)$. Also, 
$\hat B Z/Z=U\langle \mbox{all $h_\star $}   
\rangle Z / Z=U\langle h_\z
\rangle Z / Z$ 
 is isomorphic  to a Borel subgroup of
$\PSU(3,q)$.
The only relations remaining to be checked for
$G/ Z $ are those in \eqn{Steinberg relation}, and these  are    (17)  
(read modulo $ Z $).

Since  $\SU(3,q)$ is an
epimorphic image of $G $,  it follows  that
$ G \cong \SU(3,q)$  (see \cite[Corollary~6.1.7]{GLS}).
This presentation is clearly bounded.  It can be made  short using 
\eqn{Horner for polynomials} and
\eqn{base 4}.  

Finally,   Theorem~\ref{rank 1}(b) holds  when $g\in U$ by
 \eqn{Horner for polynomials}
and \eqn{base 4},  then also when  $g\in U^t$, and hence whenever $g\in
UU^t UU^t U  =\SU(3,q)$.
\qedd
\smallskip

  By  Theorem~\ref{rank 1}(b), we  also obtain a short
bounded presentation for $\PSU(3,q)$  since a generator  of
$Z(\SU(3,q))$ has length  $O(\log q)$.

\subsubsection{Suzuki groups}
\label{Suzuki case}
 We may assume that
$q>8$.  
This time we will use a remarkable result  of Suzuki \cite[p.~128]{Suz}, 
who showed that the presentation in Theorem~\ref{Steinberg presentation
rank 1} can be shorted by using (1), (2) and a single carefully chosen 
relation  (3) (in fact, this is just the relation needed to define
$\Sz(2)$).

We use the following presentation having 7 generators and 43 relations.  

\smallskip
{\noindent
\bf Generators:}
$u,w, t, h_\z, h_{\zeta+1 }, h_{\z^\theta}, h_{\z^\theta+1}$.
\smallskip

{\noindent\bf Relations:}
\begin{itemize}
\item []
\begin{itemize}
\item [(0)] All 33 relations  
(1)-(13) in Proposition~\ref{Sz Borel}.
\item [(14)] $t^2=1$. 
\item [(15)] $h_\star^t = h_\star^{-1}$ for $\star\in \{\z,\zeta+1 ,
\z^\theta,\zeta^\theta+1   \}$.
\item [(16)] $ t=u_1u_2^{\hspace{1pt} t} u_3$, a relation due to Suzuki
\cite[(13)]{Suz}.
\item [(17)] $h_\star t = u_{\star 1}u_{\star 2}^{\hspace{1pt} t} u_{\star
3}$ for $\star\in \{\z,\zeta+1 ,
\z^\theta,\zeta^\theta+1   \}$.
\end{itemize}
\noindent
Here $u_{ 1},u_{ 2} ,u_{ 3},u_{\star 1},u_{\star 2},
u_{\star 3}$
are specific elements of $U:=\langle u^{\langle h_{\zeta }\rangle }\rangle $ that depend on
the choices of $\z$ and $u$  in Section~\ref {Nonabelian $p$-groups:
Suzuki groups.}.
\end{itemize}

\smallskip

Note that, in the Suzuki group case,  we do not   
need  additional relations as in   
Section~\ref{unitary case}(19), since the power of $h_\z$ occurring in
relation  (16) is~1.

It order to show that the group $G$ presented in this manner
is isomorphic to $\Sz(q)$, first note that  $\Sz(q)$ is a
homomorphic image of $G$ such that $\hat B:=\langle u,h_\zeta ,h_{\zeta+1},
h_{\z^\theta}, h_{\z^\theta+1}   \rangle $ projects onto a Borel subgroup
$B$  of
$\Sz(q)$. (For  more explicit versions of (16) and (17), see
 \cite[(13)    and (38)]{Suz}.)

The rest of the proof is essentially the same as in 
Sections~\ref{Case $SL(2,q)$.} and \ref{unitary case}.  
 Both $U\langle h_\z \rangle/\langle h_\z^{q-1}\rangle $ and
$\langle h_\z,t\rangle /\langle h_\z^{q-1}\rangle $ are as required in
Theorem~\ref{Steinberg presentation rank 1}. Suzuki \cite[p. 128]{Suz} 
proved   that a single  relation (16) is the  {\em
only} one of the  relations
\eqn{Steinberg relation}  needed  in order to deduce  all  $q^2-1$  of
them.
Hence,
$G/\langle h_\z^{q-1}\rangle   \cong
\Sz(q)$.  Finally,
$G\cong \Sz(q)$ since $G$ has no proper perfect central extensions
(recall that  $q>8$; cf. \cite[pp. 312-313]{GLS}).

This presentation  
can be made  short 
 by using  \eqn{Horner for polynomials} and
\eqn{base 4}.  
It remains to  consider Theorem~\ref{rank 1}(b).  Each element of $U$ or
$U^t$   is a word of length $O(\log q)$ in our generators.
Now (b) follows from the fact that $\Sz(q)=UU^t UU^t U$.

\emph{This completes the proof of} Theorem~\ref{rank 1}.  \qedd

\medskip
We have not dealt with the  Ree groups in Theorem~\ref{A}.  The obstacle
to both short and bounded presentations of these groups is the fact that
all presentations presently known
 use more  than $q$ relations of the form \eqn{Steinberg relation}.



   \section{Fixed rank}
\label{Fixed rank} 
In this section we will prove the following result, which implies
Theorem~\ref{A} when the rank is bounded:

\begin{Theorem}
\label{Short presentations}
All perfect central extensions of   finite  simple groups  of Lie
type and given rank
$n,$  except perhaps  the Ree groups
$^2G_2(q),$  have  presentations
of length $O(\log|G|)$ with   $O(n^2)$ relations.
\end{Theorem}

In view of
   Theorem~\ref{rank 1}, we may assume that the   $n\ge2$. The idea
in the proof  is simple:  replace  the  relations in the
Curtis-Steinberg-Tits presentation
(Sections~\ref{Curtis-Steinberg-Tits presentation  } and
\ref{Curtis-Steinberg-Tits presentation}) involving rank 1 subgroups with the
short  bounded presentation in    Theorem~\ref{rank 1}.

\subsection{Curtis-Steinberg-Tits presentation}
\label{Curtis-Steinberg-Tits presentation  }
Steinberg \cite{St1} gave a presentation for the groups of Lie type
involving all roots groups corresponding to the  root
system.  Curtis \cite{Cu} and Tits \cite[Theorem~13.32]{Tits_book}
independently obtained   similar presentations involving only rank 1 and
2 subsystems of the root system.  It is the latter type of presentation
that we will use, since less data is required.  There are two slightly
different versions of the presentation: the ``high-level" one,  due to
Tits  and   based
on buildings, that does not explicitly deal with root groups, and the
more involved one that uses root groups and their commutator relations.
For now we will state a version of  the higher level one, explaining the
more familiar one with formulas in the next section.

Each group $G$ of Lie type  has a suitably defined root system
$\Phi$ and fundamental system $\Pi=\{\alpha_1,\dots,\alpha_n\}$ of
  roots.   We view the twisted  groups of  type
 ${}^2{\!} A_m(q)$ (with $m \ge 3$ and  $m=2n$ or $m=2n+1 $),  ${}^2{\!} D_m(q)$
(with $m \geq 4$),
   ${}^2 {\!}E_6(q)$  or ${}^3 {\!}D_4(q)$  as having   root systems of  type
  $C_{n}$,
$B_{m-1}$,  $ F_4$  or $G_2$, respectively. There is also a suitable
sense in which
$^2{\!}F_4(q)$ has 16 roots  \cite[Corollary~2.4.6]{GLS}.

\begin{Theorem} {\rm(Curtis-Steinberg-Tits presentation
\cite{Cu,St1,Tits_book} 
and \cite[Theorem~2.9.3]{GLS}) }
\label{Curtis-Steinberg-Tits presentation}
Let $G$ be a simple group of Lie type and rank at least $ 3$.
Let $\Pi$   be a fundamental system of roots in  the root system.
Whenever $\alpha \in \Pi$ let $L_{\alpha }$  be the corresponding rank $1$
subgroup of $G; $  whenever $ \alpha ,\beta \in\Pi,$ $  \beta
\ne\pm\alpha    ,$   let
$L_{\alpha,\beta}:=\langle L_\alpha ,L_\beta \rangle $ be the
corresponding rank $2 $  subgroup of $G$. 
Then   
 the amalgamated product of the rank $2 $  subgroups
$L_{\alpha,\beta}  ,$  amalgamated along
the rank $1$ subgroups
$L_{\alpha},$ is a central extension of $G$.
\end{Theorem}
 
  Of course, the above version of the
Curtis-Steinberg-Tits presentation  also   gives no information concerning
presentations of rank 2 groups.
For these we need to use  commutator
relations.  The commutator
version  of the Curtis-Steinberg-Tits presentation is in \cite{Cu}, as
well as  \cite[Theorem~2.9.3]{GLS}, and will be described in the next
section.
 
\subsection{Commutator relations}  \label{Commutator relations}

Bounded rank groups were considered in \cite{BGKLP}.  We will use the method of
that paper, inserting our presentations for rank 1 groups  in place  of the
standard ones used there.  In that paper, all commutator  relations are
correctly presented, as are the definitions of the toral elements
$h_\beta $.
However, there were significant errors in that
paper:  for most of the twisted groups,  the actions of
 $h_\beta $ on the root group  elements $x_\alpha(a)$  (or $x_\alpha(a,b)$
 for  odd-dimensional unitary groups) were  not correct.  The
methodology used  was, however, correct, and will be used  here without
any explicit formulas. 
The same methodology was used in \cite{KoLu}, though not for all possible
types of group and without considerations of lengths of presentations.
 
We assume that $G$ is a universal group of Lie type \cite[p.~38]{GLS}. 
We will factor out 
part or all of $Z(G)$ later in this section; this
point of view occurred in the Section~\ref{A presentation for PSL(2,q)},
and will also be used in Section~\ref{Theorems A and B}.
We will ignore the   tiny number  
of   cases in which $G$ has a perfect central extension by a group whose 
order is the characteristic
\cite[p.~313]{GLS}.

With each root $\alpha \in \Phi$
there is associated a rank 1 group
$L_\alpha \cong\SL(2,q), $   $\PSL(2,q),$  $\SL(2,q^2), $ $ \PSL(2,q^2), $
$\SL(2,q^3),$   $ \SU(3,q)$ or
$\Sz(q)$.  As indicated above, we use the presentation ${\bf P}_\alpha $
for
$ L_\alpha
$ obtained in Section~\ref{Presentations for rank 1 groups}. 
Among our generators of $L_\alpha $ we have the following:   a torus
$\langle h_\alpha \rangle $;  a $p$-element
$u_\alpha $ lying in  a root group
$U_\alpha  $  of
$L_\alpha $ and of $G$  (note that $U_\alpha =\langle
u_\alpha^{\langle h_\alpha 
\rangle}\rangle$  except when $L_\alpha=\SU(3,2)$); and a ``reflection''
$r_\alpha $ (called ``$t$'' in Section~\ref{Presentations for rank 1
groups}) that is an involution modulo
$Z(L_\alpha)$ such that
$r_\alpha $ normalizes $\langle h_\alpha  \rangle$ and
$U_\alpha^{r_\alpha}=U_{-\alpha}$ is the opposite root group (relative to
$\Phi$).   When $q$ is tiny the generic presentations in  
Section~\ref{Presentations for rank 1 groups} are not applicable, but
there certainly are short, bounded presentations using at least  the above
generators.

We have  $L_{-\alpha} = L_\alpha $.
Since $G$ is universal, so is $L_\alpha$.  
 Checking all cases, we find
that 
\begin{equation}
\label{h order}
\mbox{$|h_\alpha|= |\bar U_\alpha^\#|$,     where  
  $\bar U_\alpha^\#:=\bar
U_\alpha \backslash \{1\}$ for $\bar
U_\alpha:=U_\alpha/\Phi(U_\alpha)$}
\end{equation}
(here $\Phi(U_\alpha)$ denotes the Frattini subgroup of  $U_\alpha $).

There is a Borel subgroup $B$ associated with the positive roots in
$\Phi$.  It contains $U_\alpha$ for all positive $ \alpha$ and has  
a maximal torus  $H:=\langle h_\alpha \mid \alpha \in
\Pi\rangle $.  Moreover,  $H$  normalizes each root group $U_\alpha$,
and 
 acts on  $\bar U_\alpha $
  as multiplications  by   elements
 in an extension field
 $\F_{q^{e(\alpha)}}$ of $\F_q$,
  where $|\bar U_\alpha|=q^{e(\alpha)}$,
 $e(\alpha)\le 3 ,$ and 3 occurs only for
 $^3{\!}D_4(q)$.  



By Theorem~\ref{Curtis-Steinberg-Tits presentation} 
we need to consider all of  the subgroups 
$L_{\alpha_i,\alpha_j}$, $i\ne j$,  arising from our  
fundamental system $\Pi=\{\alpha_1,\dots,\alpha_n\}$. 
Here $L_{\alpha_i,\alpha_j}$ is  a rank 2 group whose  root system
 $\Phi_{\alpha_i,\alpha_j}$ 
has the  fundamental system
$  
\{\alpha_i,\alpha_j\}$.
Crucial for our argument (and those in \cite{BGKLP,KoLu}) is the fact
that
$H$ has boundedly many orbits on the direct product of suitable pairs
 $\bar U_\alpha$, $\bar U_\beta$.  For this we need  the following%

\begin{lemma}
\label{orbits on pairs of root groups}
\label{H-orbits}
Assume that  
$L:=L_{\alpha_i,\alpha_j}$  is  not $
\SU(5,2).$  
 Let $\alpha$ and $ \beta \ne
\pm\alpha$ be in  $\Phi_{\alpha_i,\alpha_j}$.   If
$u_\alpha \in U_\alpha \backslash \Phi(U_\alpha)$ and  $u_\beta \in
U_\beta \backslash \Phi(U_\beta),$ then 
$C_{\alpha ,\beta}:=C_{\langle h_\alpha ,h_\beta \rangle}(\langle
u_\alpha ,u_\beta \rangle )$ is in the center of
 the rank $2$ subgroup $X_{\alpha ,\beta}:=\langle U_\alpha ,U_{-\alpha},
U_\beta , U_{-\beta}\rangle   $ of $ L $.  Moreover$,$ 
either
\begin{enumerate}
 \item [(a)] $C_{\alpha ,\beta}$ centralizes $L,$
\item  [(b)]  $q$ is odd and $X_{\alpha ,\beta}=
 \SL(2,q)\circ \SL(2,q)$ or $
\SL(2,q)\times\SL(2,q)$ with center  $C_{\alpha ,\beta},$
inside $\Spp(4,q),$ $\SU(4, q)  $  or
$G_2(q),$
\item  [(c)] $q$ is odd and   $X_{\alpha ,\beta}=\SL(2,q)\circ
\SL(2,q^3)$ with center $C_{\alpha ,\beta},$  inside $^3\! D_4(q),$  
\item  [(d)]   $q$ is odd and   $X_{\alpha ,\beta}=\SU(4,q) $ with center
$C_{\alpha ,\beta},$ or
\item  [(e)]  $X_{\alpha ,\beta}=\SL(3,q) $ with  $3|q-1,$ inside 
$L=G_2(q)$ or
$^3\!D_4(q)$.
\end{enumerate}

\end{lemma}
 
\Proof 
As already noted,   $U_\alpha=\langle u_\alpha^{\langle h_\alpha \rangle
}\rangle $ since $L$ is not $\SU(5,2)$. 
Since $H$ acts on $L_\alpha$ as a group of automorphisms fixing 
$U_\alpha$ and $U_{-\alpha} $, we then have  $C_H(u_\alpha)= C_H(U_\alpha)=
C_H(L_\alpha). $ The latter equality  can be seen by inspecting each of the
rank 1 groups
$L_\alpha$; but it is more easily understood 
by noting that $C_H(U_\alpha)$
normalizes each member of  
 $ \{U_\alpha\}\cup
(U_{-\alpha})^{U_\alpha}=(U_\alpha)^{L_\alpha}$.

Consequently,  $C_{\alpha ,\beta}$
centralizes  the rank $2$ subgroup $X_{\alpha ,\beta}; $  we may assume
that $C_{\alpha ,\beta}$ does not centralize $L$.   There  are
very few orbits of  unordered  distinct  pairs 
$\{\pm\alpha\},\{\pm\beta\} $   under the action of the  Weyl
group of
$L$ on $\Phi_{\alpha_i,\alpha_j}$. 
 When $G$ does not have type  $^2\!F_4(q)$ it  
straightforward to calculate $X_{\alpha ,\beta}$ in each case and
then to check that it  
appears in (b)--(e). 
The case $G=L=\,^2\!F_4(q)$ cannot occur in view of the   orders of the 
centralizers of semisimple elements obtained in  
\cite{Shi}.
\qedd 
\medskip

We now describe the  version of the Curtis-Steinberg-Tits presentation
  using commutator relations  \cite{Cu} (cf.
\cite[Sec. 2.9]{GLS}).  
For all roots
$\alpha ,  \beta \ne \pm\alpha $ lying in the same rank 2 subsystem 
spanned by a pair  of members of
$\Pi$ we use the following two types of relations:

\smallskip
\begin{itemize}
\item[]
\begin{itemize}
\item[$(P_\alpha)$] A presentation ${\bf P}_\alpha $ of  $L_\alpha $,
including the   elements $u_\alpha$, $h_\alpha$ and    $r_\alpha$
described above.  Let
$U_\alpha $   be as above. 
\item[$(B_{\alpha ,\beta})$] $[ x_{\alpha }, x_{\beta }  ]= $ product of
elements of root groups
$U_{i\alpha+j\beta}$ with $i\alpha+j\beta \in \Phi$, $i,j>0$, for each
nontrivial  $x_\alpha \in U_\alpha ,x_\beta \in U_\beta $.
\end{itemize}
\end{itemize}
\smallskip
Here and below, the terms of the product are assumed to be exactly the
ones that occur in the given group $G$.  Precise formulas for these
products are given in \cite[Sec.~10]{St1}, \cite{Gr}   and
\cite[Theorem~2.4.5]{GLS}.

This gives the standard presentation for  $G$. 
We will soon  shorten this presentation with the help of one  further
type  of  relations for the same pairs  $\alpha ,\beta $:%
\smallskip
\begin{itemize}
\item[]
\begin{itemize}
\item[$(H_{\alpha ,\beta})$] $x_\alpha^{h_\beta}=\star_\alpha $ for each
nontrivial $x_\alpha \in U_\alpha $, where
$\star_\alpha \in U_\alpha $ depends on both $x_\alpha $ and $h_\beta $.
\end{itemize}
\end{itemize}
\smallskip

\noindent
The crucial if simple observation is that, since  
$ \langle h_\alpha,  h_\beta \rangle $  normalizes $U_\alpha$ and
 $U_\beta$ by $(P_\alpha)$, $(P_\beta)$  and  $(H_{\alpha
,\beta})$, it 
  acts by conjugation on the set
of relations  $(B_{\alpha ,\beta})$.  Hence we will use (some of) the
relations
  $(H_{\alpha ,\beta})$ to deduce all of the relations 
$(B_{\alpha ,\beta})$ from a small number  of them.
Specifying the latter subset is mostly a matter of
bookkeeping.

We break the proof of Theorem~\ref{Short presentations} into several
cases.

\smallskip{
\noindent \bf Rank $2$ case:~} 
Define $z=z(G)$ as follows:
\smallskip

\noindent\hspace{-4pt} 
\begin{tabular}{cccccccccccc} & \\[-4mm]
$G =$ &  \hspace{-8pt}  $\SL(3,q)$ & $\Spp(4,q)$ &$G_2(q)$ & $\SU(4,q)$ &
$\SU(5,q)$ &
$^3\!D_4(q)$ &  $^2F_4(q)$%
 \vspace{2pt}
\\
  $z\hspace{2pt} =$ & \hspace{-8pt}  $(3,q-1)$   &$(2,q-1)$& $(3,q-1)
$&$(4,q+1)$&$(20,q+1)$ &$(3,q-1)$&1%
\end{tabular}
%

\smallskip
\noindent
so that   $z\le 20$. 

Recall that $H=\langle h_{\alpha_1},  h_{\alpha_2} \rangle $ acts on each
$\bar U_\alpha $ as a group of 
field multiplications
in  an extension field   $\F_{q^{e(\alpha)}}^*$.
Let  $x_{\alpha \mu}\in U_\alpha ,$ $1\le \mu\le z$,  denote elements
corresponding to coset representatives  in  
$\F_{q^{e(\alpha)}}^*$   of  the $z$th   powers of all  elements of
$\F_{q^{e(\alpha)}}^*$ (or projecting onto such coset representatives 
in $\bar U_\alpha$ for
nonabelian $U_\alpha $ in the $\SU(5,q)$ and $^2\! F_4(q)$  cases).  

We now replace the preceding presentation by the following
one for a group $J$ (using the same pairs $\alpha ,\beta $  as above).
\smallskip
\begin{itemize}
\item[]
\begin{itemize}
\item[$(P_\alpha) $] As before.  This also provides us with 
the various $x_{\alpha \mu}$, which we may assume are among the
generators of
$L_\alpha$.

\item[$(B_{\alpha ,\beta}')$]   $[ x_{\alpha \mu}, x_{\beta \nu}  ]= $
product of elements of root groups
$U_{i\alpha+j\beta}$ with $i\alpha+j\beta \in \Phi$, $i,j>0$, for all
$\mu,\nu$.

\vspace{2pt}

\item[$(H_{\alpha ,\beta}')$] $x_{\alpha \mu}^{h_\beta}=\star_\alpha $
for all $\mu$, where 
$\star_\alpha \in U_\alpha $ depends on both $x_{\alpha \mu}$ and $h_\beta
$.
\end{itemize}
\end{itemize}

\smallskip

{\em We will show that $J\cong G$.}  
Since $G$ is a homomorphic image of $J$, for each root $\alpha$
we have a subgroup  $L_{\alpha}$ behaving as in \eqn{h order}.

Once again,
  $U_\alpha=\langle x _{\alpha }^{\langle h_\alpha \rangle
}\rangle $ if $x_{\alpha }\in U_\alpha\backslash
 \Phi(U_\alpha)$.  
The action of  $H$  on each $\bar U_\alpha $ as
a group of 
field multiplications is encoded in $(P_\alpha)$, $(P_\beta)$  and
 $(H'_{\alpha ,\beta})$.

For some perpendicular  pairs  $\alpha,\beta$  the relations $(B'_{\alpha
,\beta})$ state that   $[ x_{\alpha \mu}, x_{\beta \nu}  ]= 1$. 
For such a  pair  $\alpha,\beta$  we only need one  such relation with $x
_{\alpha \mu}
\in U_\alpha\backslash
 \Phi(U_\alpha)$ and $x _{\beta \nu}\in U_\beta\backslash
\Phi(U_\beta)$, since  
conjugating  by all elements of $\langle h_\alpha, h_\beta \rangle$  
then implies that 
$[U_\alpha,  U_\beta]=1$, and hence that  all relations 
$(B_{\alpha ,\beta})$ hold.%

We deal with $(B'_{\alpha ,\beta})$ for  all remaining pairs   
$\alpha,\beta$ in a somewhat similar manner.  
 Lemma~\ref{H-orbits} implies that $|C_{\alpha ,\beta}| $ 
\emph{divides} $z$. 
By  \eqn{h order},  each orbit of $\langle h_\alpha ,h_\beta \rangle$  on 
${\bar U_\alpha^\#\times \bar U_\beta^\#}$ has size  
$|\bar U_\alpha^\#\times \bar U_\beta^\# |/| C_{\alpha , \beta}|$.
Hence,   our pairs $( x_{\alpha \mu}, x_{\beta
\nu})$ include  representatives for all
$\langle h_\alpha ,h_\beta \rangle$--orbits on 
$\bar U_\alpha^\#\times \bar U_\beta^\#$.


Conjugating   the relations  $(B_{\alpha ,\beta}')$ by all  elements of 
$H$, and using 
$(H_{\alpha ,\beta}')$,  we  obtain all relations
$(B_{\alpha ,\beta})$ if $U_\alpha $  and  $U_\beta $ are
abelian. When $U_\alpha $ or $U_\beta $  is nonabelian,
 for each pair 
$( x_{\alpha }\Phi(U_\alpha), x_{\beta }\Phi(U_\beta))$
of cosets
we still obtain a  relation of the form $(B _{\alpha ,\beta})$. 
 By using 
 the elementary identity $[x,uv]= [x,v] [x,u]^v$
in both  $G$ and  the presented group  $J$, we see that  our conjugates
of the relations
$(B_{\alpha ,\beta}')$   imply all relations $(B_{\alpha ,\beta} )$.

Thus, all
relations  $(P_{\alpha })$  and $(B_{\alpha ,\beta})$ required in the
Curtis-Steinberg-Tits presentation hold for  $J$, so that
$J\cong G$.  

 There are at most $20^2$  relations $(B_{\alpha ,\beta}')$
and $20$  relations $(H_{\alpha ,\beta}')$  for each choice of the roots
$\alpha ,\beta
$. The length of this presentation is $O(\log q)$: each factor appearing
in  
 $(B_{\alpha ,\beta}')$, and each $\star_\alpha$, has  that length by 
Theorem~\ref{PSL2q length}(b).

\smallskip
{\noindent \bf General case {\rm (excluding $\SU(2n+1,2)$)}:} 
Let $z(G) $ be the least common multiple of 
the integers $z(L_{\alpha_i,\alpha_j}) $ 
for the various rank 2 groups  $L_{\alpha_i,\alpha_j} $ just considered.
Introduce elements   $x_{\alpha \mu}\in U_\alpha ,$ $1\le \mu\le z(G)$,
that behave as before.

 Theorem~\ref{Curtis-Steinberg-Tits presentation} uses
all
$L_{\alpha_i,\alpha_j},$  including ones of type
$A_1\times A_1$ occurring when  the roots $\alpha_i$  and  $\alpha_j$ are
perpendicular.  Therefore we need to consider the corresponding
relations $(B_{\alpha ,\beta}')$ and
$(H_{\alpha ,\beta}')$ for  $\alpha $ and $\beta$ in the root system
spanned by $\alpha_i$ and $\alpha_j$.   As above, in the $A_1\times A_1$ 
case one   relation
$(B'_{\alpha ,\beta})$ with    
 $x _{\alpha \mu}
\in U_\alpha\backslash
 \Phi(U_\alpha)$ and $x _{\beta \nu}\in U_\beta\backslash
\Phi(U_\beta)$
 implies all relations $(B_{\alpha ,\beta})$.


As  before we see that all $(P_{\alpha } )$, $(B_{\alpha
,\beta}')$
and $(H_{\alpha ,\beta}')$  imply all relations $(B_{\alpha
,\beta} )$.  Once again, $J\cong G$.

\smallskip
{\noindent \bf The case $\SU(2n+1,2)$:}
For each long root $\alpha $ lying in the span of two fundamental roots, 
we introduce two elements
$x_{\alpha \mu}$ that   generate the quaternion 
 group $U_\alpha $. 
 We can now proceed exactly  as above.

\smallskip
{\noindent \bf Summary:}
Each relation $(B_{\alpha ,\beta}')$  involves a bounded number of 
elements of various 
rank 1 groups $L_\gamma$; each of these elements has length $O(\log q)$ in
the generators of $L_\gamma$   by Theorem~\ref{PSL2q length}(b).
A similar statement holds for $(H_{\alpha ,\beta} ')$.
We needed to consider   $O(n^2)$ pairs $\alpha ,\beta $.

\smallskip
{\noindent \bf Centers:}
While this is  a presentation for the universal  group $G$,
   the center can be  killed exactly as in \cite[Sec.~5.2]{BGKLP}.  For
example, for groups of type  $A_n$ there is a standard   product  of the
form
$\prod _1^{n-1} h_{\alpha_i}(\xi^i)$,  $|\xi|=(q-1,n+1)$, 
  that generates the
center.
By Theorem~\ref{PSL2q length}(b), this expression has length
$O(n\log q)$ in our generators.
Thus, within our length
requirements we can factor out all or part of the resulting  cyclic
group  (a very different approach to this is used later in
Section~\ref{Theorems A and B}).  For the remaining types
it is also easy to write   the required central elements as short
words in our generators.

 We have now handled all cases of
Theorem~$\ref{Short presentations}$. ~\qedd
\smallskip
\smallskip

Note that  this presentation often  involves
 more elements $x_ {\alpha\mu}$ and relations 
 $(B_ {\alpha, \beta} ')$
than are actually needed.  For example, for odd-dimensional unitary
groups  we never need 20 coset representatives:  we only need at most 4
or 5 for any pair $\alpha ,\beta $ of roots (cf. Lemma~\ref{H-orbits}).
Different examples occur in  \cite{GKKL2}.

\Remark
\label{Elements of SLn}    \rm
For future reference we  mention   two elements of $\SL(n-1,q)$ that
can be readily found using the above generators.
By
Theorem~\ref{PSL2q length}(b),
\vspace{-4pt}
$$  
c_{23}:=r_{\alpha_2}^2r_{\alpha_1} =
\begin{pmatrix}
     \, 1 &  0 & 0  \\
     \, 0 & 0 &  1 &O \\
     \, 0 &  \!\!\!\!-1 & 0  \\
        &  O &  & I  \\
     \end{pmatrix}
\begin{pmatrix}
     \,\,\,\,0 &  1 & 0  \\
      -1 & 0 &  0  &O \\
       \,\,\,\,0 &  0 & 1  \\
        &  O &  & I  \\
     \end{pmatrix}^{\!\!2}
=
\begin{pmatrix}
     -1 &  0 & 0  \\
      \,\,\,\, 0 & 0 &  1 &O \\
      \,\,\,\, 0 &  1 & 0  \\
        &  O &  & I  \\
     \end{pmatrix}
    $$
is an element of $L_{\alpha_2}L_{\alpha_1}$, and hence  has length 3 in
our generators.  In the  monomial action  on the standard basis
$e_1,\dots, e_n$,
   we have  $c_{23}=(e_1,-e_1) ( e_2,e_3)$.  Hence,
$$c_{2345}:=(e_1,-e_1) ( e_2,e_3 )\cdot
(e_1,-e_1) (  e_3,  e_4 )\cdot
(e_1,-e_1) ( e_4,e_5) =(e_1,-e_1) ( e_2,e_3,  e_4,e_5) $$
   has length 9.  Compare Proposition~\ref{short words}.

\subsection{Word lengths}
\label{Word lengths}
The preceding remark is a very special case of the   
 following observation:

\begin{Proposition}
\label{short words}
If $n$ is bounded in {\rm Theorem~\ref{Short presentations},}
then every element of every central extension of $G$ is a word of length
$O(\log q)$ in the generators used in the  theorem.
\end{Proposition}

\Proof
The elements    $r_{\alpha_i}$   
 generate the Weyl group
$W$ modulo $H$.  By hypothesis, $|\Phi|$
and
$|W|$ are bounded.
By the Bruhat decomposition, $\tilde G=\bigcup _{w\in W}BwB$ with
$B=UH$ a Borel subgroup. Here $U$ and $H$ are products of fewer
than $  | \Phi|$ subgroups of
$L_\alpha ,
\alpha \in
\Phi$.  Since each $L_\alpha $ is a $W$-conjugate of some $L_{\alpha_i}$,
Theorem~\ref{rank 1}(b) implies  the result.~\qedd

\section{Theorem~\ref{A} }
\label{Bounded presentations for groups of Lie type}
\label{Theorems A and B}
We  now turn to Theorem~\ref{A}. 
By Theorem~\ref{Short presentations},    we only need to consider the
classical groups of rank greater than $  8$.

   \subsection{{A presentation for \boldmath
$\SL(n,q)$} }
\addcontentsline{toc}{subsection}
{\protect\tocsubsection{}{\thesubsection}{A presentation for
$\SL(n,q)$}}
\addtocontents{toc}{\SkipTocEntry}

\label{Bounded presentation for SL(n,q)}
\label{A presentation for SL(n,q)}
\label{SL(n,q)}
   We begin with the case  $\SL(n,q)$.  We would like to use the Weyl group
$S_n$ in our presentation, but $\SL(n,q)$ does not have a   natural
subgroup  $S_n$  when $q$ is odd.  There are various ways around this
difficulty, such as using a subgroup $(2^n\semi S_n)\cap \SL(n,q)$
or the alternating group (compare \cite{GKKL2}).   We
have chosen to use $S_{n-1}$ as an adequate substitute for
$S_n$.

\begin{Theorem}
\label{Bounded SLn}
All groups  $\SL(n,q) /Z, $  where  $Z\le Z(\SL(n,q)),$
have
bounded presentations of length
   $O(\log n +\log q)$.
\end{Theorem}

Clearly the most interesting cases are $Z=1$ or $Z(\SL(n,q))$.  However,
later we will also need the case $|Z|=2$.
\medskip

\Proof
We start with two  bounded  presentations of length  $O(\log n +\log q)$
given in Theorems~\ref{All symmetric groups} and
\ref{Short presentations}:
\smallskip

   \begin{itemize}
   \item{}
$T=\langle  X\mid R(X)\rangle$
of $S_{n-1}$ (acting on $\{ 1,2,\ldots,n    \}$, fixing 1),  and
   \item{}
   $F=\langle \tilde X\mid\tilde R\rangle$ of
$\SL(5,q)$.

The notation $R(X)$ is used here since later we will need to use 
a second copy $\langle \bar  X\mid R(\bar X)\rangle$  of $S_n$.

    We assume that     $X$ and $\tilde X$ are disjoint, and that these
 presentations
satisfy the following additional  conditions for some $ X_1,X_2\subseteq
X$ and $\tilde Y\subset \tilde X$:

\begin{itemize}
   \item [\rm(i)] $\langle X_1\cup X_2  \rangle$ projects onto
the stabilizer $T_2 $  of 2.
   \item [\rm(ii)]   $X_1$ projects into
   $A_{n-1}$,  and each element of $X_2 $
projects outside
$  A_{n-1}$.
%
%
   \item [\rm(iii)] $x_{ (2,3) },$  $ x_{ (2,3,4,5)} ,$  $x_{ (6,7)},$  $
x_{ (5+\delta,\dots, n) } $ are  words of length  $ O(\log n)$ in $X$ that
project onto $(2,3) ,$  $ (2,3,4,5) ,$  $(6,7), (5+\delta,\dots, n)$,
respectively,  where  $\delta=(2,n-1)$.
   \item [\rm(iv)]
$x_\sigma $ is a word of length  $ O(\log n)$ in $X$ and projects onto
$\sigma =(2,\dots,n  )\,$ (this is needed only for handling the center of
$\SL(n,q)$).

   \item [\rm(v)]  The elements $u,  t,   h_\z$ in \eqn{generators}
and  \eqn{SL2 generators},
that generate a subgroup $L=\SL(2,q)$, are in $\tilde X$. 
When $q$ is odd, we also assume that $h_2\in \tilde X$.  
(Here we are using the subgroup
$
   \begin{pmatrix}
        \SL(2,q) & O \\
      O & I \\
     \end{pmatrix}
    $
of $\SL(5, q)$.)
%
   \item [\rm(vi)] $c_{ 23}$ and  $ c_{ 2345}   $
are words of length $O(\log q)$ in $\tilde X$ acting as in
Remark~\ref{Elements of SLn}.
   \item [\rm(vii)] $\langle \tilde Y\rangle =\SL(4,q)$,   $|\tilde Y|=2$
and the   members of $\tilde Y$  have length $O(\log q)$ in $\tilde X$.
(Here we are using the subgroup
$
   \begin{pmatrix}
       1 & O \\
      O &  \SL(4,q) \\
     \end{pmatrix}
    $
of $\SL(5, q)$.)
%

 \item [\rm(viii)]   $X$ contains a set of generators of the stabilizer
$T_{23 }=T_{2 }\cap T_{ 3}$ of both  2 and 3.  (This will only be needed   later
when we deal with some orthogonal groups  in Section~\ref{Generic case}, Case
2.)
   \end{itemize}
   \end{itemize}
\smallskip

\emph{Existence}:
    In Remark~\ref{needed permutations}(i) we constructed
the permutations in (iii) and (iv);  in Remark~\ref{needed permutations}(ii,
iii) we  noted that conditions (i) and (ii) hold;
and    Remark~\ref{needed permutations}(iii) also takes care of (viii).
   Remark~\ref{Elements of SLn} finds the elements in (vi).  Note that the
group $L =L_{\alpha_1}$, and a  presentation for it  using the 
generators   
  in (v) among others, were essential
ingredients in Section~\ref{Commutator relations} when $q>9$. For   (vii)
see Proposition~\ref{short words}.

Now our presentation is as follows:
\smallskip

{\noindent \bf Generators:}  $X\cup \tilde X$  (we are
thinking of  $T$ and  $F$ as
embedded in $\SL(n,  q)$ as
$
   \begin{pmatrix}
        \pm1 & O\\
        \,\,\,\,O & $permutations$   \\
     \end{pmatrix}
    $ %
   using permutation matrices,
and
$
   \begin{pmatrix}
        \SL(5,q) & O \\
      O & I  \\
     \end{pmatrix}
   \! ,$ respectively).

\smallskip
{\noindent \bf Relations:}
\begin{itemize}
\item []
\begin{enumerate}
    \item $R(X)\cup \tilde R$.
   \item $ u^{x_1}=u, \bar u^{x_1}=\bar u,  $ for all $x_1 \in X_1$
(where we have abbreviated $\bar u: = u^t$).
    \item $u^{x_2}=   u^{-1}$, $\bar u^{x_2}=\bar u^{-1}$, for all $x_2 \in
X_2$.
    \item $[h_\z,X_1\cup X_2]=1$.
     \item $x_{ (2,3) }=c_{23}, x_{ (2,3,4,5)}=c_{2345} .$
   \end{enumerate}
\end{itemize}

\smallskip
We will  show that {\em the group $G$ defined by this presentation
is isomorphic to  $\SL(n,q).$}   There is a natural surjection
$\pi\colon G \to \SL(n,q)$. (For, in view of (ii), relation  (3)  
express  the fact  that each odd permutation maps the first basis vector
to its negative, while the even permutations in  (2)  fix  that
vector, as do the elements of $\tilde Y$.)

By (1) and  Lemma~\ref{It's a subgroup},
$G$ has subgroups we can identify with $T=\langle  X \rangle $ and
$F=\langle  \tilde X \rangle $.

We have   
$ \langle u, \bar u, h_\z  , x_{ (2,3,4,5)}\rangle=F $
since  $\langle u, \bar u, h_\z\rangle$ is $L=\SL(2,q)$.
By (i), (2), (3)  and (4),   
$X_1\cup X_2$ acts on 
$ \{ u^{\pm1}, \bar u^{\pm1}, h_\z \} $,
so that  the elements
  $x_{ (6,7) }$ and  $x_{ (5+\delta,\dots, n) }$ of 
 $T_2=\langle X_1\cup X_2\rangle$  
act     on $F$.  
Moreover,   $\langle   x_{ (2,3) } , x_{ (2,3,4,5)} 
\rangle =
\langle   c_{23}, c_{2345} \rangle <F$  by (5).  Thus,   
$$
|N_T(F) | \ge | \langle   x_{ (2,3) } , x_{ (2,3,4,5)}, x_{ (6,7) },
x_{ ( 5+\delta,\dots, n) }\rangle |=
|S_4 \times S_{n-5}|,$$
 while
$|F^T|\ge|\pi(F^T)| =
\binom{n-1}{4}$.
It follows that
   $T$ acts on $F^T$ as it does on the set of all
4-sets in  $\{2,\dots,n\}$.

   The  subgroups  $U_{12}=\langle u^{\langle h_\z\rangle  }\rangle $ and
$U_{21}=\langle  \bar u^{\langle h_\z\rangle  }\rangle $
of $F$  are root groups  order $q$.   
Once again, by (i),   (2), (3) and (4), 
$T_2=\langle X_1\cup X_2\rangle$
normalizes each of them.
As above, it follows that $T$ acts on both $(U_{12} )^T $
and $(U_{21} )^T $    as it does on
$\{2,\dots,n\}$.

Both  $U_{31} = (U_{21})^{c_{23}}$ and
   $U_{23}=[ U_{12} ,  U_{31}]$ are  also root groups  of $F$.
Since    $(T_2)^{c_{23}} =T_3$ normalizes $ (U_{21})^{c_{23}} =U_{31}$,
it follows that
$T_{2}\cap  T_{3}$ normalizes $ U_{23}$.   As above, we find that
   $| (U_{23}) ^T|=(n-1)(n-2)$.

 Thus, $T$  acts on both $(U_{12} )^T $ and $(U_{21} )^T $
as it
does on singletons, and on   $(U_{23}) ^T $ as it
does on  ordered pairs from  our $(n-1)$-set.
By the 4-transitivity of $T$,  any given  pair from
$ (U_{12} )^T\cup (U_{21} )^T\cup(U_{23}) ^T $ can be conjugated  into $F$ by
a single element of
$T$.

Consequently,  $N:=\langle L^T  \rangle\cong \SL(n, q)$ by the
Curtis-Steinberg-Tits presentation
(as in  Section~\ref{Commutator relations}  or
Lemma~\ref{3-transitive part 2}).
Moreover, $N\unlhd G$, and $G/N$ is a homomorphic image of $\langle  T
\rangle \cong S_{n-1}$ in which, by (5),  a transposition is sent to 1.
Thus, $G/N=1$.

{\em This proves the theorem when $Z=1$.}

Now that we have a presentation for $\SL(n,q)$, we need to factor out
an arbitrary subgroup
$Z$ of its center.  However,  a  generator of $Z$ probably  cannot 
 be written as a
short word  in our present generators. In order to deal with this
obstacle we will use the following lemma to increase our  generating
set:


\begin{lemma}
\label{MARTIN'S  LEMMA}
Suppose that $G$ is a finite group   containing the wreath product\break 
$ H\wr S_m  =S_m \ltimes H^m $ for some group $H$.   Let
   $\langle Y \mid S \rangle$ be a presentation for $G$ of length $<l$ and
$\langle \bar{Y} \mid  \bar{S} \rangle$
    a presentation for $S_m$ of length $<l$.
   Fix $\epsilon  \in H$ and let
$h_{i,j}  $    be the
element in $H^m $   with coordinates $1$ everywhere aside from $i$th
coordinate
$\epsilon $ and $j$th coordinate $\epsilon ^{-1}$.
Assume  that
\begin{itemize}
\item [\rm(a)] The following elements of 
the standard subgroup  $S_m< S_m \ltimes H^m $
can be written as   words of length $<l$ in $Y$$:$
   $~\sigma
=(1,2,\dots,m), \,
\mu=(2,\dots,m) $   and   $\tau_2=(2,3);$
\item [\rm(b)]
   The  elements
$\bar\sigma $ and   $\bar\tau_2$
in $\langle \bar{Y} \mid  \bar{S} \rangle$
corresponding to the above ones
   can be written as  words of length $<l$  in $\bar{Y};$ and
\item [\rm(c)] $h_{1,2}$ can be written
as a   word of length $<l$  in $Y$.
\end{itemize}

Then the following is a presentation for $G,$
of length $ < 13 l $$:$
\smallskip

{\noindent\bf Generators:} $Y \cup \bar{Y} \cup \{d\}$.
\smallskip

{\noindent \bf Relations:}
\begin{itemize}
\item []
\begin{enumerate}
    \item $S\cup \bar{S}$.
    \item $\bar \sigma = \sigma d$.
    \item $\bar \tau_2 = \tau_2 h_{2,3}$.
    \item $[d,\tau_2] =  [d,\mu] = [d,h_{2,3}] =1 $.

\end{enumerate}

\end{itemize}
\smallskip

\noindent
Moreover$,$ $d$ maps to the element $(\epsilon ^{1-m}, \epsilon ,
\ldots, \epsilon ) \in H^m < G$.
   \end{lemma}

Note that $Y$ and $\bar Y$ are unrelated. 
  It is the subgroups 
$\langle \sigma,\tau\rangle$ and 
$\langle \bar \sigma,\bar \tau\rangle$, both isomorphic to $S_m$, that 
are related via the isomorphism ``bar''; and we will see that they are
conjugate in  $S_m \ltimes H^m $.
Also note that $\epsilon$ is used to define the element $h_{2,3}$
and hence is essential for (3) and (4).
\smallskip

\Proof
If $\pi\in S_m$ and $(a_1,\dots,a_m)\in H$, then
 in $S_m \ltimes H^m$ we have
$$
  (a_1,\dots,a_m)\pi = \pi  (a_{\pi(1)},\dots,a_{\pi(m)}).
$$
In particular,
$(a_1,\dots,a_m)^{\sigma^{-1}}=  (a_m,a_1,\dots,a_{m-1})$
is ``pushing to the right".

  Let $J$ be the group presented above.
It is straightforward to check  that $J$ surjects onto
$G$, by taking $d$ to be  $(\epsilon ^{1-m}, \epsilon ,
\ldots, \epsilon )$    and letting
the new copy of $S_m$   be conjugate to the standard copy
via   $ (\e^{-1},\epsilon^{m-2},  \ldots,\e^{2},\e,
1) \in H^m$.
It is also easy to check that this presentation has length  $ < 13 l$.

By Lemma~\ref{It's a subgroup},
$J$ has    subgroups  we can identify with  $ G= \langle Y \rangle$
and $S_m=\langle \bar Y\rangle $.
It suffices to prove that $J\le  G$.

 Relations (4)   imply
that $d$ commutes with $\tau_k=(k,k+1) $  and
$h_{k,k+1}$ whenever $2\leq k <m $.  Hence, by (2) and (3),
$$
\bar \tau_3 = (\bar \tau_2)^{\bar \sigma^{-1}} =
(\tau_2 h_{2,3})^{ d^{-1} \sigma^{-1} } =
\left(\tau_2^{  d^{-1} } h_{2,3}^{ d^{-1}  }\right)^{\sigma^{-1}} =
\left(\tau_2 h_{2,3}\right)^{\sigma^{-1}}= \tau_3 h_{3,4}.
$$
Similarly, induction gives
$
\bar \tau_k = \tau_k h_{k,k+1}
$ whenever $2\leq k <m $.
A similar equation also holds for $\tau_m:=(m,1)$$:$
$$
\bar\tau_m =(\bar\tau_{m-1} )^{\bar \sigma^{-1}}  =
(\tau_{m-1}h_{m-1, m} )^{ d^{-1} \sigma^{-1} } =
(\tau_{m-1}h_{m-1, m} )^{  \sigma^{-1} } =
\tau_m h_{m , 1}.
$$
Consequently, $\bar\tau_{2},\dots,\bar\tau_{m }\in S_m \ltimes H^m \le
G$.  It follows that our second copy  $ \langle \bar Y\rangle $ of $S_m$ lies in $G$,
and hence so do    $ \bar \sigma $,  $\bar Y $   and $d$ (by (2)).
Thus,  $J =\langle Y \cup \bar{Y} \cup \{d\}\rangle  \le  G$, as required.

In order
to explain the motivation used here, we include the
following direct calculation that
  $d$   behaves as  desired:
\begin{align*}
\bar \sigma & =
\bar \tau_{2} \bar \tau_{3} \cdots \bar \tau_{m-1} \bar \tau_{m} \\
&  =
\tau_{2} h_{2,3} \tau_{3} h_{3,4} \cdots \tau_{m-1} h_{m-1,m} \tau_{m}
h_{m,1}\\ &=
\tau_{2}\tau_{3}\cdots \tau_{m-1} \tau_{m}
h_{2,1} h_{3,1} \cdots h_{m-1,1} h_{m,1}
\\
d& =  \sigma^{-1}\bar\sigma = h_{2,1} h_{3,1} \cdots h_{m-1,1}
h_{m,1}.~\qedd
\end{align*}

\medskip
When we use this lemma, the bound $l$ will be taken to be $O(\log n + \log
q)$.

\medskip
\noindent
\emph{Completion of the proof of} Theorem~\ref{Bounded SLn}:
Let
$Z=\langle \e I\rangle $.  In the preceding lemma let $m=n-1$ and let $H
= \langle \epsilon\rangle$ be the subgroup of
$\F_q^*$ of order $|\epsilon|$.  Embed $ H^{n-1}$ into
$  G=\SL(n,q)$ as all $\diag(\e_1,\dots,\e_n)$ with $\e_i\in\langle \e\rangle $ and $\prod
\e_i=1$.   We already used a presentation
$ \langle X\mid R(X)\rangle$ for
$S_{n-1}$  above in (i), where  $S_{n-1}$  permutes the
last
$n-1$ coordinates. We  use the
corresponding  presentation $\langle \bar  X\mid R(\bar X)\rangle$   for
a second copy of $S_{n-1}$  (with $X$ and $\bar X$ disjoint).

First we need to verify conditions (a)-(c) of the lemma.
 We temporarily use the notation in
Remark~\ref{parities}(i):
$\sigma=({\bf 2},{\bf 3} ,\dots,{\bf n}) $ and
$z=({\bf 2},{\bf 3}) $ have the required length, hence so do
$\mu= ({\bf 2},{\bf 3})\sigma$  and  $ ({\bf 3},{\bf 4})=({\bf 2},{\bf
3})^{\sigma^{-1}}$. It follows that Lemma~\ref{MARTIN'S  LEMMA}(a) holds,
and hence so does  (b).
 It remains to consider  (c).
The group $L=\SL(2,q)$ was  defined  above in (v).
Use Theorem~\ref{rank 1}(b) to write the element
$d'=\diag(\e^{   },\e^{  -1},1,\dots,1) \in L$  as a word of length
$O(\log q)$ in our generators   (cf.~(v)).
Then
$h_{2,3}= d' {}^{\sigma^{-1}}$ has the  length required  in (c).

The
lemma provides a new bounded presentation for $G$
of  length $O(\log
n +\log q)$, including a new generator $d$
representing the diagonal
matrix $\diag(1,\e^{2-n},\e,\dots,\e)$.
 Then $Z=\langle dd'\rangle $, so that
the additional
relation
$dd'=1$ produces the desired factor group.
\qedd

\subsection{Generic case}
\label{Generic case}
In this
section we obtain  short bounded presentations for the
universal
central extensions of the simple groups of Lie type.  This is 
significantly simpler
than dealing with the simple groups, which involves  
factoring out the centers of
the universal extensions.  However, in
some cases the latter is easy: it is a
matter of choosing slight
variations on the subgroups  we amalgamate.

We begin with some
general properties of
these universal central
extensions:

\begin{Lemma}
\label{Properties of universal central
extensions}
Let $G$ be a simple classical group defined on a vector
space of dimension
$>8 ,$ and let
$\hat G$ denote its  universal
central extension.
\begin{enumerate}
\item  If $G=\PSL(n,q)$ then
$\hat G=\SL(n,q)$.
\item  If $G=\PSU(n,q)$ then $\hat
G=\SU(n,q)$.
\item  If $G=\PSp(2n,q)$ then $\hat G=\Spp(2n,q)$.
\item
If $G=\Omega(2n+1,q)$ then $|Z(\hat G)|=(2,q-1) $.
\item  If $G=
\PO^+(2n,q) $ then $|Z(\hat G)|=(4,q^n-1),$
and $Z(\hat G)$ is cyclic
unless $n$ is even.
\item  If $G= \PO^-(2n,q) $ then $ Z(\hat G)$ is
cyclic of order
$(4,q^n+1) $.
\item Suppose that $G=\Omega(V)$ is an
orthogonal group over a field of odd
characteristic.  Then $Z(\hat
G)$ has an involution$,$ the {\rm
  spin involution,}  lying
in
$Z(\hat H )$  whenever $H=\Omega(U)$  for   a  nondegenerate
subspace $U$
of  $V$  of dimension at least
$3$.
\end{enumerate}
\end{Lemma}
\Proof
See  \cite[pp. 312-313]{GLS}
for (1)-(6) and  \cite[Proposition 6.2.1(b)]{GLS}
for
(7).~\qedd

\smallskip
Note that there are various exceptions to
(1)-(6)  when the
dimension is at most 8  \cite[p. 313]{GLS}.

We
will use the following crucial theorem as well as  variations on
its
proof.

\begin{Theorem}
\label{Presentations of universal central
extensions}
\label{universal}
All universal central extensions of
groups of Lie type of rank $n$ over
$\F_q$  have bounded
presentations of length $O(\log n +\log q)$.
\end{Theorem}

\Proof
By Theorem~\ref{Short presentations}, we may assume that $G$ is a
classical simple group of rank $n >8$.  We
will use the  root system   of  $G$  and 
the
Curtis-Steinberg-Tits  presentation
(Theorem~\ref{Curtis-Steinberg-Tits presentation}). For each root $\alpha $
there is a corresponding
subgroup
$L_\alpha $ that is a  central
extension of $ \PSL(2,q ), $
$ \PSL(2,q^2) $ or
$\PSU(3,q)$. We may
assume that $\Pi=\{ \alpha_1, \ldots, \alpha_{n} \}$ 
is a fundamental system  of roots  (in the standard order), 
where $\alpha_1, \ldots,
\alpha_{n-1}$ have the same length.
\smallskip
\smallskip

{\noindent \bf Case 1: $G$ does not have
type $D_n$.}
Let 
$$\mbox{
$\Pi_1=\{\alpha_1, \ldots, \alpha_{n-1}
\},~
\Pi_2=\{\alpha_1,
\ldots, \alpha_{n-2}, \alpha_n\}$ and
$\Pi_3=\{\alpha_{n-1},
\alpha_{n}\}$.}
$$
$$\mbox{  $\stackrel{1} {\bullet }  \!\! \!
  \rule[.7mm]{8mm}{.2mm}
    \hspace{-4pt}
\stackrel{2} {\bullet } 
   \hspace{3pt}
   \cdots
   \hspace{-3pt}
   \stackrel{n-2} {\bullet }
   \hspace{-10pt}
   \rule[.7mm]{8mm}{.2mm}
    \hspace{-10pt}
    \stackrel{n-1} {\bullet }
     \hspace{-10pt}
   \raisebox{.63ex}{$  \underline{  \rule[.1mm]{8mm}{.14mm}  } $}$
   \hspace{-9pt}
 $ \stackrel{n} {\bullet } $}$$

\Remark
\label{3
subsystems}
  These sets of roots have the following
properties\rm:
\begin{itemize}
\item[(a)] For  $i=1,2,3,$ either
$|\Pi_i| = 2 $ or  $G_i:=\langle
L_{\beta}\mid  \beta \in \Pi_{i} \rangle$
is  of type $\SL_{d},$
$\SL_{d} \times \SL_2$ or $\SL_{d} \times
\SU_3 $ for some $d$;  and
\item[(b)] Each pair from $\Pi$ lies in
some $\Pi_i$.
\vspace{2pt}
\end{itemize}
\smallskip

Namely, (b) is
clear, and $G_1$ is a homomorphic image of $
\SL(n,q)$
or
$\SL(n,q^2)$ (as is seen from the commutator relations
in
Section~\ref{Fixed rank} -- or more precisely  in the case
of
odd-dimensional unitary groups, from the explicit relations
in
\cite{Gr,BGKLP}).  Moreover, $G_2$ has type $\SL_{n-1} \times
\SL_2$ or
$\SL_{n-1}
\times \SU_3$, while $|\Pi_3|=2$, which proves
the remark.
\smallskip

The root $\alpha_i$ ($1\le i<n)$ can be identified with the ordered pair
$(i,i+1)$ of elements of the  set $\{1,\dots,n\}$ on which the Weyl group
$S_n$  of $G_1$  acts.  The corresponding root groups $X_{\pm \alpha_i}$
are just groups of elementary matrices 
of $G_1$, and generate a subgroup
$L_{  \alpha_i}=L_{ i,i+1}\cong\SL(2,q)$ or $\SL(2,q^2)$  acting on the span
of  two of the standard basis vectors
  of  the usual $n$-dimensional module for $G_1$.

  We use the following groups  and
presentations.

\begin{itemize}
\item {}
$G_1 $  is the universal
central extension
obtained in Theorem~\ref{Bounded SLn} (but {\em
reversing} the order of
the set
$\{1,\dots,n\}$), given with a
bounded
  presentation
$ \langle X_1\mid R_1 \rangle$
of length
$O(\log n +\log q)$ including
\begin{itemize}
\item [(i)] generators
inside $X_1$ for  $T=S_{n-1}$  acting on
  $\{1,\ldots,n\}$ and
fixing $n$, and  for 
$T_{n-1}'$, the
stabilizer  of $n-1$ in $T'=A_{n-1}$;
\item
[(ii)]
a    bounded presentation
$\langle X(n-1)\mid R(n-1) \rangle $
of length $O( \log q)$ for the
subgroup  $L_{\alpha_{n-1}} $;
and
\item [(iii)] a set  $X(n-2)$  of generators for the subgroup
$
L_{\alpha_{n-2}}$  of
$G_ 1$.
\end{itemize}

\smallskip
\emph{Existence}:
Note that Remark~\ref{parities}(iii)  provides short generators for
$T'_{n-1}$, as required in  (i).
The   
groups $ L_{\alpha_{n-1 }}$ and $ L_{\alpha_{n- 2}} $
were used in the presentation for $F$ in the proofs
of
Theorems~\ref{Short presentations}  and
 \ref{Bounded SLn},
so that   we also already have $X(n-2)\cup
X(n-1)\subseteq X_1$  and
$ R(n-1)\subseteq
R_1$, as required in  (ii) and (iii).

Then   $ \langle L_{\alpha_{n-2}}, T'_{n-1} \rangle
\cong\SL(n-1,q) $  or $\SL(n-1,q^2) $,
since $T'_{n-1}$   fixes the corresponding standard basis vector
(whereas each odd permutation   in  $T _{n-1}$ sends that vector to its
negative; cf. Section~\ref{SL(n,q)}).

\item
{}
$G_3=\langle X_3\mid R_3 \rangle  $ was obtained
in
Section~\ref{Fixed rank}  by  using the root system $\Pi_3$
and
bounded presentations of length $O( \log q)$ for two rank 1
groups,
namely
\begin{itemize}
\item [(iv)]  the   {\em same\/}
presentation  $\langle X(n-1)\mid
R(n-1)
\rangle
$ used above for
$L_{\alpha_{n-1}}$;
and
\item [(v)]  a   bounded
presentation
$\langle X(n )\mid R(n ) \rangle $
of length $O( \log
q)$
  for  $ L_{\alpha_n} $.
\end{itemize}
 Hence,    we already have  $X(n-1)\cup X(n)\subseteq
X_3$  and  $R(n-1)\cup R(n)\subseteq
R_3$.
\end{itemize}
\smallskip
Our presentation is as
follows:
\smallskip

{\noindent\bf
Generators:}   $X_1\cup   X_3.$

{\noindent\bf
Relations: }
\begin{itemize}
\item
[]
\begin{enumerate}
   \item {}
$R_1\cup  R_3.$

\item
{}
Identify  $  L_{\alpha_{n-1} }$ inside $G_1$ and
inside $G_3$ using
the identity map.
\item {}  $ [ L_{\alpha_{n-2}},  L_{\alpha_n}]=
[ T'_{n-1}, L_{\alpha_n}]=1$.
   \end{enumerate}
\end{itemize}

\smallskip
Then 
  $   G_2 :=
   \langle  L_{\alpha_{n-2}}, T'_{n-1} \rangle\cdot
L_{\alpha_n}
=\langle L_{\beta}\mid  \beta \in \Pi_{1}\cap \Pi_{2}  \,$  or  $
\,\Pi_{2}\cap
\Pi_{3}  \rangle $ is a central product of an $\SL_{n-1}$ and an 
$\SL_{2}$ or  $\SU_{3}$.
In view of  Remark~\ref{3
subsystems}(b), we obtain   the desired
universal central extension
$\hat G$  by
Theorem~\ref{Curtis-Steinberg-Tits
presentation}.  
This presentation is clearly short and bounded.
\smallskip
\smallskip

  {\noindent\bf Case 2: $G$
has type $D_n$.}  We assume that $\alpha_{n-1} $
and  $\alpha_n$ are
connected to
  $\alpha_{n-2}$ in the Dynkin diagram.
This time  we
use
$$\mbox{
$\Pi_1=\{\alpha_1, \ldots, \alpha_{n-1} \}\!,~$
$
\Pi_2=\{\alpha_1,\ldots, \alpha_{n-3}, \alpha_n\}$ and
$\Pi_3=
\{\alpha_{n-2},\alpha_{n-1},
\alpha_{n}\}$. }
$$
  Once again
Remark~\ref{3 subsystems}  holds,
this time  with $G_1,G_2$ and $G_3$
the universal central extensions
$\SL(n,q),$  
$ \SL
(n-2,q)\times
\SL(2,q)$ and $
\SL(4,q),$ respectively.
This time we
use the following presentations.

\smallskip
\begin{itemize}
\item
{}
$G_1 $  is the universal
central extension
obtained in Theorem~\ref{Bounded SLn} (but once again {\em
reversing} the order of
the set
$\{1,\dots,n\}$), given with a
bounded
  presentation
$ \langle X_1\mid R_1 \rangle$
of length
$O(\log n +\log q)$ including

\begin{itemize}

\item [(i)] a    bounded
presentation
$\langle X(n-1)\mid R(n-1) \rangle $ of length $O( \log q)$
  for
$L_{\alpha_{n-1}}  $;
\item [(ii)]
a    bounded
presentation
$\langle X(n-2)\mid R(n-2) \rangle $
of length $O( \log q)$
  for
$L_{\alpha_{n-2}} $;
\item [(iii)]
a
bounded
presentation $\langle X(n-3)\mid R(n-3) \rangle $
of length
$O( \log q)$
  for  $L_{\alpha_{n-3}} $;
  and
\item [(iv)]
generators inside $X_1$ for  $T=S_{n-1}$  acting on
$\{1,\ldots,n\}$ and  fixing $n$,  and  for the
stabilizer
$T'_{n-1, n-2}$
of both $n-1$ and $ n-2$ in $T'=A_n$.
\end{itemize}
The  groups $L= L_{\alpha_{n-1 }}, L_{\alpha_{n
-2}},
L_{\alpha_{n-3 }} $ were used in the presentation for $F$ in the
proof of
Theorem~\ref{Bounded SLn}, so that
$\cup _{n-3}^{n-1}X(i)
\subseteq X_1$  and
$\cup _{n-3}^{n-1}R(i) \subseteq R_1$.
  Generators for   $T'_{n-1 , n-2} $
were obtained  in  Remark~\ref{parities}(iii).
(Recall that the numbering in that remark has been reversed in this
proof.)  This time
$ \langle L_{\alpha_{n-3}}, T'_{n-1 , n-2} \rangle \cong\SL(n-2,q) $.

\smallskip
\item {}  $ \langle X_3\mid R_3 \rangle $   is a
bounded
presentation of length $O( \log q)$ for $G_3  =\SL(4,q)$,
obtained
as in Section~\ref{Fixed rank}
using the presentations (i) and
(ii), together with  the root system $\Pi_3$
and
\begin{itemize}
\item
[(v)] a    bounded
presentation $\langle X(n)\mid R(n) \rangle $ of
length $O( \log q)$
  for  $L_{\alpha_{n}} $.
\end{itemize}
\smallskip

\noindent
Note that
$\cup _{n-2}^{n}X(i) \subseteq X_3$  and
$\cup _{n-2}^{n}R(i)
\subseteq R_3$.

\end{itemize}

\smallskip

Our
presentation is as follows:
\smallskip

{\noindent\bf
Generators: }
$X_1\cup   X_3.$

  {\noindent\bf
Relations: }
\begin{itemize}
\item
[]
\begin{enumerate}
   \item {}
  $R_1 \cup R_3.$
{}
\item {}
Identify
$\langle   L_{\alpha_{n-2} },   L_{\alpha_{n-1} } \rangle \cong \SL(3,q)$ inside
$G_1$ and  inside $G_3$ using the identity map.
\item {}  $ [ L_{\alpha_{n-3}},  L_{\alpha_n}]=  \
[ T'_{n-1, n-2 }, L_{\alpha_n}]= 1 $.
    \end{enumerate}
\end{itemize}
\smallskip

 Once again     $   G_2 : =
    \langle  L_{\alpha_{n-3}}, T'_{n-1, n-2} \rangle\cdot
L_{\alpha_n}  =\langle L_{\beta}\mid  \beta \in \Pi_{1}\cap \Pi_{2}  \,$  
or  $\,\Pi_{2}\cap
\Pi_{3}  \rangle $  
behaves as required.
In view of  Remark~\ref{3 subsystems}(b), we obtain   the desired
universal central extension  $\hat G$   by
Theorem~\ref{Curtis-Steinberg-Tits presentation}. \qedd

\subsection{Symplectic groups}
\label{Symplectic groups}
We will need to factor out the centers of the various groups in
Theorem~\ref{universal}.  The following simple observation  will make this
easy in many cases:

\begin{Lemma}
\label{-1}  Suppose that $q$ is odd.
\begin{itemize}

\item [(i)]
If $n$ is even then  $\SL(n,q)   $  has a  bounded presentation of
length   $O(\log n +\log q)$ in which the   involution $z$ in the center
of
$\SL(n,q)
$  is a word of length  $O(\log n +\log q)$  in the generators.

\item [(ii)]
If $n$ is odd then
   $\SL( n,q)   $  has a  bounded presentation of length
   $O(\log n +\log q)$ in which the central involution $z$ in any given
Levi subgroup $\SL(n-1,q)$   is a word of length  $O(\log n +\log q)$  in
the generators.
\end{itemize}
\end{Lemma}

\Proof  Let $G=\SL(n,q) $.

(i)  In the completion of the proof of  Theorem~\ref{Bounded SLn}, we found a
presentation for $G$ in which $z=dd'$ behaves as stated
when $\epsilon$ is chosen to be $-1$.

(ii) Here we proceed exactly  as in the proof of Theorem~\ref{universal}
for the present group $G$ instead of  the  other classical groups $G$
considered in that theorem,
    using the presentation for $G_1=\SL(n-1,q)$ obtained in~(i).~\qedd

\medskip 
Note that (ii) also can be proved using Lemma~\ref{MARTIN'S
LEMMA}.

   \begin{Corollary}
$\PSp(2n,q)   $  has a  bounded presentation of
length   $O(\log n +\log q)$.
\end{Corollary}

\Proof
The proof of Theorem~\ref{universal} for $\hat G=\Spp(2n,q)$ involved a
presentation for
$G_1=\SL(n,q)$.  This time  we will  use the presentation for
$ \SL(n,q)$ in the preceding lemma.

If
$n$ is even then the additional relation $z=1$
 in Lemma~\ref{-1}(i)  produces the desired
presentation.

Suppose that $n $ is odd.
Let $e_1,\dots,e_n,f_1,\dots,f_n$ be a hyperbolic basis
of the   $\F_q$-space underlying the group  
  $\hat G$.
 Choose this basis   so that    $L_{ \alpha_i}$
has support
$\langle e_i,e_{i+1}, f_i,f_{i+1} \rangle$ for $i<n$ and
$\langle e_n, f_n \rangle$ for $i=n$.
  Let $H = \langle X_{\pm\alpha_i}\mid i\le n-2 \rangle  = \SL(n-1,q)$, with support
the orthogonal complement $\langle e_1,\dots,e_{n-1},f_1,\dots,f_{n-1} \rangle $
of  the support $\langle e_n, f_n \rangle$  of
$L_{\alpha_n}\cong\Spp(2,q)$. By   Lemma~\ref{-1}(ii),   the
involution
$z\in Z(H)$ is  a word of length $O(\log n +\log q)$ in our generators.
On the other hand, the central involution $z'$ of $L_{\alpha_n}$ has length
$O(\log q)$ by Theorem~\ref{rank 1}(b).  Hence,
the additional relation $z=z'$ produces the desired presentation. \qedd



\subsection{Orthogonal groups}
\label{Orthogonal groups}

\begin{Theorem}
\label{Bounded On}
All perfect central extensions of simple orthogonal  groups of dimension
$N$ over $\F_q$  have
bounded presentations of length
   $O(\log N +\log q)$.
\end{Theorem}
\Proof
We may assume that   $N>8$ and  the center of $\hat G$ is nontrivial.
By
Lemma~\ref{Properties of universal central extensions}, $q$ is odd and
there
are various quotient groups to consider.
\subsubsection{Factoring
out the spin involution}
\label{Factoring out the spin
involution}
For each simple orthogonal group $G$  we need to factor
out $\langle s\rangle $  from $\hat G$, where
$s$ is the spin involution in
Lemma~\ref{Properties of universal central
extensions}(7).

Assume that $G$ does not have type $D_n$.
By Lemma~\ref{Properties of universal central extensions}(7), $s$ is the
central involution   in $L_{\alpha_n} \cong\SL(2,q)$ or
$ \SL(2,q^2)$ (where the central quotient groups are
the orthogonal groups
$\PO(3,q)$ or
$\PO^-(4,q)$).
By  Theorem~\ref{rank 1}(b),  $s$ has length $O(\log q)$ in our
generators.  Thus, the additional relation $s=1$ produces the desired
presentation.

If  $G$ has type $D_n$ we use
$\langle L_{\alpha_{n-1}},L_{\alpha_{n }}\rangle \cong \SL(2,q)\times  \SL(2,q)$.
This time  $s$ is the central involution lying in neither $\SL(2,q)$
factor.  Once again $s$ has length $O(\log q)$ in our
generators, and  the additional relation $s=1$ produces the desired
presentation.

By Lemma~\ref{Properties of universal central extensions},
 we have now provided a short  bounded presentation for $\hat G/\langle
s\rangle $ in all cases; this proves Theorem~\ref{Bounded On} when $G$ is
$\Omega(2n+1,q)$, $\PO^-(2m,q)$ with $m$ even, or
$\PO^\e(2m,q)$ with $m$ odd and $q\equiv -\e 1$ (mod 4),
where $\e=\pm$.
\subsubsection{$D_n(q)$ with $n$ even}
\label{$D_n(q)$ with $n$ even}
   Here  we merely repeat the
argument in Case 2 of the proof of
Theorem~\ref{universal},
this time using  the presentation for
$G_1=\SL(n,q) $ provided by Lemma~\ref{-1}(i).  The additional relation
$z=1$ produces the desired presentation.

At this point we have a presentation for $G= \PO^+(2n,q)$, as well
 as  for
its universal cover $\hat G $. The center of  $\hat G$ is
elementary abelian of order  4.
Starting with $G$, use Proposition~\ref{Alex'}
twice (with $p=2$) in order to obtain a presentation for $\hat G$  in
which we have words of  length  $O(\log n +\log q)$   for the
generators
$a,b$ of $Z(\hat G )$.  Then each group $\hat G/\langle a\rangle $, 
$\hat G/\langle b\rangle $,
$\hat G/\langle a,b\rangle $ can be obtained by adding one or two   
relations of length $O(\log n +\log q)$. 
(Of course, one of these groups is $\hat G/\langle s\rangle $,
which was dealt with in Section~\ref{Factoring out the spin
involution}.)

\smallskip
Note that earlier   we   passed from
the universal central extension of a simple group  to quotients of  it.
The above argument went in the  opposite direction:  from a simple group
to all perfect central extensions,  in the case where  the
center of the universal central extension was bounded.  This same idea
already occurred in Corollary~\ref{extensions of alternating}.
\subsubsection{Factoring out $-1$}
It remains
to consider  the   orthogonal group   $\PO^\e(2m,q)$ where
$m$  is odd and
$q\equiv \e 1$ (mod 4). 

It is convenient to rename $G$ as  
 $\hat G/\langle s\rangle =\Omega^\e(2m,q)$.  We will obtain a new 
  presentation for $G$ in order  to   factor out
$\langle -1\rangle $.

\smallskip
   {\noindent\bf Case 1: $G=\Omega^-(2m,q)$ with $m$ and $q$ odd.}
Here the rank is $n=m-1$,  which is even.
Use
Lemma~\ref{-1}(i) to obtain a   bounded presentation
of length $O(\log n +\log q)$ for
$G_1=\SL(n ,q) $  such that  its central involution  
 $z $ has length
$O(\log n +\log q)$ in the generators.

Use this presentation for $G_1$ in
   Case 1 of the proof of Theorem~\ref{universal} in order to obtain
a presentation for  $\hat G$ of length $O(\log n + \log q)$.
Pass modulo the spin involution as in Section~\ref{Factoring out the
spin involution}  in order to obtain a presentation for $G$.
 We still  need  a short bounded presentation for $G/\langle -1\rangle $.

Let $e_1,\dots,e_m,f_1,\dots,f_m$ be a basis
of the    $\F_q$-space underlying $G$,
where $e_1,\dots,$ $e_n,f_1,\dots,f_n$ is a hyperbolic basis of its span
  and is  perpendicular to the anisotropic 2-space $\langle e_m, f_m
\rangle . $     We can choose these bases so that  
$L_{\alpha_i}\!$ has support
$\langle e_i, \hspace{-1pt}  e_{i+1},\hspace{-1pt}  f_i,\hspace{-1pt} f_{i+1} \rangle$ for each $i \le  n$.
Then $G_1=\langle L_{\alpha_1} ,\dots ,L_{\alpha_{n-1}} \rangle\cong
\SL(n,q)$ has support $\langle e_1,\dots,e_{n},$ $f_1,\dots,f_{n}
\rangle$,  and $z$  is  $-1$ on this subspace.
Note  that $\langle e_1,\dots,e_{n},f_1,\dots,f_{n}
\rangle^\perp =   \langle e_m, f_m
\rangle$  
has type $\Omega ^-(2,q)$, and  
$-1\in \GL( \langle e_m, f_m
\rangle)$ lies in  $\Omega ^-(2,q)$  since $q\equiv 3$ (mod 4).

The group $ \langle L_{\alpha_{n-1}}, L_{\alpha_n}\rangle  \cong
\Omega^-(6,q)$ has an element $z'$ 
with support  
$\langle e_m, f_m
\rangle$ and acting as
$-1$ on this 2-space.  Use Proposition~\ref{short words} to write $z'$ as
a   word of length $O(\log q)$  in our generators.  Now add the relation
$z=z'$ in order to obtain  the desired
presentation for $G/\langle -1\rangle   =\PO^-(2m,q)$.
\smallskip

   {\noindent\bf Case 2: $G=\Omega^+ (2n,q)$ with $n$ and $q$ odd.}
Here the rank is $n$.
Let 
$e_1,\dots,$ $e_n,f_1,\dots,f_n$ be a hyperbolic basis of $V$ such   
$L_{ \alpha_i}$ has support
$\langle e_i,e_{i+1}, f_i,f_{i+1} \rangle$ for each $i \le  n$.
Let $G_1=\SL(n,q)  <G$ preserve each of the subspaces 
$\langle e_1,\dots, e_n\rangle$ and
$\langle f_1,\dots, f_n\rangle$.
Let $H=\SL(n-1,q)$ preserve 
$\langle e_1,\dots, e_{n-1}\rangle$ and
$\langle f_1,\dots, f_{n-1}\rangle$, inducing 1 on 
$\langle e_1,\dots, e_{n-1} ,f_1,\dots, f_{n-1}\rangle ^\perp
= \langle e_n, f_n\rangle$.

 Use this subgroup  $H$ in Lemma~\ref{-1}(ii) in order
to obtain a presentation for $G_1$ such that the central
involution
$z$ of $H$ has length $O(\log n + \log q)$.
Use this presentation for $G_1$ in
   Case 2 of the proof of Theorem~\ref{universal} in order to obtain
a presentation for  $\hat G$ of length $O(\log n + \log q)$.
Pass modulo the spin involution as in Section~\ref{Factoring out the
spin involution}  in order to obtain a presentation for $G$.

The element 
$-1\in \GL( \langle e_n, f_n
\rangle)$ lies in  $\Omega ^+(2,q)$ since $q\equiv 1$ (mod 4).
This time  $ \langle  L_{\alpha_{n-1}},
  L_{\alpha_n}\rangle  \cong \Omega^+(4,q)$ has an element $z'$  
 with support 
$\langle e_n, f_n
\rangle$ and acting as
$-1$ on this 2-space.  
Use Theorem~\ref{rank 1}(b) to write  
$z'  $   as a   word
of length $O(\log q)$ in
our generators.  Now add the relation
$z=z'$ in order to obtain the desired
presentation for $G/\langle -1\rangle  = \PO^+ (2n,q)$.   \qedd

\smallskip

\subsection{Unitary groups}
\label{Unitary groups}
\label{unitary}

We  handle unitary groups as in  Section~\ref{Generic case}, once again
choosing a presentation  more carefully.

\begin{Theorem}
\label{Bounded SUn}
All groups  $\SU(N,q) /Z, $  where  $Z\le Z(\SU(N,q)),$
have
bounded presentations of length
   $O(\log N +\log q)$.
\end{Theorem}

\Proof
As above we rename $G$  as  $\SU(N,q) $, and let 
  $V$ be the underlying  $\F_{q^2}$-space.
Let $n=[N/2]$, and let the vectors
$e_1,\dots,e_n,$ $f_1,\dots,f_n$ be a  hyperbolic basis for the subspace
they span;  moreover,  let these be perpendicular to a final basis vector
$v$ if
$N=\dim V$ is odd.  We will focus on the case where $N = 2n+1$ is odd, the
even case possibly being slightly easier.
 We write matrices with our basis ordered
$e_1,\dots,e_n,f_1,\dots,f_n,v$.

We can choose the basis so that   $L_{ \alpha_i}$ has support
$\langle   e_i,e_{i+1}, f_i,f_{i+1}\rangle $ for $1\le i\le n-1$,  and
$\langle   e_{n },f_{n}, v \rangle $ when $i=n$.  Then the support of $G_1
=\SL(n,q^2)$ is
   $\langle  e_1,\dots,e_n,$ $f_1,\dots,f_n\rangle $,
while the support of $G_3 \cong\SU(5,q)$ is
$\langle    e_{n-1 },e_{n },f_{n-1}, f_{n}, v \rangle $.

In Case 1 of the proof of Theorem~\ref{universal}
  we used the presentation for $G_1 $
appearing in Theorem~\ref{SL(n,q)}. This time we use the presentation
for $ G_1$ deduced from it via Lemma~\ref{MARTIN'S  LEMMA}, using
$m=n-1 $, $Z=\langle \e I\rangle $ and  $H=\langle \e\rangle \le 
\F_{q^2}^*$.  Here 
$H^m=H^{n-1}$ is embedded in $\SL(n,q)$ as the diagonal matrices
$\diag(\e_1,\dots, \e_n,\e_1,\dots, \e_n,1)$
with
$\e_i \in H$ and
$\prod_i\e_i=1$ (note that $\bar\e^{-1}=\e$:  this matrix is, indeed, an
isometry);  while
$S_m=S_{n-1}$ is embedded  as all matrices
$$
   \begin{pmatrix}
         \begin{matrix}
        \pi &  \,\,\,\,0 & \\
         0 &  \pm1
     \end{matrix}  & O  &O \\
        O &  \begin{matrix}
       \pi^{-1}  &  \,\,\,\,0 & \\
         0 & \pm1
     \end{matrix}  &O  \\
       O&O&   1
     \end{pmatrix}
    $$
   using permutation matrices $\pi$ of sign $\pm 1$.
Lemma~\ref{MARTIN'S  LEMMA} produces a new
generator
$d$ in  a slightly  modified short bounded  presentation for $G_1$; the
matrix for
$d$ is
$\diag(\e,\dots,\e,\e^{2-n},1,\e,\dots,\e,\e^{2-n},1,1  )$
(where we have reordered the  coordinates used in that lemma).

Since $\e\in Z(G)$ we  have $\e^{2n+1}=1$, and  there is an element $d'\in
G_3$ with matrix
$$\diag(1,\dots,1,\e^{n-1},\e,1,\dots,1,\e^{n-1},\e,\e  )$$
such that $d d' = \e I$.  Express $d'$ as a word of length $O(\log q)$ in
the generators for $G_3$  by  using Proposition~\ref{short words}.
The additional relation $dd'=1$ produces the desired presentation
for the quotient group
$G/\langle \e I\rangle =G/Z$.~\qedd

\subsection{Perfect central extensions}

The Curtis-Steinberg-Tits presentation produces the universal central
extension of a group of Lie type.  Once we obtained a version of this
presentation, we  had to factor out all or part of the center.  In the
process we  proved the following additional

\begin{Theorem}
\label{central extensions}
Each perfect central extension of a finite  simple group  $G$  of rank
$n$ over   $\F_q $ has  a bounded presentation  of length
   $O(\log n +\log q)$.
\end{Theorem}
\Proof
In the preceding parts of   Section~\ref{Theorems A and B} we  proved
this theorem for the various classical groups.  We dealt with the
alternating groups in Corollary~\ref{extensions of alternating}.
    \qedd

\section{Theorems~\ref{B} and~\ref{B'}}
   \label{Theorems B and B'}

In this section, we prove Theorems~\ref{B}  and~\ref{B'} (Holt's
Conjecture
for simple groups).   We recall the following
well-known
observation (cf. \cite[Lemma~1.1]{Holt}): 

\begin{Lemma}
If  a  finite  group   $G$ has
a presentation $F/R$ with
$F$ free and $R$ normally generated by $r$
elements$,$ then
$\dim_{\,\Fb_p \!} H^2(G,M) \le r \dim _{\,\Fb_p \!}M$ for any
$\F_pG$-module
$M$.

\end{Lemma}

  \Proof
  A  short exact sequence
$$
1
\rightarrow M \rightarrow E \rightarrow F/R \rightarrow 1
$$
  is
determined (up to equivalence) by a $G$-homomorphism
   $R\to M$.
 Since
$R$ can be generated by $r$ elements (as a  normal
subgroup), any
homomorphism is determined by  the images of the $r$
generators.
Thus,  there are most $|M|^r$ such homomorphisms and
therefore at
most that many inequivalent  extensions.   Consequently,
$\dim
_{\,\Fb_p \!}H^2(G,M) \le r \dim_{\,\Fb_p \!} M$.
\qedd

\medskip

\noindent \emph{Proof of} Theorem~\ref{B'}:
By  the
lemma, Holt's Conjecture follows from Corollary~\ref{A'}   except for the
case  $^2G_2(3^{2e+1} )$, which we now handle.

Note that, if $D$ is
a subgroup of $G$ that contains a Sylow
$p$-subgroup of
$G$, then the
restriction map   $H^2(G,M) \rightarrow H^2(D,M)$
is injective  
\cite[p.~91]{Gru}. 

In order to use this, we consider the cases
$p=2,$  $ p=3$ and $p > 3$
separately.

If $p > 3$, a Sylow
$p$-subgroup   $D$ of $G$ is cyclic, and hence  can
be presented with
one generator and one relation.
Applying the preceding lemma to  $D$
yields the result  for    $G$
with constant $1$.

A Sylow
$2$-subgroup of $G$ is elementary abelian of order $8$,
and hence has
a  presentation
  with $3$ generators and $6$ relations.

Finally,
for $p=3$, by Proposition~\ref{Ree Borel}
  a Borel subgroup of $G$  has
a  presentation with a bounded number of
relations, whence the  above
lemma produces a bound in characteristic $3$ as  well.~\qedd

\medskip
We will provide a different    argument for  this Ree case in 
\cite{GKKL}.

Working in the profinite category, Theorems~\ref{B}
and~\ref{B'} are equivalent. 
Namely, in  \cite{Lub2} there is  a formula
for the minimal number
$\hat  r(G)$ of  relations needed for a
profinite presentation   of a
finite group
$G$:
\begin{equation}
\label{r(G)}
\hat r(G)= \sup_p \sup_M \Big(
\Big\lceil  \frac{\dim
H^2 (G, M) - \dim H^1 (G, M)}{\dim M} \Big\rceil  + d(G) - \xi_M
\Big) ,
\end{equation}
    where $d(G)$ is the minimum number of  generators for
$G$,
$p$ runs over all primes, $M$ runs over all irreducible
${\F}_p [G]$-modules, and $\xi_M = 0$ if $M$ is the trivial
module and 1 if not.
Note that $\dim H^1(G,M)\le d(G) \dim M$: every element of
$Z^1(G,M)$ is a map $G\to M$ that is completely   determined on a
generating set of
$G$.  Hence, \eqn{r(G)} implies that
\begin{equation}
\label{hat r(G) inequalities}
h (G) \leq \hat r(G) \leq   h(G)+ d(G) ,  
\end{equation}
where
$$
h(G) :=\mathop{\sup}\limits_p \mathop{\sup}\limits_M \frac{\dim
H^2(G,M)}{\dim M} .
$$
Since $d(G)=2$ for every finite simple group,  \eqn{hat r(G)
inequalities} proves our assertion that Theorems~\ref{B} and~\ref{B'} are
equivalent. 

\section{Concluding remarks}
\label{Concluding remarks}

\noindent
1.
One of the purposes of \cite{BGKLP} was to provide presentations for use
in Computational Group Theory.
This has turned out to be essential, for example  in \cite{KS1,KS2}.
When used  in these references, those presentations led to efficient
algorithms  to test whether or not a given matrix group actually  is
isomorphic to a specific simple group.  Short presentations are also
essential in \cite{LG,KS2} for gluing together presentations in a normal
series (essentially a chief series) in order to obtain a presentation for
a given matrix group.

We expect that   versions  of  many of the
presentations in the present paper
or  \cite{GKKL2}
will have both practical and theoretical
significance in Computational Group Theory.

\medskip
\noindent
2. As mentioned in the Introduction,
an elementary counting argument shows that our $O(\log n +\log q)$
bound in Theorem~\ref{A} is  optimal in terms of $n$ and $q$.  However, in
terms  of $n,$  $p$ and $e$ (where $q=p^e)$ one might hope for presentations
providing a slightly better bound:
$O(\log n +\log p  +\log e)$.  We have no idea whether  such
presentations exist  even for $\PSL(2,p^e)$. 

\medskip
\noindent
3.
A standard topological interpretation of Corollary~\ref{A'} states that all
finite simple groups (except perhaps $^2G_2(q)$) are fundamental
groups of 2--dimensional CW-complexes having a bounded number of cells. 

\medskip
\noindent
4.
As suggested in the Introduction, it is not difficult to check that
fewer than 1000 relations are required in Theorem~\ref{A}.   
It would be
interesting to have far better constants in all of the theorems. 
On the other hand, it would also be interesting to have presentations in
Theorem~\ref{A} having only 2 generators, even if the number of relations
grew somewhat (while remaining bounded); compare \cite{GKKL2}.

Earlier we mentioned Wilson's conjecture  that the universal central
extension of every finite simple  group has a presentation with 2
generators and 2 relations
\cite{Wil}.  In view of a standard property of  the Schur multiplier
$M(G)$  \cite{Sch}, this would be optimal for a bounded presentation
and  would imply that  every nonabelian
finite simple group $G$ has a presentation  with 2 generators
   and $2+d(M(G))$ relations.   The only  infinite family  for 
which the conjecture is presently proven consists of the groups 
$\PSL(2,p)$  \cite{CR2}.
Wilson essentially proved the conjecture for $\Sz(q)$,
$q>8$, in the category of  profinite presentations.

\medskip

\noindent
5.  Finally,  
in a similar vein, we note that there
does not appear to be a  known  presentation of
$\SL(2,p)$, when  $p$ is  prime,  having  $k$
generators and $k+1$ relations for some $k$,
and also having length $O(\log p)$.  
In Section~\ref{Congruence Subgroup Property}   we   already  observed
that the presentation of Sunday \cite{Sun}   is not this short.

\medskip
\noindent
\emph{Acknowledgement}:  We are grateful to the referee for
many helpful comments.
We also thank  Ravi Ramakrishna    and   Joel Rosenberg
 for assistance with  Lemma~\ref{a,b}, and
  James Wilson for     checking our 
computations in Table~\ref{139-irr} using {{\sf GAP}}. 

\appendix

\setcounter{table}{0}
\renewcommand{\thetable}{\thesection.\arabic{table}}

\section{Field lemma}
\label{Field lemma} 
\noindent
\emph{Proof  of} Lemma~\ref{a,b}.
Assume that $ab \not = 0$.
Then \eqn{a,b equation}  is equivalent to
\begin{equation}
\label{a,b equation again}
\bar a^2 b + \bar b^2 a = ab,
\quad
\bar aa + \bar bb =1.
\end{equation}
Taking  norms of both sides of the first equation and
inserting the second  yields
\begin {align*}
a^2 \bar a^2 b \bar b + \bar a^3 b^3 + a^3 \bar b^3 + a \bar a b^2 \bar
b^2  &= a\bar a b\bar b(\bar aa + \bar bb)
\\
&=  a^2 \bar a^2 b \bar b  + a \bar a b^2 \bar b^2 .
\end {align*}
Therefore, $(a\bar b)^3 \in \F_q \theta $
with  $\theta ^q= - \theta \ne0$.

Consequently, $a\bar b= k\theta \omega$ with $k\in \F_q,$ $ \omega\in F$
and  $\omega^3=1=\omega\bar \omega$, where $\omega=1$ if
$q\not\equiv-1$ (mod 3).
Substituting $\bar b = k \theta  \omega /a $ into \eqn{a,b equation again}
gives
$$
- \bar a^2 k \theta   \omega^{-1}/\bar a + a (k  \theta   \omega
/a)^2  = - a
k\theta  \omega^{-1}/\bar a,\quad
a\bar a - k^2 \theta ^2/ (a\bar a) =1,
$$
which imply that
$$
k \theta  = a (\bar a - a/\bar a), \quad k^2 \theta ^2 = a\bar a(1-a\bar a
).
$$
Then $-k^2 \theta ^2=k \theta k \bar\theta = a\bar a (\bar a - a/\bar a) (
a -
\bar a/ a)$, which simplifies to
\begin{equation}
\label{the curve}
2(a\bar a)^2= a^3 + \bar a^3.
\end{equation}

Any solution of \eqn{the curve} yields
$(3,q+1)$ solutions of
\eqn{a,b equation again} corresponding to the different possibilities
for $\omega$.   For example, 
using $k : = a (\bar a - a/\bar a)/ \theta$ and $\bar b:=k\theta
\omega/a$ with
$\omega^3=1$, we find that $k\in \F_q$ since 
$$
\bar k= - \bar a (  a - \bar a/ a)/ \theta
=   a (  \bar a -  a/\bar a)/ \theta =k
$$
by \eqn{the curve}; and hence  
\begin{align*}
\bar a ^2b  +a \bar b ^2 - ab& =
b \cdot  \bar a   \bar b / \omega + a  \bar b^2
\\
&=
 \bar b
[ \bar a  \bar \omega\cdot  (a- \bar a /a) \bar \omega  + a( \bar a -a/
\bar a ) \omega]
\\
&=
 \bar b  \omega[ 2 a^2  \bar a ^2 -  \bar a ^3 - a^3]/ a \bar a   =0.
\end{align*}

  If $q$ is odd then $\F_{q^2} =\F_q[\theta]$, so we can
write $a= x + y \theta$ with $x,y\in \F_q$ and  (\ref{the curve})
becomes
$$
2(x^2 - y^2 \theta^2)^2 = 2(x^3 + 3xy^2 \theta^2).
$$
The left side is homogeneous of degree $4$ and the right side is
homogeneous of degree $3$.
This implies that this equation defines a  curve is of genus $0$, and a
parametrization can be built using the fact that  each line $x=ty$ intersects
the curve at
$(0,0)$ and at most one other point.
Letting $x=t f(t)$ and $y= f(t)$ we obtain the equation
$$
2f(t)^4(t^2 - \theta^2)^2 = 2 f(t)^3 t (t^2 + 3 \theta^2),
$$
which has solutions
$$
f(t) = t\frac{(t^2+ 3\theta ^2 )}{(t^2-\theta ^2 )^2}
\quad
\mbox{and}
\quad
f(t) = 0.
$$
Therefore
$$
a = t \frac{(t^2+ 3\theta ^2 )}{(t^2-\theta ^2 )^2} (t +  \theta )
$$
(with a related formula for $b$) is a parametrization of the curve.


The map $t  \mapsto a\bar a$ is a rational
function over
$\F_q$  of degree
$6$, therefore $a\bar a$ does not generate $\F_q$  for at most
$\sum_{k|e,k<e} 6 p^k$ values of $t$.  Moreover, $a^3$  (or, equivalently,
$a^{2q-1}$) is in
$\F_q$ for at most 5 choices $t$, because
$$
a^3 = t^3 \frac{(t^2+ 3\theta ^2 )^3}{(t^2-\theta ^2 )^6}
 \left[(t^2 + 3 \theta^2) +  \theta(3t^2 + \theta^2 )\right],
$$
and the equation
$t^3 (t^2+ 3\theta ^2 )^3(3t^2 + \theta^2 )=0$ has at most $5$ distinct solutions.
Hence,  the number of $t$ such that $a ^3 \notin \F_q$ and
$a\bar a $   generates
$\F_q$ is at  least
$$
q-5 -6 \sum_{k|e,k<e} p^{k} \ge  q-5 -6(e-1) q^{1/2}.
$$
Hence there exist at least 4 such    $a$ whenever $q\ne 3,5,9, 25.$ 
Since $a\bar a= a^{2q-1}\bar a^{2q-1}$
generates  $\F_q $ over $\F_p$,  it follows that  $a^{2q-1}$ generates $F$
over $\F_p$.
Direct verification shows that \eqn{a,b equation again} has a solution with
the desired properties for $q=9$ and $q=25$.

If $q$ is even then \eqn{the curve} is  a union of   $(3,q+1)$ lines:
$\bar a=a$, together with $\bar a=\omega a$ and $\bar a=\bar\omega a$
if  some $\omega\notin \F_q$  satisfies  $\omega^3=1$.
If $a$ is any generator of
$\F_q  ^*$ then  the pair $(a,1-a)$ is solution to \eqn{a,b equation
again} such that  
 $\F_q=\F_p[a\bar
a]=\F_p[a^{2q-1}]$.
\qedd

\medskip

Our original approach to Lemma~\ref{a,b} was very different.
First we observed that if $(3 ,q+1)=1$ then \eqn{a,b equation}  is equivalent
to
$$
A+B =1 ,  \quad A^{q+1}+B^{q+1}=1
$$
and hence, after substitution of $B=1-A$ in the second equation, to
$$
A + \bar A = 2 A \bar A,
$$
and this has $q$ solutions in $F$, because it defines a curve of genus
0 over $\F_q$  (if we view $A=x+ y \theta$  with $x,y \in \F_q$).
An easy counting argument shows that unless $q$ is even, there is always a
solution which generates
$F$ over $\F_p$.

If $3|q+1$ then  $F=\F_q[\omega]$, where $\omega^2 + \omega +1=0$.
This allows us to write
$
a=x + y \omega 
$,
$b=u+ v \omega
$
where $x,y,u,v \in \F_q$. Substituting in \eqn{a,b equation},  we obtain
a system of
$3$ equations in $4$ variables over $\F_q $, of degrees $3$, $3$ and $2$.
Since the coefficients in these equations turn out to be  integers these
equations define a curve
$C$ (or several curves)
over $\Q$. Using the bounds on the degrees of the equations one can
obtain an upper bound
for the genus $g$ of the curve.
It is known that  the number of points of $C$  over $\F_q$ is at least
$q - 2g \sqrt{q}$, which for large enough  $q$ is more than the number
of elements in
$\F_q$ which do not generate this field.

We started studying the above system   in order to obtain a nice
bound for the genus $g$
and  to  find how large $q$ has to be for the above argument to work.
We expected to get a curve of genus $9$ or $10$ but after doing
computations we
were surprised to find  that actually we had  $3$ curves of genus
$0$, which allowed
us to parametrize all solutions of the system.
(N.B. If we change one of the constants on the right in  \eqn{a,b
equation}  to a generic element in $\F_p$, we  indeed obtain a curve of
genus $9$.)
 
\section{Suzuki triples}
\label{Suzuki triples}
\noindent
\emph{Proof  of} Lemma~\ref{A Sz}.
Every vector in $\Z^4$ that cannot be reduced using   (a)
and (b) must have one of the following 25 types for some integers $a,b\ge0$.
(For,  by   Example~\ref{Original BCRW triples},    the first and  third
coordinates must be one of
$(-a,0)$, $(0,-a)$, $(a,a)$, $(a,a+1)$, $(a+1,a)$,
and similarly for the second and fourth.)

{\small
$$
\begin{array}{lllllllll}

\quad{\rm Case~ 1}&\quad{\rm Case~ 2} &\quad{\rm Case~ 3}&\quad {\rm
Case~ 4} &\quad{\rm Case~ 5}
\\
(-a,-b,0,0)&(-a,0,0,-b) &(a,-b,a,0)&(0,b,-a,b)  &(a,b,a,b)
\\
 &&(a,-b,a+1,0)& (0,b+1,-a,b)&(a,b,a+1,b)

\\
(0,0,-a,-b)&(0,-b,-a,0)&(a+1,-b,a,0)&(0,b,-a,b+1)&(a+1,b,a,b)

\\
&&& &(a,b,a,b+1)
\\
&&(a,0,a,-b)
 &(-a,b,0,b)&(a,b,a+1,b+1)
\\
&&
  (a+1,0,a,-b)&(-a,b,0,b+1)&(a+1,b,a,b+1)
\\
&& (a,0,a+1,-b) &(-a,b+1,0,b) &(a,b+1,a,b)
\\
&&&&(a,b+1,a+1,b)
\\
&&&&(a+1,b+1,a,b)
\end{array}
$$
}

Since $\mathcal{M}_\Sz$  is   invariant under the permutation  $(1,3) (2,4)$ of
coordinates,  there are only  5   cases that need to be considered
on this  list:
 the ones on the second row above,  in the last 3 columns decorated with
some
$+1$'s as seen on the remainder of the last three columns.
That is, we only need to consider the five cases
$v=(-a,-b,0,0), (-a,0,0,-b),
(a+\alpha ,b+\beta ,a+\gamma, b+\delta)$,
$(a+\alpha , -b,a+\gamma,0)$,
$(0,b +\beta , -a ,b+\delta)$
 with $a,b\ge0$ and
$\alpha ,\beta ,\gamma,\delta \in\{0,1\}$.
We will present the second and third of these cases, leaving the
remaining  high school algebra to the reader.

\medskip
\noindent {\bf Case 2:}
Let   $v=(-a,0,0,-b)$ with $a,b\ge0$.
Triple (c) has  three translates containing $v$:
$$  
\begin{array}{llllllll}
1 & (-a,0,0,-b)&& (1-a,-1,0,-b )&&(-a,0,1,-1-b )
\\
2& (-1-a,1,0,-b)&& (-a,0,0,-b)&& (-1-a ,1,-1-b)
\\
3&( -a,0,-1,1-b )&& ( 1-a,-1,-1, 1 -b)&& ( -a,0,0,-b)
\end{array}
$$
Triple (d) has  three translates containing $v$:
$$ \hspace{-9pt}
\begin{array}{llllllll}
4 & (-a,0,0,-b)& (-2-a,1,0,-b )&(-a,0,-2, 1-b )
\\
5& (2-a,-1,0,-b)& (-a,0,0,-b)& (2-a,-1 ,-2, 1-b)
\\
6&(  -a,0,2,-1-b )& ( -2-a, 1,2,-1-b)& ( -a,0,0,-b)
\end{array}
$$

 Here $ v $ is  (c)-reducible  via 1 iff

$a^2+2b^2 >  (1-a)^2+2+2b^2$ and

$a^2+2b^2 > a^2+1+2(1-b)^2$,

\noindent
hence iff  $2a > 2 $ and $ 0>2+4b$, which never occurs.

\medskip
Also $ v $ is  (c)-reducible  via 2 iff

$a^2+2b^2 >  (-1-a)^2+2+2b^2$ and

$a^2+2b^2 > (-1-a)^2+2+1+ 2(-1-b) ^2$,

\noindent
hence iff  $0>2a+3 $ and $0>2a+4b+6$, which never occurs.

\medskip
Also $ v $ is  (c)-reducible  via  3 iff

$a^2+2b^2 >  a^2+1+2( 1-b)^2$ and

$a^2+2b^2 >  ( 1-a)^2+2+1+2(1-b)^2$,

\noindent
hence iff  $ 4b>3$ and $ 2a+4b> 6 $.

\medskip

Thus, {\em the only way $v$ can be  {\rm(c)}-irreducible is to have
$4b\le 3$  or $2a+4b\le 4$.
Hence $b\le 1$.}

\medskip

Next,  $ x $ is  (d)-reducible  via 4 iff

$a^2+2b^2 > (-2-a)^2+2+2 b^2$ and

$a^2+2b^2 > a^2+4 +2(1-b)^2$,

\noindent
hence iff  $ 0>4a+6$ and $4b>6$, which never occurs.

\medskip
Also $ v $ is  (d)-reducible  via 5 iff

$a^2+2b^2 >  ( 2-a)^2+2+2b^2$ and

$a^2+2b^2 >   ( 2-a)^2+2+4+2(1-b)^2$,

\noindent
hence iff  $ 4a>6$ and $ 4a+4b>12 $.
\medskip

Also $ v $ is  (d)-reducible  via 6 iff

$a^2+2b^2 >  a ^2+4+2(-1-b)^2$  and

$a^2+2b^2 >   (-2-a)^2+2+4+2(-1-b)^2$,

\noindent
hence iff  $0>4b+6 $ and $ 0>4a+4b+12$, which never occurs.

\medskip

Thus, {\em the only way $v$ can be  {\rm(d)}-irreducible is to have
$4a\le 6$  or  $4a + 4b \le  12$.
Hence $a\le 3$.}

Consequently, all  irreducible $v$ of the specified sort have
$a\le 3$,
$b\le 1$ and hence
$|v|^2\le  3^2+2
$.

\medskip

\noindent {\bf Case 3:}
The  9 types in the last column have the form $v=(a+\alpha ,b+\beta ,a+\gamma, b+\delta)$
 with $a,b\ge0$ and
$\alpha ,\beta ,\gamma,\delta \in\{0,1\}$.
Triple (c) has  three translates containing $v$:
$$\hspace{-2pt}
\begin{array}{llllllll}
\scriptstyle\hspace{-4pt} 1\hspace{-6pt} &
\scriptstyle (a+\alpha ,b+\beta ,a+\gamma, b+\delta)&
\scriptstyle (1+a+\alpha ,-1+b+\beta ,a+\gamma, b+\delta)&
\scriptstyle (a+\alpha ,b+\beta ,1+a+\gamma, -1+b+\delta)
\\
\scriptstyle \hspace{-4pt}2\hspace{-6pt}&
\scriptstyle (-1+a+\alpha ,1+b+\beta ,a+\gamma, b+\delta)&
\scriptstyle (a+\alpha ,b+\beta ,a+\gamma, b+\delta)&
\scriptstyle \hspace{-21pt}(-1+a+\alpha ,1+b+\beta ,1+a+\gamma, -1+b+\delta)%
\\
\scriptstyle \hspace{-4pt}3\hspace{-6pt}&
\scriptstyle ( a+\alpha , b+\beta ,-1+a+\gamma, 1+b+\delta)&
\scriptstyle ( 1+a+\alpha ,-1+b+\beta ,-1+a+\gamma,1+ b+\delta)&
\scriptstyle (a+\alpha ,b+\beta ,a+\gamma, b+\delta)
\end{array}
$$ 

\noindent
Triple (d) has  three translates containing $v$:
$$\hspace{-2pt}
\begin{array}{llllllll}
\scriptstyle \hspace{-4pt}4 \hspace{-6pt}&
\scriptstyle (a+\alpha ,b+\beta ,a+\gamma, b+\delta)&
\scriptstyle (-2+a+\alpha , 1+b+\beta ,a+\gamma, b+\delta)&
\scriptstyle  (a+\alpha ,b+\beta ,-2+a+\gamma, 1+b+\delta)
\\
\scriptstyle \hspace{-4pt}5\hspace{-6pt} &
\scriptstyle (2+a+\alpha ,-1+b+\beta ,a+\gamma, b+\delta)&
\scriptstyle (a+\alpha ,b+\beta ,a+\gamma, b+\delta)&
\scriptstyle\hspace{-21pt} (2+a+\alpha ,-1+b+\beta ,-2+a+\gamma,  1+b+\delta)%
\\
\scriptstyle \hspace{-4pt}6\hspace{-6pt} &
\scriptstyle ( a+\alpha , b+\beta ,2+a+\gamma, -1+b+\delta)&
\scriptstyle (-2+ a+\alpha ,1+b+\beta ,2+a+\gamma,-1+ b+\delta)&
\scriptstyle (a+\alpha ,b+\beta , a+\gamma,  b+\delta)
\end{array}
$$

Here  $ v$ is  (c)-reducible  via 1 iff
\begin{align*}
\scriptstyle (a+\alpha)^2+2(b+\beta)^2+(a+\gamma)^2+ 2(b+\delta)^2 &\scriptstyle >
(1+ a+\alpha)^2+2(-1+ b+\beta)^2+(a+\gamma)^2+ 2(b+\delta)^2
\quad \scriptstyle \mbox{\scriptsize and}
\\
\scriptstyle (a+\alpha)^2+2(b+\beta)^2+(a+\gamma)^2+ 2(b+\delta)^2 &\scriptstyle >
(a+\alpha)^2+2(b+\beta)^2+(1+ a+\gamma)^2+ 2(-1+ b+\delta)^2,\quad \quad{}
\end{align*}
hence iff  $4(b +\beta)>2(a+\alpha) +3$ and $4(b+ \delta)>2(a+\gamma) +3$.

\medskip
Also $ v$ is  (c)-reducible  via 2 iff
\begin{align*}
\scriptstyle
(a+\alpha)^2+2(b+\beta)^2+(a+\gamma)^2+ 2(b+\delta)^2&\scriptstyle>
(-1+ a+\alpha)^2+2(1+ b+\beta)^2+(a+\gamma)^2+ 2(b+\delta)^2
\quad \mbox{\scriptsize and}
\\
\scriptstyle
(a+\alpha)^2+2(b+\beta)^2+(a+\gamma)^2+ 2(b+\delta)^2&\scriptstyle>
(-1 +a+\alpha)^2+2(1+b+\beta)^2+(1+ a+\gamma)^2+ 2(- 1+ b+\delta)^2,
\end{align*}
hence iff
 $2(a+\alpha)>4(b+\beta) +3$  and $ 0 >2\gamma-2\alpha+ 4\beta-4\delta
+6  $, which never occurs.

\medskip
Also $ v $ is  (c)-reducible  via 3  iff
\begin{align*}
\scriptstyle
(a+\alpha)^2+2(b+\beta)^2+(a+\gamma)^2+ 2(b+\delta)^2&\scriptstyle>
(a+\alpha)^2+2( b+\beta)^2+(-1+a+\gamma)^2+ 2(1+b+\delta)^2
\quad \mbox{\scriptsize and}
\\
\scriptstyle
(a+\alpha)^2+2(b+\beta)^2+(a+\gamma)^2+ 2(b+\delta)^2&\scriptstyle>
( 1+a+\alpha)^2+2(- 1+b+\beta)^2+(-1+ a+\gamma)^2+ 2( 1+ b+\delta)^2,
\end{align*}
 hence iff  $2(a+\gamma)>4(b+\beta) +3$  and $0 >  2\alpha
-2\gamma+4\delta-4\beta+6$,  which never occurs.

Thus, {\em the only way $v$ can be  {\rm(c)}-irreducible is to have
 $4(b +\beta)\le 2(a+\alpha) +3$ or $4(b+ \delta)\le 2(a+\gamma) +3,$
and hence to have $4b \le 2a+  5$.}

\medskip
Next,  $ v $ is  (d)-reducible  via 4  iff
\begin{align*}
\scriptstyle
(a+\alpha)^2+2(b+\beta)^2+(a+\gamma)^2+ 2(b+\delta)^2&\scriptstyle>
( -2 + a+\alpha)^2+2( 1+b+\beta)^2+(a+\gamma)^2+ 2(b+\delta)^2
\quad\mbox{\scriptsize and}
\\
\scriptstyle
(a+\alpha)^2+2(b+\beta)^2+(a+\gamma)^2+ 2(b+\delta)^2&\scriptstyle>
(a+\alpha)^2+2(b+\beta)^2+( -2 + a+\gamma)^2+ 2(1+b+\delta)^2,
\end{align*}
hence iff
$4(a+\alpha)>4(b+\beta) + 6$ and $4(a+\g)>4(b+\delta)+6$.

\medskip

Also $v $ is  (d)-reducible  via 5 iff
$$\scriptstyle
(a+\alpha)^2+2(b+\beta)^2+(a+\gamma)^2+ 2(b+\delta)^2>
(2+a+\alpha)^2+2(-1+b+\beta)^2+(  a+\gamma)^2+ 2( b+\delta)^2
\quad\mbox{\scriptsize  and}
$$
$$\scriptstyle
(a+\alpha)^2+2(b+\beta)^2+(a+\gamma)^2+ 2(b+\delta)^2>
(2+a+\alpha)^2+2( -1+b+\beta)^2+( -2+a+\gamma)^2+ 2( 1+b+\delta)^2,
$$
hence iff $4(b+\beta) >4(a+\alpha)+6$ and
$0>4\alpha-4\gamma -4\beta +4\delta +12$,  which never occurs.

\medskip
Also $ v $ is  (d)-reducible  via 6 iff
\begin{align*}
\scriptstyle
(a+\alpha)^2+2(b+\beta)^2+(a+\gamma)^2+ 2(b+\delta)^2&\scriptstyle>
(a+\alpha)^2+2(b+\beta)^2+(2+a+\gamma)^2+ 2( -1+ b+\delta)^2
\quad\mbox{\scriptsize and}
\\
\scriptstyle
(a+\alpha)^2+2(b+\beta)^2+(a+\gamma)^2+ 2(b+\delta)^2&\scriptstyle>
(-2+a+\alpha)^2+2(1+b+\beta)^2+(2+a+\gamma)^2+ 2(-1+b+\delta)^2,
\end{align*}
hence iff $4(b+\delta)>4(a+\gamma)+6$ and $0 >  -4\alpha + 4\gamma +
4\beta-4\delta+12$, but the latter never occurs.

\medskip
Thus, {\em the only way $v$ can be {\rm(d)}-irreducible is to have
$4(a+\alpha)\le4(b+\beta) + 6$ and $4(a+\g)\le4(b+\delta)+6  ;$
and hence to have $4a \le 4b+  1 0$.}

Consequently, {\em all irreducible $v$ of the specified sort have
$4b \le 2a+  5$  and $4a \le 4b+  1 0,$
and hence satisfy   $0\le a\le
7,0\le b\le 4 $   and $|v|^2\le 2(8^2+2\cdot 5^2)$. }

\medskip
We have now bounded the length of each irreducible vector.  The bounds
are in the following table.
$$
\begin{array}{cccccc}
\rm Case & a\le &  b\le & |v|^2\le \\
\rm 1 & 9 &  6& 153 \\
\rm 2 & 3 &  1 & 11 \\
\rm 3 & 1 & 2 &9 \\
\rm 4 & 1 &2 & 26 \\
\rm 5 & 7 & 4 & 228 \\
\end{array}
$$
Instead of considering all cases even more tediously
than above, we used a computer calculation to obtain all irreducible
vectors.
\qedd

\medskip
For completeness  we list all 139 irreducible  vectors in Table~\ref{139-irr}.
\begin{table}
\caption{The irreducible vectors for $\mathcal{M}_{\Sz}$}
\label{139-irr}
$
\begin{array}{|llllll|}
\hline
\scriptstyle (-9,-6,0,0) &
\scriptstyle (-8,-5,0,0) &
\scriptstyle (-7,-5,0,0) &
\scriptstyle (-7,-4,0,0) &
\scriptstyle (-6,-4,0,0) &
\scriptstyle (-6,-3,0,0)
\\
\scriptstyle (-5,-4,0,0) &
\scriptstyle (-5,-3,0,0) &
\scriptstyle (-5,-2,0,0) &
\scriptstyle (-4,-3,0,0) &
\scriptstyle (-4,-2,0,0) &
\scriptstyle (-4,-1,0,0)
\\
\scriptstyle (-3,-3,0,0) &
\scriptstyle (-3,-2,0,0) &
\scriptstyle (-3,-1,0,0) &
\scriptstyle (-3,0,0,0) &
\scriptstyle (-2,-2,0,0) &
\scriptstyle (-2,-1,0,0)
\\
\scriptstyle (-2,0,0,0) &
\scriptstyle (-2,1,0,0) &
\scriptstyle (-1,-2,0,0) &
\scriptstyle (-1,-1,0,0) &
\scriptstyle (-1,0,0,-1) &
\scriptstyle (-1,0,0,0)
\\
\scriptstyle (-1,0,0,1) &
\scriptstyle (-1,1,0,0) &
\scriptstyle (0,-2,1,0) &
\scriptstyle (0,-1,-1,0) &
\scriptstyle (0,-1,0,0) &
\scriptstyle (0,-1,1,0)
\\
\scriptstyle (0,0,-9,-6) &
\scriptstyle (0,0,-8,-5) &
\scriptstyle (0,0,-7,-5) &
\scriptstyle (0,0,-7,-4) &
\scriptstyle (0,0,-6,-4) &
\scriptstyle (0,0,-6,-3)
\scriptstyle\\
\scriptstyle (0,0,-5,-4) &
\scriptstyle (0,0,-5,-3) &
\scriptstyle (0,0,-5,-2) &
\scriptstyle (0,0,-4,-3) &
\scriptstyle (0,0,-4,-2) &
\scriptstyle (0,0,-4,-1)
\\
\scriptstyle (0,0,-3,-3) &
\scriptstyle (0,0,-3,-2) &
\scriptstyle (0,0,-3,-1) &
\scriptstyle (0,0,-3,0) &
\scriptstyle (0,0,-2,-2) &
\scriptstyle (0,0,-2,-1)
\\
\scriptstyle (0,0,-2,0) &
\scriptstyle (0,0,-2,1) &
\scriptstyle (0,0,-1,-2) &
\scriptstyle (0,0,-1,-1) &
\scriptstyle (0,0,-1,0) &
\scriptstyle (0,0,-1,1)
\\
\scriptstyle (0,0,0,-1) &
\scriptstyle (0,0,0,0) &
\scriptstyle (0,0,0,1) &
\scriptstyle (0,0,1,-1) &
\scriptstyle (0,0,1,0) &
\scriptstyle (0,0,1,1)
\\
\scriptstyle (0,1,-1,0) &
\scriptstyle (0,1,0,0) &
\scriptstyle (0,1,1,0) &
\scriptstyle (0,1,1,1) &
\scriptstyle (0,2,1,1) &
\scriptstyle (1,-1,0,0)
\\
\scriptstyle (1,-1,1,0) &
\scriptstyle (1,0,0,-2) &
\scriptstyle (1,0,0,-1) &
\scriptstyle (1,0,0,0) &
\scriptstyle (1,0,0,1) &
\scriptstyle (1,0,1,-1)
\\
\scriptstyle (1,0,1,0) &
\scriptstyle (1,0,1,1) &
\scriptstyle (1,0,2,-1) &
\scriptstyle (1,0,2,0) &
\scriptstyle (1,0,2,1) &
\scriptstyle (1,1,0,0)
\\
\scriptstyle (1,1,0,1) &
\scriptstyle (1,1,0,2) &
\scriptstyle (1,1,1,0) &
\scriptstyle (1,1,1,1) &
\scriptstyle (1,1,1,2) &
\scriptstyle (1,1,2,0)
\\
\scriptstyle (1,1,2,1) &
\scriptstyle (1,1,2,2) &
\scriptstyle (1,2,1,1) &
\scriptstyle (1,2,2,1) &
\scriptstyle (2,-1,1,0) &
\scriptstyle (2,0,1,0)
\\
\scriptstyle (2,0,1,1) &
\scriptstyle (2,0,2,1) &
\scriptstyle (2,1,1,0) &
\scriptstyle (2,1,1,1) &
\scriptstyle (2,1,1,2) &
\scriptstyle (2,1,2,0)
\\
\scriptstyle (2,1,2,1) &
\scriptstyle (2,1,2,2) &
\scriptstyle (2,1,3,0) &
\scriptstyle (2,1,3,1) &
\scriptstyle (2,1,3,2) &
\scriptstyle (2,2,1,1)
\\
\scriptstyle (2,2,2,1) &
\scriptstyle (2,2,3,1) &
\scriptstyle (2,2,3,2) &
\scriptstyle (2,3,3,2) &
\scriptstyle (3,0,2,1) &
\scriptstyle (3,1,2,1)
\\
\scriptstyle (3,1,2,2) &
\scriptstyle (3,1,3,2) &
\scriptstyle (3,2,2,1) &
\scriptstyle (3,2,2,2) &
\scriptstyle (3,2,2,3) &
\scriptstyle (3,2,3,1)
\\
\scriptstyle (3,2,3,2) &
\scriptstyle (3,2,3,3) &
\scriptstyle (3,2,4,1) &
\scriptstyle (3,2,4,2) &
\scriptstyle (3,2,4,3) &
\scriptstyle (3,3,3,2)
\\
\scriptstyle (3,3,4,2) &
\scriptstyle (4,1,3,2) &
\scriptstyle (4,2,3,2) &
\scriptstyle (4,2,3,3) &
\scriptstyle (4,2,4,3) &
\scriptstyle (4,3,3,2)
\\
\scriptstyle (4,3,4,2) &
\scriptstyle (4,3,5,2) &
\scriptstyle (4,3,5,3) &
\scriptstyle (4,4,5,3) &
\scriptstyle (5,2,4,3) &
\scriptstyle (5,3,4,3)
\\
\scriptstyle (5,3,4,4) &
\scriptstyle (5,3,5,4) &
\scriptstyle (5,4,5,3) &
\scriptstyle (5,4,6,3) &
\scriptstyle (6,3,5,4) &
\scriptstyle (6,5,7,4)
\\
\scriptstyle (7,4,6,5) & & & & &
\\
\hline
\end{array}
$
\end{table}

A more conceptual way to prove Lemma~\ref{A Sz} is following.  The set
$\mathcal{M_U}$ is the product of  two copies of $\mathcal{M}_{{\rm
BCRW}}$. Therefore any irreducible vector for $\mathcal{M_U}$ lies close
to   (at distance at most $1$)   one of $9$ quarter planes $A_i$ (each
quarter plane arises as a Cartesian product of $2$ rays). These quarter
planes are listed in the cases in the beginning of this appendix

The system  $  \mathcal{M_{-V}}$ also can be transformed  
to  the product of  two copies of
$\mathcal{M}_{{\rm BCRW}}$. This implies that 
 the irreducible vectors for
$ \mathcal{M_{-V}} $ also are close to one of $9$ quarter planes $B_ j$.
A relatively easy computation shows  that the intersection of any two quarter
planes
$A_i \cap B_j$   consists only of the origin, which implies that
the set of irreducible vectors is bounded.

\medskip

In order to significantly decrease the size of the set of irreducible vectors
we will  enlarge $\mathcal{M}_{\Sz}$. This can be done using the
following refinement of Section~\ref{Some combinatorics behind the BCRW
Trick}:
\begin{lemma}
\label{enlarge system}
Let $\mathcal{M}$ be as in \emph{Definition~\ref{irreducible}}.
Let $M_1, M_2 \in \mathcal{M}$ and $v\in \Z^k$ be such that $M_1 \cap
(v + M_2)$ is  a single point.
Let  $M = (M_1 \setminus (v + M_2)) \cup ((v+M_2) \setminus M_1)$
denote the symmetric difference of $M_1$  and  $v + M_2$.
Then any set $S\subseteq  \Z^k $ that  is $\mathcal{M}$-closed is also
$\mathcal{M}
\cup \{M\}$-closed.
\end{lemma}
\Proof
Assume that $S$ is $\mathcal{M}$-closed. For any  $u$ such that
$|S \cap (u+M)| \geq |M|-1$ we have
$|S \cap (u+ M_1)| \geq |M_1|-1 $ or  $|S \cap (u+v+ M_2)| \geq |M_2|-1 $;
we assume the first  of these.
Since $S$ is $\mathcal{M}$-closed we have
$u + M_1 \subseteq S$; in particular the intersection
$(u + M_1) \cap (u+v + M_2)$ is in $S$.
Since  $|(u + M_1) \cap (u+v + M_2)|= |M_1 \cap (v + M_2)|=1$  by
hypothesis, it follows  that $|S \cap (u+v+ M_2)| \geq |M_2| -1$.
Since  $S$ is $\mathcal{M}$-closed, we see  that  $ u+v+ M_2
\subseteq S$ and hence that
$u+ M \subseteq S$, which shows that $S$ is also closed under $\{M\}$.
\qedd

\smallskip

This lemma is natural from a group-theoretic point of view.  For, in our
uses of Definition~\ref{irreducible}, for each  $M\in \mathcal{M}$ there
is a relation of the form $\prod_{m\in M}(u^m)^{\pm1}=1$.  In the
preceding lemma  we have  two such relations
$\prod_{m\in M_1}(u^m)^{\pm1} =1= \prod_{m\in M_2}(u^{v+m})^{\pm1}$.
By eliminating the one common term in these products  we obtain a third
relation  $\prod_{m\in M}(u^m)^{\pm1}=1$.
\smallskip

Let $\mathcal{M}_{22}$ be the system consisting of the $22$ sets in the
Table~\ref{bigset}.
Each of the  sets  $S_5$--$S_{22}$ in that table is the symmetric
difference of the two sets in the last two columns, and the preceding
lemma can be applied.
\begin{table}
\caption{The set $\mathcal{M}_{22}$}
\label{bigset}
$
\begin{array}{|l|l|l|l|}
\hline
S_{1} &
\scriptstyle \{(0,0,0,0),(1,0,0,0),(0,0,1,0)\}
&  &   \\
S_{2} &
\scriptstyle \{(0,0,0,0),(0,1,0,0),(0,0,0,1)\}
&  &   \\
S_{3} &
\scriptstyle \{(0,0,0,0),(1,-1,0,0),(0,0,1,-1)\}
&  &  \\
S_{4} &
\scriptstyle \{(0,0,0,0),(-2,1,0,0),(0,0,-2,1)\}
&  &  \\
S_{5} &
\scriptstyle \{(0,0,0,0),(0,0,1,0),(0,1,0,0),(0,1,1,-1)\}
& \scriptstyle S_{1} & \scriptstyle (0,1,0,0)+S_{3} \\
S_{6} &
\scriptstyle \{(0,0,0,0),(1,0,0,0),(0,0,0,1),(1,-1,0,1)\}
& \scriptstyle S_{1} & \scriptstyle (0,0,0,1)+S_{3} \\
S_{7} &
\scriptstyle \{(1,0,0,0),(0,0,1,0),(2,-1,0,0),(2,-1,-2,1)\}
& \scriptstyle S_{1} & \scriptstyle (2,-1,0,0)+S_{4} \\
S_{8} &
\scriptstyle \{(1,0,0,0),(0,0,1,0),(0,0,2,-1),(-2,1,2,-1)\}
& \scriptstyle S_{1} & \scriptstyle (0,0,2,-1)+S_{4} \\
S_{9} &
\scriptstyle \{(0,0,0,0),(0,0,1,0),(1,0,2,-1),(-1,1,2,-1)\}
& \scriptstyle S_{1} & \scriptstyle (1,0,2,-1)+S_{4} \\
S_{10} &
\scriptstyle \{(0,0,0,0),(1,0,0,0),(2,-1,1,0),(2,-1,-1,1)\}
& \scriptstyle S_{1} & \scriptstyle (2,-1,1,0)+S_{4} \\
S_{11} &
\scriptstyle \{(0,0,0,0),(0,0,0,1),(0,1,-1,1),(1,0,-1,1)\}
& \scriptstyle S_{2} & \scriptstyle (0,1,-1,1)+S_{3} \\
S_{12} &
\scriptstyle \{(0,0,0,0),(0,1,0,0),(-1,1,0,1),(-1,1,1,0)\}
& \scriptstyle S_{2} & \scriptstyle (-1,1,0,1)+S_{3} \\
S_{13} &
\scriptstyle \{(0,0,0,0),(0,0,1,-1),(1,-1,2,-1),(-1,0,2,-1)\}
& \scriptstyle S_{3} & \scriptstyle (1,-1,2,-1)+S_{4} \\
S_{14} &
\scriptstyle \{(0,0,0,0),(1,-1,0,0),(2,-1,1,-1),(2,-1,-1,0)\}
& \scriptstyle S_{3} & \scriptstyle (2,-1,1,-1)+S_{4} \\
S_{15} &
\scriptstyle \{(0,0,0,0),(0,0,0,1),(1,0,-2,1),(1,0,-1,1),(2,0,0,0)\}
& \scriptstyle S_{2} & \scriptstyle (1,0,-2,1)+S_{9} \\
S_{16} &
\scriptstyle \{(0,0,0,0),(0,1,0,0),(-2,1,1,0),(-1,1,1,0),(0,0,2,0)\}
& \scriptstyle S_{2} & \scriptstyle (-2,1,1,0)+S_{10} \\
S_{17} &
\scriptstyle \{(0,0,0,0),(0,1,0,0),(1,-1,0,0),(1,0,0,0),(0,0,1,0)\}
& \scriptstyle S_{2} & \scriptstyle (1,-1,0,0)+S_{12} \\
S_{18} &
\scriptstyle
\Big\{  \hspace{-4pt}
\begin{array}{l}
\scriptstyle (1,0,0,0),(0,0,1,0),(0,0,2,-1), \\
\scriptstyle (-1,0,2,-1),(0,0,3,-2),(0,0,1,-1)
\end{array} \hspace{-4pt}
\Big\} \hspace{-4pt}
& \scriptstyle S_{8} & \scriptstyle (-2,1,2,-1)+S_{14}\\
S_{19} &
\scriptstyle
\Big\{ \hspace{-4pt}
\begin{array}{l}
\scriptstyle (1,0,0,0),(0,0,1,0),(2,-1,0,0), \\
\scriptstyle (2,-1,-1,0),(3,-2,0,0),(1,-1,0,0)
\end{array} \hspace{-4pt}
\Big\} \hspace{-4pt}
& \scriptstyle S_{7} & \scriptstyle (2,-1,-2,1)+S_{13}\\
S_{20} &
\scriptstyle
\Big\{ \hspace{-4pt}
\begin{array}{l}
\scriptstyle (0,0,0,0),(0,0,1,0),(1,-1,0,0), \\
\scriptstyle (-1,0,1,0),(0,0,1,0),(1,-1,2,0)
\end{array} \hspace{-4pt}
\Big\}
& \scriptstyle S_{1} & \scriptstyle (1,-1,0,0)+S_{16}\\
S_{21} &
\scriptstyle
\Big\{  \hspace{-4pt}
\begin{array}{l}
\scriptstyle (0,0,1,0),(1,0,2,-1),(-1,1,2,-1), \\
\scriptstyle (0,0,1,-1),(1,-1,2,-1),(-1,0,2,-1)
\end{array} \hspace{-4pt}
\Big\}
& \scriptstyle S_{9} & \scriptstyle  S_{13} \\
S_{22} &
\scriptstyle
\Big\{ \hspace{-4pt}
\begin{array}{l}
\scriptstyle (1,0,0,0),(2,-1,1,0),(2,-1,-1,1), \\
\scriptstyle (1,-1,0,0),(2,-1,1,-1),(2,-1,-1,0)
\end{array} \hspace{-4pt}
\Big\}
& \scriptstyle S_{10} & \scriptstyle  S_{14} \\

\hline
\end{array}
$
\end{table}

By  Lemma~\ref{enlarge system},  any set which is closed
under $\mathcal{M}_{\Sz}$ is also closed under $\mathcal{M} _{22}$.
Consequently, the set of irreducible vectors for $\mathcal{M} _{22}$ is a
subset of the  set of  irreducible vectors for $\mathcal{M}_{\Sz}$.
A computer calculation 
  gives that this system has only the $33$ 
irreducible vectors listed in Table~\ref{31-irr} (this involves
checking    that  106 vectors in Table~\ref{139-irr} reduce
under the new tuples in  $\M_{22}$).%
\begin{table}
\caption{The irreducible vectors for $\mathcal{M}_{22}$}
\label{31-irr}
$
\begin{array}{|lllllll|}
\hline
\scriptstyle (-2,-1,0,0) &
\scriptstyle (-1,-1,0,0) &
\scriptstyle (-1,0,0,0)  &
\scriptstyle (0,-1,0,0)  &
\scriptstyle (0,-1,1,0)  &
\scriptstyle (0,0,-2,-1) &
\scriptstyle (0,0,-1,-1) \\
\scriptstyle (0,0,-1,0)  &
\scriptstyle (0,0,-1,1)  &
\scriptstyle (0,0,0,-1)  &
\scriptstyle (0,0,0,0)   &
\scriptstyle (0,0,0,1)   &
\scriptstyle (0,0,1,0)   &
\scriptstyle (0,0,1,1)   \\
\scriptstyle (0,1,0,0)   &
\scriptstyle (0,1,1,0)   &
\scriptstyle (1,0,0,-1)  &
\scriptstyle (1,0,0,0)   &
\scriptstyle (1,0,0,1)   &
\scriptstyle (1,0,1,0)   &
\scriptstyle (1,0,1,1)   \\
\scriptstyle (1,1,0,0)   &
\scriptstyle (1,1,1,0)   &
\scriptstyle (1,1,1,1)   &
\scriptstyle (1,1,2,0)   &
\scriptstyle (1,1,2,1)   &
\scriptstyle (2,0,1,1)   &
\scriptstyle (2,1,1,1)   \\
\scriptstyle (2,1,2,1)   &
\scriptstyle (2,2,3,1)   &
\scriptstyle (3,2,3,2)   &  
\scriptstyle (-2,-2,0,0)  & \scriptstyle (0,0,-2,-2) & &\\
\hline
\end{array}
$
\end{table}

By slightly perturbing the length function we will further
shrink the  set of irreducible vectors
without introducing any new irreducible vectors.
With respect to the  length
given by
$|(x_1,x_2,x_3,x_4)|^2=1.01 x_1^2+2x_2^2+1.011x_3^2+2.002x_4^4$ 
there are only $16$  
irreducible vectors  for  $\M_{22}$, listed in
Table~\ref{14-irr}. 
(Even without leaving the original collection $\M_{\Sz}$, 
a simpler calculation shows that this new length
decreases the number of irreducible vectors from 139 to 69.)
Therefore, using this length we have:
\begin{table}
\caption{The set $S_0$ of irreducible vectors for $\mathcal{M}_{22}$,
with respect to the perturbed length}
\label{14-irr}
$
\begin{array}{|lllll|}
\hline
(-1,-1,0,0) &
(-1,0,0,0)  &
(0,-1,0,0)  &
(0,0,-1,-1) &
(0,0,-1,0)  \\
(0,0,0,-1)  &
(0,0,0,0)   &
(0,1,0,0)   &
(1,0,0,0)   &
(1,0,1,0)   \\
(1,1,0,0)   &
(1,1,1,0)   &
(1,1,1,1)   &
(2,1,1,1)   & \\
(-2,-2,0,0) & (-2,-1,0,0)& &&\\
\hline
\end{array}
$
\end{table}

\begin{proposition}
\label{the 14}
The  only  $\mathcal{M}_{\Sz}$-closed subset of $\Z^4$ which contains the
 set $S_0$ of $16$ vectors in \emph{Table~\ref{14-irr}}   is
$\Z^4$.
\end{proposition}

\section{Ree triples and quadruples}
\label{Ree triples} 

We will sketch the ideas involved in the proof of
Proposition~\ref{Ree Borel}.  We use the notation  for $U$ and $h_\z$
introduced before that proposition.
Let $V$ denote the subgroup $x(0,*,*)$ of $U$, and let $W$ be the subgroup
$x(0,0,*)$.
As in 
\eqn{$Z^4$ action}, if $(i,j,k,l,m,n)\in \Z^6$  and  $x\in  U$ we write
\begin{equation}
\label{$Z^6$ action}
\mbox{
$x^{(i,j,k,l,m,n)}:= x^{  h}$ with $h=
h_{\z} ^i
h_{\z^\theta}^j
h_{\z+1} ^k
h_{\z^\theta+1} ^l
h_{\z+2} ^m
h_{\z^\theta+2} ^n
$.}
\end{equation}

The following trivial  identities in  $F$  
$$
\begin{array}{lllll}
(\z) + 1   =    (\z +1)
&&
(\z)   =  1+ (\z +2)
\\
(\z+1) + 1     =  (\z +2)
&&
(\z^\theta) + 1    =   (\z^\theta +1)
\\
(\z^\theta)  = 1+  (\z^\theta +2)
&&
(\z^\theta+1) + 1    =   (\z^\theta +2)
\end{array}
$$
 become  relations in $U/V$ by \eqn{Suzuki conjugation}: 
$$
u^{(0,0,0,0,0,0)} u^{(1,0,0,0,0,0)} \equiv  u^{(0,0,1,0,0,0)} 
\qquad
u^{(1,0,0,0,0,0)} \equiv  u^{(0,0,0,0,0,0)} u^{(0,0,0,0,1,0)} 
$$ 
$$
u^{(0,0,0,0,0,0)} u^{(0,0,1,0,0,0)} \equiv  u^{(0,0,0,0,1,0)} 
\qquad
u^{(0,0,0,0,0,0)} u^{(0,1,0,0,0,0)}  \equiv  u^{(0,0,0,1,0,0)} 
$$
$$
u^{(0,1,0,0,0,0)} \equiv  u^{(0,0,0,0,0,0)} u^{(0,0,0,0,0,1)}  
\qquad
u^{(0,0,0,0,0,0)} u^{(0,0,0,1,0,0)} \equiv  u^{(0,0,0,0,0,1)}  
$$ 
for any  $u\in U$. These relations correspond to the following collection
  $\mathcal{M_U}  $ of
triples  from
$\Z^6$:
\vspace{2pt}
 {\footnotesize
$$
\begin{array}{llll}
\{(0,0,0,0,0,0),(1,0,0,0,0,0),(0,0,1,0,0,0)\} &
\{(0,0,0,0,0,0),(1,0,0,0,0,0),(0,0,0,0,1,0)\}
\\
\{(0,0,0,0,0,0),(0,0,1,0,0,0),(0,0,0,0,1,0)\}  &
\{(0,0,0,0,0,0),(0,1,0,0,0,0),(0,0,0,1,0,0)\}
\\
\{(0,0,0,0,0,0),(0,1,0,0,0,0),(0,0,0,0,0,1)\} &
\{(0,0,0,0,0,0),(0,0,0,1,0,0),(0,0,0,0,0,1)\}  .
\end{array}
$$
}

Since $\theta^2=3$, we have 
\begin{align*}
\left(\z^\theta\right )^{\theta-1} \z^{\theta-1} &=
\z^{(\theta+1)(\theta-1)} = \z^{\theta^2-1} = \z^2
\\
\left(\z^\theta\right )^{\theta-1}
\left(\z^{\theta-1} \right)^{3}
&=
\z^{(\theta-1)(\theta+3)} = \z^{\theta^2 + 2\theta -3} =
\left(\z^\theta\right)^2\!.
  \end{align*} 

The relations $1 + \z^2 + (\z+1)^2+ (\z+2)^2=0$ and
$ \z^2(\z+1)^2 + (\z+1)^2(\z+2)+
 (\z+2)^2\z^2=1$
imply 
$$
1 + (\z)^{\theta+1}(\z^\theta)^{\theta+1} +
(\z+1)^{\theta+1}(\z^\theta+1)^{\theta+1} +
(\z+2)^{\theta+1}(\z^\theta+2)^{\theta+1}
= 0,
$$
\begin{align*}
(\z)^{\theta+1}(\z^\theta)^{\theta+1} (\z+1)^{\theta+1}(\z^\theta+1)^{\theta+1} &+
(\z+1)^{\theta+1}(\z^\theta+1)^{\theta+1} (\z+2)^{\theta+1}(\z^\theta+2)^{\theta+1}
\\
& +~(\z+2)^{\theta+1}(\z^\theta+2)^{\theta+1} (\z)^{\theta+1}(\z^\theta)^{\theta+1}
= 1
\end{align*}
and two further identities obtained by applying $\theta$, 
which by \eqn{Ree conjugation} imply   relations in $V/W$: 
$$
\begin{array}{llll}
v^{(0,0,0,0,0,0,0)} v^{(1,1,0,0,0,0)}v^{(0,0,1,1,0,0)}v^{(0,0,0,0,1,1)}
\equiv 1
\\
v^{(0,0,0,0,0,0,0)} \equiv
v^{(1,1,1,1,0,0)}v^{(0,0,1,1,1,1)}v^{(1,1,0,0,1,1)}
\\
v^{(0,0,0,0,0,0,0)}
v^{(3,1,0,0,0,0)}v^{(0,0,3,1,0,0)}v^{(0,0,0,0,3,1)}\equiv 1
\\
v^{(0,0,0,0,0,0,0)} \equiv
v^{(3,1,3,1,0,0)}v^{(0,0,3,1,3,1)}v^{(3,1,0,0,3,1)}
\end{array}
$$
for any  $v\in V$. These relations correspond  to the following
  collection $\mathcal{M_V}$ of quadruples  from $\Z^6$:
 {\footnotesize
$$
\begin{array}{llll}
\{(0,0,0,0,0,0),(1,1,0,0,0,0),(0,0,1,1,0,0),(0,0,0,0,1,1)\}
\\
\{(0,0,0,0,0,0),(1,1,1,1,0,0),(0,0,1,1,1,1),(1,1,0,0,1,1)\}
\\
\{(0,0,0,0,0,0),(3,1,0,0,0,0),(0,0,3,1,0,0),(0,0,0,0,3,1)\}
\\
\{(0,0,0,0,0,0),(3,1,3,1,0,0),(0,0,3,1,3,1),(3,1,0,0,3,1)\}.
\end{array}
$$
}

Similarly we have 
$$
\left(\z^\theta\right )^{2-\theta} \left(\z^{2-\theta}
\right)^{2} =
\z^{(2-\theta)(\theta+2)} = \z^{4-\theta^2} = \z
$$ 
and
\begin{align*}
(\z^2)^{2-\theta} (\z^\theta)^{2-\theta} + 1 &  =     ((\z
+1)^2)^{2-\theta}(\z^\theta +1)^{2-\theta}
\\
(\z^2)^{2-\theta} (\z^\theta)^{2-\theta}  &  = 1+     ((\z
+2)^2)^{2-\theta}(\z^\theta +2)^{2-\theta}
\\
((\z +1)^2)^{2-\theta}(\z^\theta +1)^{2-\theta} + 1 &  =
((\z +1)^2)^{2-\theta}(\z^\theta +1)^{2-\theta},
\end{align*}
together with additional similar identities;
 by \eqn{Ree conjugation} these  
  imply  the following   relations in $W$: 
$$\vspace{1pt}
w^{(0,0,0,0,0,0,0)} w^{(2,1,0,0,0,0)} = w^{(0,0,2,1,0,0)}
\qquad
w^{(2,1,0,0,0,0)} = w^{(0,0,0,0,0,0,0)} w^{(0,0,0,0,2,1)}
$$
$$
w^{(0,0,0,0,0,0,0)} w^{(0,0,2,1,0,0)} = w^{(0,0,0,0,2,1)}
\qquad
w^{(0,0,0,0,0,0,0)} w^{(3,2,0,0,0,0)} = w^{(0,0,3,2,0,0)}
$$
$$
w^{(3,2,0,0,0,0)} = w^{(0,0,0,0,0,0,0)} w^{(0,0,0,0,3,2)}
\qquad
w^{(0,0,0,0,0,0,0)} w^{(0,0,3,2,0,0)} = w^{(0,0,0,0,3,2)}
$$ 

\noindent
for any  $w\in W$. These relations correspond  to the following
 collection   $\mathcal{M_W} $   of triples  from  $\Z^6$:
 {\footnotesize
$$
\begin{array}{llll}
\{(0,0,0,0,0,0),(2,1,0,0,0,0),(0,0,2,1,0,0)\}   &
\{(0,0,0,0,0,0),(2,1,0,0,0,0),(0,0,0,0,2,1)\}  
\\
\{(0,0,0,0,0,0),(0,0,2,1,0,0),(0,0,0,0,2,1)\}    &
\{(0,0,0,0,0,0),(3,2,0,0,0,0),(0,0,3,2,0,0)\} 
\\
\{(0,0,0,0,0,0),(3,2,0,0,0,0),(0,0,0,0,3,2)\}    &
\{(0,0,0,0,0,0),(0,0,3,2,0,0),(0,0,0,0,3,2)\} .
\end{array}
 $$
}

We will use the fact that all six systems
$$
\mathcal{M_U} \cup \mathcal{M_{-U}},~
\mathcal{M_U} \cup \mathcal{M_{-V}},~
\mathcal{M_U} \cup \mathcal{M_{-W}},~
\mathcal{M_V} \cup \mathcal{M_{-V}},~
\mathcal{M_V} \cup \mathcal{M_{-W}}, ~
\mathcal{M_W} \cup \mathcal{M_{-W}}
$$
are of finite type with respect to  suitable length functions  
in order
to show that $[U,U] \equiv 1$   (mod~$V$), $[U,V]\equiv1$   (mod $W$),
$[U,W]=1$,
$[V,V]\equiv1$   (mod $W)$,
$[V,W]=1$ and $[W,W]=1$, respectively. The first system contains 
 several 
copies of $\mathcal{M}_{\pm {\rm BCRW}}$, therefore it is of finite type. The
last system can be transformed to the first one by a linear transformation,
therefore 
 is also of finite type for the transformed length function. The next
lemma gives the same result for the remaining  four
systems;  we sketch a proof in
  Appendix~\ref{Ree triples}.
As in Section~\ref{Nonabelian $p$-groups: Suzuki groups.},
 we can avoid using the first and  last sets because we
can use Lemma~\ref{BCRW1 trick}  to show that  the group
$\langle u^{\langle h_\z\rangle }\rangle $ is abelian (mod $V$), and then impose additional relations
to make it invariant under all $h_\star$.

\begin{Lemma}
\label{ReeSystems}
Each of the systems
$\mathcal{M}_{{\rm Ree }_1}:=\mathcal{M_U} \cup \mathcal{M_{-V}},$
$\mathcal{M}_{{\rm Ree }_2}:=\mathcal{M_U} \cup \mathcal{M_{-W}},$
$\mathcal{M}_{{\rm Ree }_3}:=\mathcal{M_V} \cup \mathcal{M_{-V}}$ and
$\mathcal{M}_{{\rm Ree }_4}:=\mathcal{M_W} \cup \mathcal{M_{-V}}$
  in $\Z^6$ is of finite type  using suitable length functions for the
various systems. 
\end{Lemma}
 
\proof

As in the case of Lemma~\ref{A Sz},
 again the proof is long and tedious.
We use two $3$-dimensional analogues $\mathcal{M}_{3\BCRW}$
and $\mathcal{M}_{3\BCRW'}$ of BCRW:
$\mathcal{M}_{3\BCRW}$ is
$$
\begin{array}{llll}
\{(0,0,0),(1,0,0),(0,1,0)\}
\\
\{(0,0,0),(1,0,0),(0,0,1)\}
\\
\{(0,0,0),(0,1,0),(0,0,1)\}   ,
\end{array}
$$
consisting of   3   copies of $\mathcal{M}_{\BCRW}$;
and $\mathcal{M}_{3\BCRW'}$ is
$$
\begin{array}{llll}
 \{(0, 0, 0), (1, 0, 0), (0, 1, 0), (0, 0, 1)\}
\\
\{(0, 0, 0), (1,
1, 0), (0, 1, 1), (1, 0, 1)\} .
\end{array}
$$
We use the Euclidean  length defined by  $|(x,y,z)|^2 = {x^2+y^2+z^2}$.

\begin{lemma}
\begin{enumerate}
\item The irreducible vectors for $\mathcal{M}_{3\BCRW}$ are as follows
$($up to permutation of the coordinates$)$$:$
$$
(-a,0,0), (a,a,a), (a,a,a+1), (a,a,a-1)
$$
 for integers    $a \geq 0.$

\item  The irreducible vectors for $\mathcal{M}_{3\BCRW'}$ are as follows
$($up to permutation of the coordinates$)$$:$
$$
(-a,0,0), (-a,0,-1), (-a,0,1), (a,a,a), (a,a,a+1), (a,a,a-1),
(a,a+1,a-1)
$$
 for integers    $a \geq 0.$
\end{enumerate}
\end{lemma}

\Proof
(1) $\mathcal{M}_{3\BCRW}$ contains one copy of
$\mathcal{M}_{\BCRW}$ for each pair of coordinates,
therefore the projection of an irreducible vector to
any coordinate plane is an irreducible vector for
$\mathcal{M}_{\BCRW}$, i.e., is one of $(-a,0), (a,a), (a,a+1)$ for some
non-negative integer $a$.
It is easy to see that the
vectors with this property are exactly the ones in the statement of the lemma.

(2) Here we have to consider several different cases depending
on the signs of the coordinates.
We will do just one case and leave the remainder to the reader.

Let $(x,y,z)$ be an irreducible
vector with $x\geq y \geq z \geq 0$.
This vector reduces via $ (x-1,y,z),$ $(x,y,z),$ $(x-1,y+1,z),$ $(x-1,y+1,z)$
unless $x < y+2$. Similarly it reduces via $(x-1,y-1,z),$ $(x,y,z),$
$(x-1,y,z+1),$
$(x,y-1,z+1)$ unless $y < z+2$.
Therefore the only irreducible vectors of this form are the ones with
$$
x=y=z
\mbox{ or }
x=y=z+1
\mbox{ or }
x=y+1=z+1
\mbox{ or }
x=y+1=z+2.
$$
The other three  cases $x\geq y \geq 0 > z$;
$x\geq 0 > y \geq z$ and
$0> x\geq y \geq z$
are similar.~\qedd

\smallskip
\smallskip
\noindent\emph{Proof of} Lemma~\ref{ReeSystems}. 
As in the  preceding argument,  $\mathcal{M_U}$ can be
transformed
 by a  linear transformation
to a product of two copies of $\mathcal{M}_{3\BCRW}$, and $\mathcal{M_V}$
can be transformed to $\mathcal{M}_{3\BCRW'}$.
We   choose the length function so that our transformation 
turns it into   the Euclidean length function. 

In order to find all irreducible vectors for the system
$\mathcal{M}_X \cup \mathcal{M}_{-Y}$, where $X,Y\in
\{\mathcal{U}, \mathcal{V}, \mathcal{W}\}$, we can proceed as
follows: First we apply a linear transformation and pick a length function
such that $\mathcal{M}_X$ becomes a product of two copies of
$\mathcal{M}_{3\BCRW}$ or $\mathcal{M}_{3\BCRW'}$ with the Euclidean norm.
Any vector which is irreducible for $\mathcal{M}_X$ lies close
(a bounded distance from) one of $16$ quarter planes $A_i$,
so
we only have to check which of these vectors are irreducible
for $\mathcal{M}_{-Y}$. After doing the computations in each case
one   obtains a bound for the length of any irreducible vector,
which proves that the system is of finite type.
Unfortunately, the computations are very long and  tedious and  (as in
Appendix~\ref{Suzuki triples}) do not provide
any insight.

Equivalently one can argue that a vector which is
irreducible  for  $\mathcal{M}_{-Y}$ lies close to one of
$16$ quarter planes $B_j$,
and check that the intersection $A_i \cap B_j$ is trivial in all
cases.
\qedd

\smallskip

We are now ready to provide the presentation needed in
Proposition~\ref{Ree Borel}.
\smallskip

{\bf Generators:} $u, v, w, h_\star $ for $\star \in \{
  \z,  {\z+1},  {\z+2}, {\z^\theta},  {\z^\theta+1}, {\z^\theta+2}
\}$.

Define
$
x^{(i,j,k,l,m,n)}
$
as in \eqn{$Z^6$ action}.
\smallskip

{\bf Relations:}
\begin{itemize}
\item []
\begin{enumerate}
\item $[h_\star, h_{\bullet}] =1$ for  $\star \in \{
  \z,  {\z+1},  {\z+2}, {\z^\theta},  {\z^\theta+1}, {\z^\theta+2}
\}$.
\item $w^3=1$.
\item $ w^{(0,0,1,0,0,0)} =w w^{(1,0,0,0,0,0)} $.
\item $ w w^{(0,0,0,0,1,0)} =w^{(1,0,0,0,0,0)} $.
\item $[w,w^{(1,0,0,0,0,0)} ]=1$.
\item $[[w^{m _\z(x)}]]_{h_\z}=1$.
\item $w^{h_{\z^\theta}} =[[w^{x^{3^{k+1}}}]]_{h_\z}$,
where   $x^{3^{k+1}}$ is  reduced mod $ m_\z(x)$ in order to obtain  a
short relation.
\vspace{2pt}

\item $ w^{(0,0,0,1,0,0)} =w w^{(0,1,0,0,0,0)}   $\!.

\item $ w w^{(0,0,0,0,1,0)} = w^{(0,1,0,0,0,0)} $\!.

\item $v^3 =w_1$.

\item $v v^{(1,1,0,0,0,0)} v^{(0,0,1,1,0,0)}
v^{(0,0,0,0,1,1)}=w_2$.

\item $v v^{(3,1,0,0,0,0)} v^{(0,0,3,1,0,0)}
v^{(0,0,0,0,3,1)}=w_3$.
\item
$  v^{(1,1,1,1,0,0)} v^{(0,0,1,1,1,1)}
v^{(1,1,0,0,1,1)}
=vw_4$.

\item
$  v^{(3,1,3,1,0,0)} v^{(0,0,3,1,3,1)}
v^{(3,1,0,0,3,1)}
=vw_5$.

\item
  $[w^{S_1},v]=1$, where $S_1$ is the set of  irreducible vectors for
$   \mathcal{M_W} \cup \mathcal{M_{-V}}   $.

\item
$[v^{S_2},v]=w_6$, where $S_2$ is the set of  irreducible vectors for
$  \mathcal{M_V} \cup \mathcal{M_{-V}}$.

\item$
[[v^{m_{\z^2}(x)}]]_{h_\z h_{\z^\theta}} =w_7 $.

\vspace{2pt}
\item
  $v^{h_\star} = [[v^{f_{\star^{\theta -2} ;\z^2} (x)  }]]_
{   {h_\z}  {h_{\z^\theta} }}
\, w_\star$  for $\star \in \{
  \z,  {\z+1},  {\z+2}, {\z^\theta},  {\z^\theta+1}, {\z^\theta+2}
\}$.

\item
$u^3 =v_1$.

\item
$u^{(0,0,2,1,0,0)}=u u^{(2,1,0,0,0,0)} v_2$.

\item
$u u^{(0,0,0,0,2,1)} = u^{(2,1,0,0,0,0)} v_3$.

\item
$u^{(0,0,3,2,0,0)}  = u u^{(3,2,0,0,0,0)} v_4$.

\item
$u u^{(0,0,0,0,3,2)} = u^{(3,2,0,0,0,0)} v_5$.

\item
$[u,u^{(2,1,0,0,0,0)}]=v_6$.

\item
$[u^{S_3},w] =1$, where $S_3$ is the set of irreducible vectors for
$\mathcal{M_U} \cup \mathcal{M_{-W}}$.

\item
$[u^{S_4},v] =w_9$, where $S_4$ is the set of irreducible vectors for
$   \mathcal{M_U} \cup \mathcal{M_{-V}}$.

\item
$[[u^{m_{\z}(x)}]]_{h_\z^2 h_{\z^\theta}} = v_8 $.

\vspace{2pt}
\item
$u^{h_\star} = [[u ^{f_{\star^{\theta -2} ;\z}(x)}]] _{   h_{{\z}^2}
h_{{\z^\theta}} } v_\star$
for $ \star \in \{
  \z,  {\z+1},  {\z+2}, {\z^\theta},  {\z^\theta+1}, {\z^\theta+2}
\}$.

\end{enumerate}
\end{itemize}
Here,
$w_1,\dots $ and
$v_1,\dots $ are suitable elements of   $W:=\langle  w^{\langle  \text{all
$h_\star$}\rangle }
\rangle $  and $V:=\langle  v^{\langle  \text{all
$h_\star$}  \rangle } ,W \rangle $, respectively.  These
depend on our initial choices
 $\z, u,v,w$ as well as on the length functions
 used to determine
the sets $S_i$.

\smallskip

The proof that the above is a presentation for an infinite central
extension of $B$  is similar to the Suzuki case.
We  omit the details.

There are $46+ \sum_i|S_i|$ relations in this presentation, which would
probably be somewhat unmanageable     in practice. \qedd

\end{document}